\documentclass[hidelinks,10pt]{elsarticle}
\usepackage[paperwidth=7.0in, paperheight=10.0in, margin=.875in]{geometry}
\usepackage[utf8]{inputenc}
\usepackage{graphicx}
\usepackage{float}
\usepackage{color}
\usepackage{microtype}
\usepackage{hyperref}

\usepackage{subfigure}
\usepackage{subfigmat}
 \usepackage{amssymb}
\usepackage{amsmath}
 \usepackage{amsthm}

\begin{document}

\begin{frontmatter}



\title{Use of Jordan forms for convection-pressure split Euler solvers}

 \author[1]{Naveen Kumar Garg\fnref{a}}
 \fntext[a]{Current address: Post Doctoral Fellow, TIFR Center for Applicable Mathematics, Bangalore, India}
 \ead{garg.naveen70@gmail.com, naveen@tifrbng.res.in}
 \author[2]{S.V. Raghurama Rao}
 \ead{raghu@aero.iisc.ernet.in}
 \author[3]{M. Sekhar}
 \ead{muddu@civil.iisc.ernet.in}
 \address[1]{Department of Mathematical Sciences, Indian Institute of Science (IISc), Bangalore, India}
 \address[2]{Department of Aerospace Engineering, IISc, Bangalore, India}
 \address[3]{Department of Civil Engineering, IISc, Bangalore, India}
\begin{abstract}
 In this study, we analyze convection-pressure split Euler flux functions which contain genuine weakly hyperbolic convection subsystems.  A system is said to be a genuine weakly hyperbolic if all eigenvalues are real with no complete set of linearly independent (LI) eigenvectors.  To construct an upwind solver based on flux difference splitting (FDS) framework, we require to generate complete set of LI eigenvectors. This can be done through addition of generalized eigenvectors which can be computed from theory of Jordan canonical forms. Once we have complete set of LI generalized eigenvectors, we construct upwind solvers in convection-pressure splitting framework. Since generalized eigenvectors are not unique, we take extra care to ensure no direct contribution of generalized eigenvectors in the final formulation of both the newly developed numerical schemes. First scheme is based on Zha and Bilgen type splitting approach, while second is based on Toro \& V\'azquez  splitting.  Both the schemes are tested on several bench-mark test problems on 1-D and one of them is tested on some typical 2-D test problems which involve shock instabilities.  The concept of generalized eigenvector based on Jordan forms is found to be useful in dealing with the genuine weakly hyperbolic parts of the considered Euler systems.  
\end{abstract}

\end{frontmatter}


%


\section{Introduction}
\label{}
Numerical algorithms based on compressible Euler systems are of great importance and are frequently used in simulations. These algorithms can be broadly divided into two major categories, namely, central discretization methods and upwind discretization methods. In this study we mainly focus on upwind methods and, in particular, on the Flux Difference Splitting (FDS) based upwind schemes.  Another class of upwind methods based on Flux Vector Splitting (FVS) schemes, like Steger \& Warming \cite{Steger_&_Warming} and van Leer \cite{Leer} schemes, are constructed using eigenvector structure of Euler system.  But these schemes are quite diffusive and can't capture isolated contact discontinuities crisply.  Similarly, schemes based on FDS framework, such as an approximate Riemann solvers of Roe \cite{Roe} and Osher \cite{Osher}, also depend on the eigen-structure (both eigenvalues and eigenvectors) but are quite accurate.  Osher scheme captures expansion waves well but is computationally expensive.  Roe's approximate Riemann solver is accurate and can capture steady discontinuities either exactly or with a single  interior point, without being as expensive.  Because of less numerical diffusion, Roe scheme tends to produce unphysical expansion shocks \cite{wesseling}, post-shock oscillations \cite{Stiriba_Donat_postshock_oscillations} and it is non-trivial to avoid instability problems \cite{Quirk} in its application. Similarly, HLL scheme \cite{HLL_SIAM_1983} and HLLC scheme \cite{Toro_Spruce} are special upwind schemes mainly dependent on structure of eigenvalues. HLLC scheme is a modified version of HLL scheme and can capture an isolated contact discontinuity exactly but suffers from infamous carbuncle phenomena and some other shock instabilities in 2-D, as shown in \cite{Huang_Wu_Yan}, \cite{Shen_Yan}.

There is a different class of upwind schemes based on splitting of the Euler flux function, with several possible splittings.  We consider here the popular splittings of the Euler flux function as is usually done in three distinct ways. In the first category, we have the AUSM family of schemes \cite{Meng_Sing_Liou_Review}, \cite{Liou_&_Steffen_1991} and CUSP scheme of Jameson \cite{Jameson_1}, in which flux vector is split into a convection and a pressure part.  Liou and Steffen first introduced the category of AUSM family of schemes, such that in their convection-pressure splitting, the pressure term is alone present in the pressure part of the split flux vector.  The Second category of convection-pressure split schemes were initially proposed by Steger and Warming in \cite{Steger_&_Warming}, but were throughly explored by Jameson \cite{Jameson_2}, Zha \& Bilgen \cite{Zha_Bilgen}, Balakrishnan \& Deshpande \cite{Balakrishnan_and_Deshpande} and by Raghurama Rao \& Deshpande \cite{Rao_PVU}.  The structure of this splitting is such that pressure term of momentum equation and the term containing the product of pressure and velocity in the energy equation constitute the pressure part of the convection-pressure split Euler flux function.  The main feature of this type of convection-pressure splitting is that the eigenvalues corresponding to Jacobian of pressure subsystem become completely free from fluid velocity $u$, unlike in the splitting utilized by Meng-Sing Liou and others.  Another interesting and more recent convection-pressure splitting is proposed by Toro \& V\'azquez \cite{Toro_Vazquez}.  In this category, the authors split convection and pressure parts in such a way that convection flux becomes completely free from pressure terms.  For all the three splittings, convection part always corresponds to a genuine weakly hyperbolic system.  As the each convection part is weakly hyperbolic, we can utilize the theory of Jordan Canonical forms to recover complete set of linearly independent (LI) generalized eigenvectors.  In contrast, each pressure part corresponds to strict or non-strict hyperbolic system.  

In the first category of convection-pressure splitting, although the convection part contains contribution of two different eigenvalues namely, $u$ and $\gamma u$, with $u$ as repeated eigenvalue, the eigenvalues of pressure part do not contain any contribution of acoustic speed and this may result in an unstable scheme if used in FDS framework.  As pressure part of other two splittings contain the contribution of acoustic speeds as well, we propose two numerical schemes namely, Zha and Bilgen split - flux difference splitting scheme (ZBS-FDS)  and Toro and V\'azquez split - flux difference splitting TVS-FDS scheme.  The main idea of each scheme is first to construct traditional FDS scheme for pressure sub-system and then utilize the resulting averaged values of all variables together with theory of Jordan forms for convection part.  Our motivation is to develop efficient and workable flux difference split schemes based on convection-pressure splitting, together with the use of Jordan forms for convection subsystems, rather than focusing on reproducing an ideal  approximate Riemann solver in this framework.  Both schemes are tested on various shock tube problems in 1-D and are found to require no entropy fix for sonic points and in strong expansion regions. Both the schemes capture isolated and steady contact discontinuities exactly.  Out of the two, ZBS-FDS scheme is extended to two dimensions and is further tested on a variety of shock instability problems, including shock diffraction around a corner, flow over over a half cylinder and reflection of a plane shock from a wedge.     

\section{Convection-Pressure splittings for Euler flux function}  
Consider the one-dimensional inviscid Euler system
\begin{equation}
\frac{\partial\boldsymbol{U}}{\partial t}  \ + \ \frac{\partial \boldsymbol{F} \left( \boldsymbol{U} \right)} {\partial x} \ = \ \boldsymbol{0}, \ \ \ \ (x, t) \in {\rm I\!R}\times [0,\infty),
\end{equation}
where $\boldsymbol{U}:({\rm I\!R}\times{\rm I\!R^{+}})\longmapsto \Omega\subseteq{\rm I\!R^{3}}$, $\Omega$ is an open subset, is the conserved variable vector and $\boldsymbol{F}:\Omega\longmapsto {\rm I\!R^{3}}$ is the flux vector defined by
\begin{equation} 
 \boldsymbol{U} = \begin{bmatrix} 
     \rho \\[0.3em]  
		 \rho u \\[0.5em] 
		 \rho E 
		\end{bmatrix} \ \mbox{and} \ 
 \boldsymbol{\boldsymbol{F} \left( \boldsymbol{U} \right)} = \begin{bmatrix} 
        \rho u \\[0.3em] 
				p + \rho u^{2} \\[0.3em] 
				p u + \rho u E 
				\end{bmatrix} 
\end{equation} 
Here the total energy $E$ is defined as the sum of internal energy ($e$) and kinetic energy ($\frac{1}{2}u^{2}$), given as: $E = e + \frac{1}{2} u^{2} = \frac{p}{\rho \left(\gamma - 1\right)} + \frac{1}{2} u^{2}$.  Till now, three distinct convection-pressure splittings have been proposed, which are described in the following.  
\subsection{Liou and Steffen splitting procedure}
 Liou and Steffen, in formulating their upwind scheme \cite{Liou_&_Steffen_1991}, introduced a unique convection-pressure splitting by taking out pressure part from momentum equation of full Euler system.  
 \begin{equation}
 \boldsymbol{F} \ =\ \boldsymbol{F}_{c}^{\boldsymbol{LS}}  + \boldsymbol{F}_{p}^{\boldsymbol{LS}}
\end{equation}
where
\begin{equation}
 \boldsymbol{F}_{c}^{\boldsymbol{LS}}  = \begin{bmatrix}
         \rho u \\[0.3em]
         \rho u^{2} \\[0.3em]
         \rho u E + p u
         \end{bmatrix} \ \mbox{and} \
 \boldsymbol{F}_{p}^{\boldsymbol{LS}}  = \begin{bmatrix}
         0  \\[0.3em]
         p  \\[0.3em]
         0 
        \end{bmatrix}
\end{equation}
 Let us split the system $(1)$ into Liou and Steffen type convection and pressure subsystems, for gaining better insight by analyzing each part separately.  
\begin{equation}
\frac{\partial \boldsymbol{U}}{\partial t} \ + \ \frac{\partial \boldsymbol{F}_{c}^{\boldsymbol{LS}} \left(\boldsymbol{U} \right)} {\partial x} \ =  \ \boldsymbol{0}
\end{equation}
\mbox{and} \
\begin{equation}
 \frac{\partial \boldsymbol{U}}{\partial t}  \ + \ \frac{\partial \boldsymbol{F}_{p}^{\boldsymbol{LS}} \left( \boldsymbol{U} \right)} {\partial x} \ =  \ \boldsymbol{0}
\end{equation}
Both subsystems can also be written in quasilinear form as follows.  
\begin{equation}
\frac{\partial \boldsymbol{U}}{\partial t} \ + \ \boldsymbol{A}_{c}^{\boldsymbol{LS}} \frac{\partial \boldsymbol{U}} {\partial x} \ =  \ \boldsymbol{0}
\end{equation}
\begin{equation}
\frac{\partial\boldsymbol{U}}{\partial t}  \ + \  \boldsymbol{A}_{p}^{\boldsymbol{LS}} \frac{\partial \boldsymbol{U}} {\partial x} \ =  \ \boldsymbol{0}
\end{equation}
Here $\boldsymbol{A}_{c}^{\boldsymbol{LS}}$ and $\boldsymbol{A}_{p}^{\boldsymbol{LS}}$ are Jacobian matrices for convection and pressure parts respectively and are given by
\begin{equation*}
 \boldsymbol{A}_{c}^{\boldsymbol{LS}} = \begin{bmatrix}
        \ 0  &&  1  &&   0   \\[0.3em]
       \ -u^{2}  &&  2 u  &&  0  \\[0.3em]
       \ -\gamma u E + (\gamma - 1) u^{3}  &&    \gamma E -\frac{3}{2}(\gamma - 1) u^{2} && \gamma u
        \end{bmatrix} \
        \end{equation*}
        
        \mbox{and} \
 \begin{equation*}
   \boldsymbol{A}_{p}^{\boldsymbol{LS}} = \begin{bmatrix}
         \ 0  && 0 && 0  \\[0.3em]
         \frac{1}{2} (\gamma -1) {u^2}  &&  -(\gamma -1) u  &&  (\gamma -1) \\[0.3em]
         \ 0  &&  0  &&  0 
        \end{bmatrix}
\end{equation*} 
Eigenvalues corresponding to convective Jacobian matrix $\boldsymbol{A}_{c}^{\boldsymbol{LS}}$ are $ \lambda_{c,1}^{\boldsymbol{LS}}   =  \gamma u, \ \ \lambda_{c,2}^{\boldsymbol{LS}} = \lambda_{c,3}^{\boldsymbol{LS}}  =  u$ and algebraic multiplicity (AM) of the eigenvalue $u$ is 2. Similarly, eigenvalues corresponding to pressure Jacobian matrix $\boldsymbol{A}_{p}^{\boldsymbol{LS}}$ are $ \lambda_{p,1}^{\boldsymbol{LS}}   =  -(\gamma - 1) u, \ \ \lambda_{p,2}^{\boldsymbol{LS}} = \lambda_{p,3}^{\boldsymbol{LS}}  =  0$.  Since AM of $u$ is 2, so we have to find its eigenvector space to see whether $\boldsymbol{A}_{c}^{\boldsymbol{LS}}$ has complete set of linearly independent eigenvectors or not. The analysis of matrix $\boldsymbol{A}_{c}^{\boldsymbol{LS}}$ shows that convective subsystem is weakly hyperbolic as there is no complete set of linearly independent eigenvectors. Indeed, its eigenvectors are 
      \begin{equation}
      \boldsymbol{R}_{c,1}^{\boldsymbol{LS}}  =  \begin{bmatrix}
                  \ 0    \\[0.3em]
                  \ 0     \\[0.3em]
                  \ 1 
                 \end{bmatrix}  ~~\mbox{and}~~ \
      \boldsymbol{R}_{c,2}^{\boldsymbol{LS}} =  \begin{bmatrix}
                  \ 1    \\[0.3em]
                  \ u     \\[0.3em]
                  \ \frac{1}{2}u^2 
                 \end{bmatrix} \          
  \end{equation}
Similarly, eigenvectors corresponding to $ \lambda_{p,1}^{\boldsymbol{LS}}   =  -(\gamma - 1) u, \ \ \lambda_{p,2}^{\boldsymbol{LS}} = \lambda_{p,3}^{\boldsymbol{LS}}  =  0$ of pressure subsystems are: 
\begin{equation}
      \boldsymbol{R}_{p,1}^{\boldsymbol{LS}}  =  \begin{bmatrix}
                  \ 0    \\[0.3em]
                  \ 1     \\[0.3em]
                  \ 0 
                 \end{bmatrix}   \ , \
       \boldsymbol{R}_{p,2}^{\boldsymbol{LS}} =  \begin{bmatrix}
                  \ 1    \\[0.3em]
                  \ 0     \\[0.3em]
                  \ -\frac{1}{2}u^2 
                 \end{bmatrix}   \ , \
       \boldsymbol{R}_{p,3}^{\boldsymbol{LS}} = \begin{bmatrix}
                  \ 0    \\[0.3em]
                  \ 1     \\[0.3em]
                  \ u  
                 \end{bmatrix}
   \end{equation}
Convection subsystem turns out to be weakly hyperbolic,  and Jordan theory can be applied to explore it further, whereas pressure subsystem is non-strict hyperbolic.  Apart from eigenvectors, traditional FDS solvers depend heavily on eigenvalues also, but for present case all eigenvalues are either $u$  or constant times $u$.  In other words, there is no direct or indirect contribution of acoustic speed $a$ as an eigenvalue for both subsystems. This is a serious issue as $u$ frequently goes to zero or near to zero in a flow field which results in zero or near zero diffusion at some parts of the flow. Thus, the scheme may generate near zero diffusion which effectively reduces the scheme to forward in time and central in space (FTCS) framework, and as FTCS doesn't preserve the hyperbolicity, the solution $blows$-$up$. In fact, we constructed FDS scheme for present splitting but unfortunately, it led to blow-up of the solution for almost all problems.  Note that we are only considering the application of flux difference splitting to Liou and Steffen splitting here and not their alternative  upwinding procedure.  
\subsection{Zha and Bilgen splitting procedure}
Another type of flux splitting is given by Zha and Bilgen \cite{Zha_Bilgen}, in which they split the full Euler flux function into convection and pressure fluxes in such a way that eigenvalues corresponding to Jacobian of pressure flux $\boldsymbol{A_{p}^{ZB}}$ contains no contribution of fluid velocity $u$, unlike in Liou and Steffen splitting.  Their convection-pressure splitting is as follows.  
\begin{equation}
 \boldsymbol{F} \ =\ \boldsymbol{F}_{c}^{\boldsymbol{ZB}}  + \boldsymbol{F}_{p}^{\boldsymbol{ZB}}
\end{equation}
where
\begin{equation}
 \boldsymbol{F}_{c}^{\boldsymbol{ZB}}  = \begin{bmatrix}
         \rho u \\[0.3em]
         \rho u^{2} \\[0.3em]
         \rho u E
         \end{bmatrix} \ \mbox{and} \
 \boldsymbol{F}_{p}^{\boldsymbol{ZB}}  = \begin{bmatrix}
         0  \\[0.3em]
         p  \\[0.3em]
         pu
        \end{bmatrix}
\end{equation}
As done earlier, we split system $(1)$ into convection and pressure subsystems, using Zha and Bilgen type flux splitting, separately as    
\begin{equation}\label{conservation_ZB_c_part}
\frac{\partial \boldsymbol{U}}{\partial t}  + \frac{\partial \boldsymbol{F}_{c}^{\boldsymbol{ZB}} \left(\boldsymbol{U} \right)} {\partial x} \ =  \ \boldsymbol{0}
\end{equation}
\mbox{and} \
\begin{equation}\label{conservation_ZB_p_part}
 \frac{\partial \boldsymbol{U}}{\partial t}  + \frac{\partial \boldsymbol{F}_{p}^{\boldsymbol{ZB}} \left( \boldsymbol{U} \right)} {\partial x} \ =  \ \boldsymbol{0}
\end{equation}
Again, both subsystems can also be written in quasilinear form as follows.  
\begin{equation}\label{quasi_form_ZB_c_part}
\frac{\partial \boldsymbol{U}}{\partial t}  + \boldsymbol{A}_{c}^{\boldsymbol{ZB}} \frac{\partial \boldsymbol{U}} {\partial x} \ =  \ \boldsymbol{0}
\end{equation}
\begin{equation}\label{quasi_form_ZB_p_part}
\frac{\partial\boldsymbol{U}}{\partial t}  + \boldsymbol{A}_{p}^{\boldsymbol{ZB}} \frac{\partial \boldsymbol{U}} {\partial x} \ =  \ \boldsymbol{0}
\end{equation}
where, $\boldsymbol{A}_{c}^{\boldsymbol{ZB}}$ and $\boldsymbol{A}_{p}^{\boldsymbol{ZB}}$ are Jacobian matrices for convection and pressure parts respectively and are given by
\begin{equation*}
 \boldsymbol{A}_{c}^{\boldsymbol{ZB}} = \begin{bmatrix}
        \ 0  &&  1  &&   0   \\[0.3em]
       \ -u^{2}  &&  2 u  &&  0  \\[0.3em]
       \ -u E  &&    E &&  u
        \end{bmatrix} \
        \end{equation*}
        \mbox{and} \
 \begin{equation*}
   \boldsymbol{A}_{p}^{\boldsymbol{ZB}} = \begin{bmatrix}
         \ 0  && 0 && 0  \\[0.3em]
         \frac{1}{2} (\gamma -1) {u^2}  &&  -(\gamma -1) u  &&  (\gamma -1) \\[0.3em]
         \ - \frac{a^{2}u}{\gamma} + \frac{(\gamma -1)}{2} u^{3}  &&  \frac{a^{2}}{\gamma} - (\gamma - 1) u^{2} &&  (\gamma - 1) u 
        \end{bmatrix}
\end{equation*} 
Now, eigenvalues corresponding to convective Jacobian matrix $\boldsymbol{A}_{c}^{\boldsymbol{ZB}}$ are $ \lambda_{c,1}^{\boldsymbol{ZB}} = \lambda_{c,2}^{\boldsymbol{ZB}} = \lambda_{c,3}^{\boldsymbol{ZB}}  =  u$, thus algebraic multiplicity (AM) of eigenvalue $u$ is 3. Similarly, eigenvalues corresponding to pressure Jacobian matrix $\boldsymbol{A}_{p}^{\boldsymbol{ZB}}$ are $ \lambda_{p,1}^{\boldsymbol{ZB}}   =  -\sqrt{\frac{(\gamma - 1)}{\gamma}} a, \ \ \lambda_{p,2}^{\boldsymbol{ZB}} =  0 \ and \  \lambda_{p,3}^{\boldsymbol{ZB}}  = \sqrt{\frac{(\gamma - 1)}{\gamma}} a$. Since AM of $u$ is 3, so we have to find its eigenvector space to see whether $\boldsymbol{A}_{c}^{\boldsymbol{ZB}}$ has complete set of linearly independent eigenvectors or not. The analysis of matrix $\boldsymbol{A}_{c}^{\boldsymbol{ZB}}$ shows that convective subsystem is weakly hyperbolic as there is no complete set of linearly independent eigenvectors. Indeed, its eigenvectors are 
      \begin{equation}
      \boldsymbol{R}_{c,1}^{\boldsymbol{ZB}}  =  \begin{bmatrix}
                  \ 1    \\[0.3em]
                  \ u     \\[0.3em]
                  \ 0 
                 \end{bmatrix}  ~~\mbox{and}~~ \
      \boldsymbol{R}_{c,2}^{\boldsymbol{ZB}} =  \begin{bmatrix}
                  \ 0    \\[0.3em]
                  \ 0     \\[0.3em]
                  \ 1 
                 \end{bmatrix} \          
  \end{equation}
Since all  eigenvalues for pressure part are real and distinct, this makes pressure subsystem strictly hyperbolic.  Analysis of the flux Jacobian matrix for the pressure part shows complete set of eigenvectors, as given below.
\begin{equation}
      \boldsymbol{R}_{p,1}^{\boldsymbol{ZB}}  =  \begin{bmatrix}
                  \ 0    \\[0.3em]
                  \ 1     \\[0.3em]
                  \ u - \frac{a}{\sqrt{\gamma (\gamma - 1)}} 
                 \end{bmatrix}   \ , \
       \boldsymbol{R}_{p,2}^{\boldsymbol{ZB}} =  \begin{bmatrix}
                  \ 1    \\[0.3em]
                  \ u   \\[0.3em]
                  \ \frac{1}{2}u^2 
                 \end{bmatrix}   \ , \
       \boldsymbol{R}_{p,3}^{\boldsymbol{ZB}} = \begin{bmatrix}
                  \ 0    \\[0.3em]
                  \ 1     \\[0.3em]
                  \ u + \frac{a}{\sqrt{\gamma (\gamma - 1)}} 
                 \end{bmatrix}
   \end{equation}
\subsection{Toro and V\'azquez splitting Procedure}
More recently, Toro \& V\'azquez-Cend\'on \cite{Toro_Vazquez} presented a flux splitting in which convection part contains no pressure term  at all, leading to following type of splitting.  
\begin{equation}
 \boldsymbol{F} \ =\ \boldsymbol{F}_{c}^{\boldsymbol{TV}}  + \boldsymbol{F}_{p}^{\boldsymbol{TV}}
\end{equation}
where
\begin{equation}
 \boldsymbol{F}_{c}^{\boldsymbol{TV}}  = \begin{bmatrix}
         \rho u \\[0.3em]
         \rho u^{2} \\[0.3em]
        \frac{1}{2} \rho u^{3}
         \end{bmatrix} \ \mbox{and} \
 \boldsymbol{F}_{p}^{\boldsymbol{TV}}  = \begin{bmatrix}
         0  \\[0.3em]
         p  \\[0.3em]
         \frac{\gamma}{\gamma - 1} p u 
        \end{bmatrix}
\end{equation}
Let us split the system $(1)$ into convection and pressure subsystems, for gaining better insight by analysing each part separately.  
\begin{equation}\label{conservation_TV_c_part}
\frac{\partial \boldsymbol{U}}{\partial t}  + \frac{\partial \boldsymbol{F}_{c}^{\boldsymbol{TV}} \left(\boldsymbol{U} \right)} {\partial x} \ =  \ \boldsymbol{0}
\end{equation}
\mbox{and} \
\begin{equation}\label{conservation_TV_p_part}
 \frac{\partial \boldsymbol{U}}{\partial t}  + \frac{\partial \boldsymbol{F}_{p}^{\boldsymbol{TV}} \left( \boldsymbol{U} \right)} {\partial x} \ =  \ \boldsymbol{0}
\end{equation}
Again, both subsystems can also be written in quasilinear form as follows.  
\begin{equation}\label{quasi_form_TV_c_part}
\frac{\partial \boldsymbol{U}}{\partial t}  + \boldsymbol{A}_{c}^{\boldsymbol{TV}} \frac{\partial \boldsymbol{U}} {\partial x} \ =  \ \boldsymbol{0}
\end{equation}
\begin{equation}\label{quasi_form_TV_p_part}
\frac{\partial\boldsymbol{U}}{\partial t}  + \boldsymbol{A}_{p}^{\boldsymbol{TV}} \frac{\partial \boldsymbol{U}} {\partial x} \ =  \ \boldsymbol{0}
\end{equation}
Here $\boldsymbol{A}_{c}^{\boldsymbol{TV}}$ and $\boldsymbol{A}_{p}^{\boldsymbol{TV}}$ are Jacobian matrices for convection and pressure parts respectively and are given by
\begin{equation*}
 \boldsymbol{A}_{c}^{\boldsymbol{TV}} = \begin{bmatrix}
        \ 0  && 1  &&   0   \\[0.3em]
       \ -u^{2}  &&  2 u  &&  0  \\[0.3em]
       \ -u^{3} &&   \frac{3}{2}u^{2}  && 0
        \end{bmatrix} \
        \end{equation*}
        
        \mbox{and} \
 \begin{equation*}
   \boldsymbol{A}_{p}^{\boldsymbol{TV}} = \begin{bmatrix}
         \ 0  && 0 && 0  \\[0.3em]
         \frac{1}{2} (\gamma -1) {u^2}  &&  -(\gamma -1) u  &&  (\gamma -1) \\[0.3em]
         \ -\frac{u a^{2}}{(\gamma - 1)} + \frac{1}{2}\gamma u^{3}  &&  \frac{a^{2}}{(\gamma - 1)} - \gamma u^{2}  &&  \gamma u 
        \end{bmatrix}
\end{equation*} 
Eigenvalues corresponding to convective Jacobian matrix $\boldsymbol{A}_{c}^{\boldsymbol{TV}}$ are $ \lambda_{c,1}^{\boldsymbol{TV}}   =  0, \ \ \lambda_{c,2}^{\boldsymbol{TV}} = \lambda_{c,3}^{\boldsymbol{TV}}  =  u$ and algebraic multiplicity (AM) of eigenvalue $u$ is 2, so we have to find its eigenvector space to see whether $\boldsymbol{A}_{c}^{\boldsymbol{TV}}$ has complete set of linearly independent eigenvectors or not. The analysis of matrix $\boldsymbol{A}_{c}^{\boldsymbol{TV}}$ shows that convective subsystem is weakly hyperbolic as there is no complete set of linearly independent eigenvectors. Indeed, its eigenvectors are 
      \begin{equation}
       \boldsymbol{R}_{c,1}^{\boldsymbol{TV}}  =  \begin{bmatrix}
                  \ 0    \\[0.3em]
                  \ 0     \\[0.3em]
                  \ 1 
                 \end{bmatrix}  ~~\mbox{and}~~ \
       \boldsymbol{R}_{c,2}^{\boldsymbol{TV}} =  \begin{bmatrix}
                  \ 1    \\[0.3em]
                  \ u     \\[0.3em]
                  \ \frac{1}{2}u^2 
                 \end{bmatrix} \          
  \end{equation}
      
Similarly, the eigenvalues corresponding to pressure Jacobian matrix $\boldsymbol{A}_{p}^{\boldsymbol{TV}}$, when evaluated, are found to be   $\lambda_{p,1}^{\boldsymbol{TV}} = \frac{1}{2}(u - \beta),    \ \   \lambda_{p,2}^{\boldsymbol{TV}} = 0,  \ \  \lambda_{p,3}^{\boldsymbol{TV}}  = \frac{1}{2}(u + \beta) $, where $\beta$ \ = \ $\sqrt{u^2 + 4a^2}$. 
All eigenvalues for pressure part are real and distinct and this makes pressure subsystem strictly hyperbolic.  Analysis of the flux Jacobian matrix for the pressure part shows complete set of eigenvectors, as given below.  
      \begin{equation}
       \boldsymbol{R}_{p,1}^{\boldsymbol{TV}}  =  \begin{bmatrix}
                  \ 0    \\[0.3em]
                  \ 1     \\[0.3em]
                  \ u + \frac{1}{2}(\frac{u - \beta}{\gamma - 1}) 
                 \end{bmatrix}   \ , \
       \boldsymbol{R}_{p,2}^{\boldsymbol{TV}} =  \begin{bmatrix}
                  \ 1    \\[0.3em]
                  \ u     \\[0.3em]
                  \ \frac{1}{2}u^2 
                 \end{bmatrix}   \ , \
       \boldsymbol{R}_{p,3}^{\boldsymbol{TV}} = \begin{bmatrix}
                  \ 0    \\[0.3em]
                  \ 1     \\[0.3em]
                  \ u + \frac{1}{2}(\frac{u + \beta}{\gamma - 1}) 
                 \end{bmatrix}
   \end{equation}
   Since the convective subsystems for both Zha-Bilgen splitting and Toro-V\'azquez splitting have incomplete set of linearly independent (LI) eigenvectors, it will be nontrivial to construct an upwind scheme based on eigenvector structure.  
 \section{Addition of generalized eigenvectors}
First we consider Jacobian matrix $\boldsymbol{A}_{c}^{\boldsymbol{TV}}$ corresponding to Toro-V\'azquez convective subsystem, for which, our aim is to get complete set of linearly independent generalized eigenvectors. For this system, we have two different sets of eigenvalues. Here, we briefly discuss a procedure to find generalized eigenvectors for cases where resultant Jordan matrix possess exactly one Jordan block for each set of eigenvalues.
Let
\begin{equation*}
 {\boldsymbol J} = \begin{bmatrix}
        \ {\boldsymbol J}(\lambda_{1})  & \boldsymbol{0} & \cdots  &  \boldsymbol{0}     \\[0.3em]
       \  \boldsymbol{0} & {\boldsymbol J}(\lambda_{2}) & \cdots  &  \boldsymbol{0}   \\[0.3em]
       \ \vdots    &  \vdots   & \ddots   &  \vdots   \\[0.3em]
       \boldsymbol{0} & \boldsymbol{0} & \cdots & {\boldsymbol J}(\lambda_{p})
        \end{bmatrix}, \  \ \textrm{where} \  {\lambda_{1}, \lambda_{2}, \cdots , \lambda_{p}} \in \sigma(\boldsymbol{A})
\end{equation*} 
are set of distinct eigenvalues, some or all of them with arithmetic multiplicity greater than one. Moreover, assume there exists a single Jordan block for each $\lambda_{i}$. Let us focus on one such $ \lambda_{i}$, with AM equal to $m>1$. Then 
\begin{equation*}
 {\boldsymbol J(\lambda_{i})} = \begin{bmatrix}
        \   \lambda_{i}   &   1    &        \\[0.3em]
       \    &     \ddots   &   \ddots  &  \\[0.3em]
       \  &    &   \ddots   &  1   \\[0.3em]
       \  &       &        &  \lambda_{i}
        \end{bmatrix}_{m\times m} \
\end{equation*}  \\
In order to find set of generalized eigenvectors corresponding to $\lambda_{i}$, we need to focus on portion $\boldsymbol{P}^{*} = \big[\boldsymbol{X}_1, \boldsymbol{X}_2, \boldsymbol{X}_3,........., \boldsymbol{X}_m\big]$ of $\boldsymbol{P} \ = \ [... \boldsymbol{P}^{*}...]$ that corresponds to the position $\boldsymbol J(\lambda_{i})$ in $\boldsymbol J$. Now $\boldsymbol{A}\boldsymbol{P} = \boldsymbol{P}\boldsymbol{J}$ implies $\boldsymbol{A}\boldsymbol{P}^{*} = \boldsymbol{P}^{*}\boldsymbol J(\lambda_{i})$, {\em i.e.},
\begin{equation*}
 \boldsymbol{A}\big[\boldsymbol{X}_1, \boldsymbol{X}_2, \boldsymbol{X}_3,........., \boldsymbol{X}_{m}\big] \ = \ \big[\boldsymbol{X}_1, \boldsymbol{X}_2, \boldsymbol{X}_3,........., \boldsymbol{X}_{m}\big] \begin{bmatrix}
        \   \lambda_{i}   &   1    &        \\[0.3em]
       \    &     \ddots   &   \ddots  &  \\[0.3em]
       \  &    &   \ddots   &  1   \\[0.3em]
       \  &       &        &  \lambda_{i}
        \end{bmatrix}_{m\times m} \
\end{equation*}
 On equating columns on both sides, we get
 \begin{align}
 \begin{split}
  \boldsymbol{A}\boldsymbol{X}_{1}  \ &= \  \lambda_{i} \boldsymbol{X}_{1} \\
  \boldsymbol{A}\boldsymbol{X}_{2}  \ &= \  \lambda_{i} \boldsymbol{X}_{2} \ + \ \boldsymbol{X}_{1}  \\
  \boldsymbol{A}\boldsymbol{X}_{3}  \ &= \  \lambda_{i} \boldsymbol{X}_{3} \ + \ \boldsymbol{X}_{2}  \\
                                 \vdots     \\
\boldsymbol{A}\boldsymbol{X}_{m}  \ &= \  \lambda_{i} \boldsymbol{X}_{m} \ + \ \boldsymbol{X}_{m-1}
 \end{split}
\end{align}
\\
Now $u$ is a repeated eigenvalue of matrix $\boldsymbol{A}_c^{\boldsymbol{TV}}$ with AM is equal to two and other eigenvalue is zero with multiplicity one. First we have to compute ranks of matrices $(\boldsymbol{A}_{c}^{\boldsymbol{TV}}- u\boldsymbol{I})$, \  \ $(\boldsymbol{A}_{c}^{\boldsymbol{TV}}- u\boldsymbol{I})^2$, $\cdots$. It turns out that $rank(\boldsymbol{A}_{c}^{\boldsymbol{TV}}- u\boldsymbol{I})^2  \ =  \ 1  \ =  \  rank(\boldsymbol{A}_{c}^{\boldsymbol{TV}}- u\boldsymbol{I})^{3}$, which means there should be a  $Jordan$ block of order $2$. Therefore, there is single Jordan block of order two corresponding to an eigenvalue $u$. Thus, a Jordan chain of order two will be formed by matrix $\boldsymbol{A}_{c}^{\boldsymbol{TV}}$, {\em i.e.},
\begin{align}\label{gen_for_TV}
 \begin{split}
  \boldsymbol{A}_{c}^{\boldsymbol{TV}}\boldsymbol{X}_{1} \ &= \ u\boldsymbol{X}_{1} \ \textrm{and} \   \\
  \boldsymbol{A}_{c}^{\boldsymbol{TV}}\boldsymbol{X}_{2} \ &= \ u\boldsymbol{X}_{2} \ + \ \boldsymbol{X}_{1}
 \end{split}
\end{align} should hold.
From first relation we get 
\begin{equation}
 \boldsymbol{X}_{1} \ = \ \boldsymbol{R}_{c,2}^{\boldsymbol{TV}} =  \begin{bmatrix}
         \ 1 \\[0.3em]
         \ u \\[0.3em]
         \ \frac{1}{2} u^2
       \end{bmatrix}
\end{equation} 
and on using  $\boldsymbol{X}_{1}$ in the second relation of (\ref{gen_for_TV}), we can find required generalized eigenvector $\boldsymbol{X}_{2}$ which is given below.
\begin{equation}
 \boldsymbol{X}_{2} \ = \ \boldsymbol{R}_{c,3}^{\boldsymbol{TV}} =  \begin{bmatrix}
         \ x_{1} \\[0.3em]
         \ 1 + ux_{1} \\[0.3em]
         \ u + \frac{1}{2}u^{2}x_{1}
       \end{bmatrix}
\end{equation}
Here, $x_{1} \in {\rm I\!R}$ is a real constant and $det(\boldsymbol{P})$ is equal to one. If we take $\boldsymbol{P}$ equal to  
\begin{equation}
  \begin{bmatrix}
        \ 0  &&  1  &&       x_{1}   \\[0.3em]
        \ 0  &&  u  &&       1 + ux_{1}  \\[0.3em]
        \ 1  && \frac{1}{2}u^{2}  && u + \frac{1}{2}u^{2}x_{1}
        \end{bmatrix} \  \textrm{then}  \  
  \boldsymbol{P}^{-1}\boldsymbol{A}_{c}^{\boldsymbol{TV}}\boldsymbol{P} \ = \  \begin{bmatrix}
        \ 0  &&  0  &&  0   \\[0.3em]
        \ 0  &&  u  &&  1   \\[0.3em]
        \ 0  &&  0  &&  u 
        \end{bmatrix} \ \ = \ \boldsymbol{J}_1 
 \end{equation}
and if we take $\boldsymbol{P}$  equal to
\begin{equation}
 \begin{bmatrix}
        \ 1  &&  x_{1}  && 0  \\[0.3em]
        \ u  &&  1 + ux_{1}  && 0 \\[0.3em]
        \ \frac{1}{2}u^{2}  && u + \frac{1}{2}u^{2}x_{1}  && 1
        \end{bmatrix} \  \textrm{then}  \
 \boldsymbol{P}^{-1}\boldsymbol{A}_{c}^{\boldsymbol{TV}}\boldsymbol{P} \ = \ \begin{bmatrix}
        \ u  && 1  && 0  \\[0.3em]
        \ 0  && u  && 0 \\[0.3em]
        \ 0  && 0  && 0 
        \end{bmatrix} \ \ = \ \boldsymbol{J}_2
\end{equation}
\\
Next, we have to find generalized eigenvectors corresponding to Jacobian matrix
$\boldsymbol{A}_{c}^{\boldsymbol{ZB}}$ for a convective subsystem of Zha and Bilgen type splitting. As explained earlier, eigenvalues for $\boldsymbol{A}_{c}^{\boldsymbol{ZB}}$ are $u, u, u$ and set of LI eigenvectors are,
 \begin{equation}
      \boldsymbol{R}_{c,1}^{\boldsymbol{ZB}}  =  \begin{bmatrix}
                  \ 1    \\[0.3em]
                  \ u     \\[0.3em]
                  \ 0 
                 \end{bmatrix}  ~~\mbox{and}~~ \
      \boldsymbol{R}_{c,2}^{\boldsymbol{ZB}} =  \begin{bmatrix}
                  \ 0    \\[0.3em]
                  \ 0     \\[0.3em]
                  \ 1 
                 \end{bmatrix} \          
  \end{equation}
On computing ranks of matrices $(\boldsymbol{A}_{c}^{\boldsymbol{ZB}}- u\boldsymbol{I})$, \  \ $(\boldsymbol{A}_{c}^{\boldsymbol{ZB}}- u\boldsymbol{I})^2$, $\cdots$, we find  $rank(\boldsymbol{A}_{c}^{\boldsymbol{ZB}}- u\boldsymbol{I})^2  \ =  \ 0  \ =  \  rank(\boldsymbol{A}_{c}^{\boldsymbol{ZB}}- u\boldsymbol{I})^{3}$. Thus there will be one $Jordan$ block of order $2$ and since all eigenvalues are equal then there must be another Jordan block of order one. In short, Jordan matrix $\boldsymbol{J}$ corresponding to matrix $\boldsymbol{A}_c^{\boldsymbol{ZB}}$ is made up of two Jordan blocks, which clearly shows that for the present case there is no single Jordan block for given set of eigenvalues.  Thus, earlier theory may not be directly applicable for this case.  But a Jordan chain of order two should form by matrix $\boldsymbol{A}_{c}^{\boldsymbol{ZB}}$. If possible, without loss of generality, let first assume 
\begin{align}
 \begin{split}
  \boldsymbol{A}_{c}^{\boldsymbol{ZB}}\boldsymbol{R}_{c,1}^{\boldsymbol{ZB}} \ &= \ u\boldsymbol{R}_{c,1}^{\boldsymbol{ZB}} \\
  \boldsymbol{A}_{c}^{\boldsymbol{ZB}}\boldsymbol{X} \ &= \ u\boldsymbol{X} + \boldsymbol{R}_{c,1}^{\boldsymbol{ZB}}
 \end{split}
\end{align}
holds. On expanding second relation 
\begin{equation*}
 \begin{bmatrix}
         \ 0  &&  1  &&   0   \\[0.3em]
       \ -u^{2}  &&  2 u  &&  0  \\[0.3em]
       \ -u E  &&    E &&  u
        \end{bmatrix} \  \begin{bmatrix}
                  \ x_{1}    \\[0.3em]
                  \  x_{2}    \\[0.3em]
                  \  x_{3} 
                 \end{bmatrix} \ \textrm{=} \  u \begin{bmatrix}
                  \ x_{1}    \\[0.3em]
                  \  x_{2}    \\[0.3em]
                  \  x_{3} 
                 \end{bmatrix}  \ + \ \begin{bmatrix}
                  \ 1    \\[0.3em]
                  \ u     \\[0.3em]
                  \ 0 
                 \end{bmatrix}
\end{equation*}
and from first two equations, we get
\begin{align}
\begin{split} 
   x_{2} \ &= \ ux_{1} + 1
\end{split}
\end{align}
Similarly, from third equation 
\begin{align}
 \begin{split}
  -u Ex_{1} \ + \ Ex_{2} \ + \ ux_{3} \ &= \  ux_{3} \\
  \Rightarrow x_{2} \ &= \  ux_{1}.
\end{split}
\end{align}
We get two different expressions for real constant $x_{2}$, which is a contradiction. Thus eigenvector $\boldsymbol{R_{c,1}^{ZB}}$ can't form a Jordan chain of order two corresponding to matrix $\boldsymbol{A}_{c}^{\boldsymbol{ZB}}$. If possible, let us assume now $\boldsymbol{R}_{c,2}^{\boldsymbol{ZB}}$ forms a Jordan chain of order two, {\em i.e.},
\begin{align}
 \begin{split}
  \boldsymbol{A}_{c}^{\boldsymbol{ZB}}\boldsymbol{R}_{c,2}^{\boldsymbol{ZB}} \ &= \ u\boldsymbol{R}_{c,2}^{\boldsymbol{ZB}} \\
  \boldsymbol{A}_{c}^{\boldsymbol{ZB}}\boldsymbol{X} \ &= \ u\boldsymbol{X} + \boldsymbol{R}_{c,2}^{\boldsymbol{ZB}}
 \end{split}
\end{align}
holds. Again after expanding second relation, we have
\begin{equation*}
 \begin{bmatrix}
         \ 0  &&  1  &&   0   \\[0.3em]
       \ -u^{2}  &&  2 u  &&  0  \\[0.3em]
       \ -u E  &&    E &&  u
        \end{bmatrix} \  \begin{bmatrix}
                  \ x_{1}    \\[0.3em]
                  \  x_{2}    \\[0.3em]
                  \  x_{3} 
                 \end{bmatrix} \ \textrm{=} \  u \begin{bmatrix}
                  \ x_{1}    \\[0.3em]
                  \  x_{2}    \\[0.3em]
                  \  x_{3} 
                 \end{bmatrix}  \ + \ \begin{bmatrix}
                  \ 0    \\[0.3em]
                  \ 0     \\[0.3em]
                  \ 1 
                 \end{bmatrix}
\end{equation*}
and further on solving first two equations of expanded system we get
\begin{align}
\begin{split}
 x_{2} \ &= \ ux_{1}.
\end{split}
\end{align}
From third equation we have
\begin{align}
 \begin{split}
  -u Ex_{1} \ + \ Ex_{2} \ + \ ux_{3} \ &= \  ux_{3} + 1 \\
  \Rightarrow x_{2} \ &= \ \frac{1 + u Ex_{1}}{E}
\end{split}
\end{align}
If we compare both values of $x_{2}$, we get $0 = 1$ which is impossible, hence a contradiction. Therefore, neither $\boldsymbol{R}_{c,1}^{\boldsymbol{ZB}}$ nor $\boldsymbol{R}_{c,2}^{\boldsymbol{ZB}}$ helps in forming a Jordan chain of order two corresponding to matrix $\boldsymbol{A}_{c}^{\boldsymbol{ZB}}$. Thus, we need to go more deep into the theory of Jordan canonical forms to obtain proper generalized eigenvectors.  Since there will be a Jordan block of order 2, this means that we need to construct a generalized eigenvector which should help in generating a {\em Jordan chain}.    
Let $R(\boldsymbol{A})$ denote the space spanned by the columns of matrix $\boldsymbol{A}_{c}^{\boldsymbol{ZB}} -u\boldsymbol{I}$. Then 
\begin{equation}
 R(\boldsymbol{A}) \ = \ x_{1}\boldsymbol{A}_1 \ + \ x_{2}\boldsymbol{A}_2 \ + \ x_{3}\boldsymbol{A}_3
\end{equation}
where, $\boldsymbol{A}_1,\boldsymbol{A}_2,\boldsymbol{A}_3$ are column vectors of $\boldsymbol{A_{c}^{ZB}} -u\boldsymbol{I}$. Now
\begin{equation}
 R(\boldsymbol{A}) \ = \ x_{1}\begin{bmatrix}
                  \ -u    \\[0.3em]
                  \  -u^{2}    \\[0.3em]
                  \  -u E 
                 \end{bmatrix} \ + \ 
               x_{2}\begin{bmatrix}
                  \ 1    \\[0.3em]
                  \ u    \\[0.3em]
                  \ E 
                 \end{bmatrix}  \ + \ 
                 x_{3}\begin{bmatrix}
                  \ 0    \\[0.3em]
                  \ 0     \\[0.3em]
                  \ 0 
                 \end{bmatrix}
\end{equation} 
or
\begin{equation}
 R(\boldsymbol{A}) \ = \ -u x_{1}\begin{bmatrix}
                  \ 1    \\[0.3em]
                  \ u    \\[0.3em]
                  \ E 
                 \end{bmatrix} \ + \ 
               x_{2}\begin{bmatrix}
                  \ 1    \\[0.3em]
                  \ u    \\[0.3em]
                  \ E 
                 \end{bmatrix}  \ = \ (-u x_{1} + x_{2})
                 \begin{bmatrix}
                  \ 1    \\[0.3em]
                  \ u    \\[0.3em]
                  \ E
                  \end{bmatrix}
\end{equation} 
Therefore, column vector $\boldsymbol{X} \ = \ (1,u,E)^{t}$ is a range space of $R(\boldsymbol{A})$. Next, let $N(\boldsymbol{AX})$ be a null space of column vectors $(\boldsymbol{A}_{c}^{\boldsymbol{ZB}} -u\boldsymbol{I})\boldsymbol{X}$. By definition 
\begin{equation}
 N(\boldsymbol{AX}) \ = \ \big\{\ \boldsymbol{v} \in {\rm I\!R^n} ;  (\boldsymbol{AX})\boldsymbol{v} = \boldsymbol{0} \big\}
\end{equation} and dimension of $\boldsymbol{v}$ is equal to number of  entries in each row of matrix $\boldsymbol{AX}$.  
Now 
\begin{equation}
 \boldsymbol{AX} \ = \ (\boldsymbol{A_{c}^{ZB}} -u\boldsymbol{I})\boldsymbol{X}  \ = \ \begin{bmatrix}
         \ -u  &&  1  &&   0   \\[0.3em]
       \ -u^{2}  &&  u  &&  0  \\[0.3em]
       \ - u E  &&   E &&  0
        \end{bmatrix} \ 
        \begin{bmatrix}
                  \ 1    \\[0.3em]
                  \ u    \\[0.3em]
                  \ E 
                 \end{bmatrix} \ = \ 
                 \begin{bmatrix}
                  \ 0    \\[0.3em]
                  \ 0    \\[0.3em]
                  \ 0 
                 \end{bmatrix}. \ 
\end{equation}
For present case, $\boldsymbol{AX}$ is just a null vector, therefore $\boldsymbol{v}$ reduces to a scalar coefficient. By definition of null space of column vectors we have, 
  \begin{equation}
                 \begin{bmatrix}
                  \ 0    \\[0.3em]
                  \ 0    \\[0.3em]
                  \ 0 
                 \end{bmatrix}v \ = \ 
                 \begin{bmatrix}
                  \ 0    \\[0.3em]
                  \ 0    \\[0.3em]
                  \ 0 
                 \end{bmatrix} \ 
 \end{equation}
which holds for any $v \in {\rm I\!R}$. And by definition, $\boldsymbol{Xv}$ which is equal to $\boldsymbol{X}v$ should from a basis for $ R(\boldsymbol{A}_{c}^{\boldsymbol{ZB}} -u\boldsymbol{I})$ $\cap$ $N(\boldsymbol{A}_{c}^{\boldsymbol{ZB}} -u\boldsymbol{I})$. Thus, $\boldsymbol{X}_1 \ = \ (1,u,E)^{t}$ should be a generalized eigenvector and to check that we need to see whether for $\boldsymbol{X} \ = \ \boldsymbol{X}_1$, relation $\boldsymbol{A}_{c}^{\boldsymbol{ZB}}\boldsymbol{X} \ = \ u\boldsymbol{X}$ holds or not.
Now
\begin{equation}
  \boldsymbol{A}_{c}^{\boldsymbol{ZB}}\boldsymbol{X}_1  \ = \ \begin{bmatrix}
       \ 0 &&  1  &&   0   \\[0.3em]
       \ -u^{2}  &&  2u  &&  0  \\[0.3em]
       \ - u E  &&   E &&  u
        \end{bmatrix} \ 
        \begin{bmatrix}
                  \ 1    \\[0.3em]
                  \ u    \\[0.3em]
                  \ E 
                 \end{bmatrix} \ = \ 
                u \begin{bmatrix}
                  \ 1    \\[0.3em]
                  \ u    \\[0.3em]
                  \ E 
                 \end{bmatrix} \ = \ u \boldsymbol{X}_1 \
\end{equation}
This generalized eigenvector is expected to form a Jordan chain of order two corresponding to matrix $\boldsymbol{A}_{c}^{\boldsymbol{ZB}}$, {\em i.e.},
\begin{align}
 \begin{split}
 \boldsymbol{A}_{c}^{\boldsymbol{ZB}}\boldsymbol{X}_1 \ &= \ u\boldsymbol{X}_1  \\
  \boldsymbol{A}_{c}^{\boldsymbol{ZB}}\boldsymbol{X}_2 \ &= \ u\boldsymbol{X}_2 \ + \  \boldsymbol{X}_1
 \end{split}
\end{align}
Eigenvector $\boldsymbol{X}_2$ can be find from second relation and in expanded from it is written as,
\begin{equation}
 \begin{bmatrix}
         \ 0  &&  1  &&   0   \\[0.3em]
       \ -u^{2}  &&  2u  &&  0  \\[0.3em]
       \ - u E  &&   E &&  u
        \end{bmatrix} \begin{bmatrix}
                  \ x_{1}    \\[0.3em]
                  \ x_{2}    \\[0.3em]
                  \ x_{3}
                 \end{bmatrix}  \ = \ u\begin{bmatrix}
                  \ x_{1}    \\[0.3em]
                  \ x_{2}    \\[0.3em]
                  \ x_{3}
                 \end{bmatrix}        \ + \   
                  \begin{bmatrix}
                  \ 1    \\[0.3em]
                  \ u    \\[0.3em]
                  \ E 
                 \end{bmatrix} 
\end{equation}
After little algebra, $\boldsymbol{X}_2$ comes out as
\begin{equation}
 \boldsymbol{X}_2  \ = \ \begin{bmatrix}
                  \ x_{1}    \\[0.3em]
                  \ 1 + ux_{1}    \\[0.3em]
                  \ x_{3}
                 \end{bmatrix} \ = \ \boldsymbol{R}_{c,2}^{\boldsymbol{ZB}}
 \end{equation}
 where $x_1, x_3 \in {\rm I\!R}$.
 \begin{equation}
 \boldsymbol{X}_1  \ = \ \begin{bmatrix}
                  \ 1    \\[0.3em]
                  \ u    \\[0.3em]
                  \ E
                 \end{bmatrix} \ = \ \boldsymbol{R}_{c,1}^{\boldsymbol{ZB}}
 \end{equation}
 and we can take $\boldsymbol{R}_{c,3}^{\boldsymbol{ZB}}$ either equal to
\begin{equation}
\begin{bmatrix}
                  \ 0    \\[0.3em]
                  \ 0    \\[0.3em]
                  \ 1
                 \end{bmatrix} \  \ \textrm{or} \ \
                 \begin{bmatrix}
                  \ 1   \\[0.3em]
                  \ u    \\[0.3em]
                  \ 0
                 \end{bmatrix}.
\end{equation}
If we take 
\begin{equation}
 \boldsymbol{R}_{c,3}^{\boldsymbol{ZB}} \ = \ \begin{bmatrix}
                  \ 0    \\[0.3em]
                  \ 0    \\[0.3em]
                  \ 1
                 \end{bmatrix} \
\ \textrm{then} \
 \boldsymbol{P} \ = \ \begin{bmatrix}
       \ 1  &  x_{1}  &   0   \\[0.3em]
       \ u  &  1 + ux_{1}  &  0  \\[0.3em]
       \ E  &   x_{3} &  1
        \end{bmatrix} 
\end{equation}
 and $det(\boldsymbol{P}) \ = \ 1$. Further,
  \begin{equation}
        \boldsymbol{P}^{-1}\boldsymbol{A}_{c}^{\boldsymbol{ZB}}\boldsymbol{P} \ = \ \ \begin{bmatrix}
       \ u  &  1  & 0   \\[0.3em]
       \ 0  &  u  & 0  \\[0.3em]
       \ 0  &  0 &  u
        \end{bmatrix} \ = \ \boldsymbol{J}_{1}
  \end{equation}
Similarly, if we take
\begin{equation}
 \boldsymbol{P} \ = \ \begin{bmatrix}
       \ 0 & 1  &  x_{1}     \\[0.3em]
       \ 0 & u  &  1 + ux_{1}  \\[0.3em]
       \ 1 & E  &  x_{3}
        \end{bmatrix} \ \textrm{then} \ 
        \boldsymbol{P}^{-1}\boldsymbol{A}_{c}^{\boldsymbol{ZB}}\boldsymbol{P} \ = \ \ \begin{bmatrix}
       \ u  &  0 & 0   \\[0.3em]
       \ 0  &  u & 1  \\[0.3em]
       \ 0  &  0 &  u
        \end{bmatrix} \ = \ \boldsymbol{J}_{2}
\end{equation}
 Similarly, let 
 \begin{equation}
 \boldsymbol{P} \ = \ \begin{bmatrix}
       \ 1  &  x_{1}  &   1   \\[0.3em]
       \ u  &  1 + ux_{1}  &  u  \\[0.3em]
       \ E  &   x_{3} &  0
        \end{bmatrix} 
\end{equation} then, $det(\boldsymbol{P}) \ = \ -E \neq 0$.
\section{Formulation of ZBS-FDS and TVS-FDS schemes}
\subsection{ZBS-FDS scheme}
We first consider pressure subsystem of Zha and Bilgen type splitting and on comparing  (\ref{conservation_ZB_p_part}) and (\ref{quasi_form_ZB_p_part}), we get 
\begin{equation}
 d\boldsymbol{F}_{p}^{\boldsymbol{ZB}} \ = \ \boldsymbol{A}_{p}^{\boldsymbol{ZB}}d\boldsymbol{U}
\end{equation}
The finite difference analogue of the above differential relation is,
\begin{equation}\label{FDS_ZBS_p_part}
 \Delta{\boldsymbol{F}_{p}^{\boldsymbol{ZB}}} \ = \ \boldsymbol{\bar{A}}_{p}^{\boldsymbol{ZB}}\Delta{\boldsymbol{U}}
\end{equation}
where $\boldsymbol{\bar{A}}_{p}^{\boldsymbol{ZB}}$ is now a function of left and right states, {\em i.e.}, $\boldsymbol{\bar{A}}_{p}^{\boldsymbol{ZB}} = \boldsymbol{\bar{A}}_{p}^{\boldsymbol{ZB}}(\boldsymbol{U}_L,\boldsymbol{U}_R)$. Since $\Delta{\boldsymbol{U}}$ is a column vector, therefore it can be written as linear combination of LI eigenvectors.
\begin{equation}
 \Delta \boldsymbol{U}   \ = \  \sum_{i = 1}^{3} \bar{\alpha}_{p,i}^{\boldsymbol{ZB}}\boldsymbol{\bar{R}}_{p,i}^{\boldsymbol{ZB}} 
\end{equation}
On using above expression in (\ref{FDS_ZBS_p_part}),
\begin{equation}
 \\    \Delta{\boldsymbol{F}_{p}^{\boldsymbol{ZB}}} \ = \  \boldsymbol{\bar{A}}_{p}^{\boldsymbol{ZB}} \sum_{i = 1}^{3} \bar{\alpha}_{p,i}^{\boldsymbol{ZB}}\boldsymbol{{\bar{R}}}_{p,i}^{\boldsymbol{ZB}}
\end{equation}
or
\begin{equation}
\Delta{\boldsymbol{F}_p^{\boldsymbol{ZB}}}   \ = \   \bar{\alpha}_{p,1}^{\boldsymbol{ZB}} \boldsymbol{\bar{A}}_{p}^{\boldsymbol{ZB}} \boldsymbol{{\bar{R}}}_{p,1}^{\boldsymbol{ZB}}  \ + \ \bar{\alpha}_{p,2}^{\boldsymbol{ZB}} \boldsymbol{\bar{A}}_{p}^{\boldsymbol{ZB}} \boldsymbol{{\bar{R}}}_{p,2}^{\boldsymbol{ZB}}  \ + \ \bar{\alpha}_{p,3}^{\boldsymbol{ZB}} \boldsymbol{\bar{A}}_{p}^{\boldsymbol{ZB}} \boldsymbol{{\bar{R}}}_{p,3}^{\boldsymbol{ZB}}
\end{equation} 
which is further equal to
\begin{equation} 
\Delta{\boldsymbol{F}_p^{\boldsymbol{ZB}}}   \ = \   \bar{\alpha}_{p,1}^{\boldsymbol{ZB}}\bar \lambda_{p,1}^{\boldsymbol{ZB}} \boldsymbol{{\bar{R}}}_{p,1}^{\boldsymbol{ZB}}  \ + \ \bar{\alpha}_{p,2}^{\boldsymbol{ZB}} \bar \lambda_{p,2}^{\boldsymbol{ZB}} \boldsymbol{{\bar{R}}}_{p,2}^{\boldsymbol{ZB}}  \ + \ \bar{\alpha}_{p,3}^{\boldsymbol{ZB}} \bar \lambda_{p,3}^{\boldsymbol{ZB}} \boldsymbol{{\bar{R}}}_{p,3}^{\boldsymbol{ZB}}
\end{equation} 
Now, $\Delta{\boldsymbol{F}_p^{+\boldsymbol{ZB}}}$ must have contribution of positive part of eigenvalues only, {\em i.e.},
\begin{align}\label{ZBS-FDS_positive_p_part}
 \begin{split}
  \Delta{\boldsymbol{F}_p^{+\boldsymbol{ZB}}}   \ = \  \bar{\alpha}_{p,1}^{\boldsymbol{ZB}}\bar \lambda_{p,1}^{+\boldsymbol{ZB}} \boldsymbol{{\bar{R}}}_{p,1}^{\boldsymbol{ZB}}  \ &+ \ \bar{\alpha}_{p,2}^{\boldsymbol{ZB}} \bar \lambda_{p,2}^{+\boldsymbol{ZB}} \boldsymbol{{\bar{R}}}_{p,2}^{\boldsymbol{ZB}}  \ + \ \bar{\alpha}_{p,3}^{\boldsymbol{ZB}} \bar \lambda_{p,3}^{+\boldsymbol{ZB}} \boldsymbol{{\bar{R}}}_{p,3}^{\boldsymbol{ZB}} 
 \end{split}
\end{align}
Similarly,
\begin{align}\label{ZBS-FDS_negative_p_part}
 \begin{split}
  \Delta{\boldsymbol{F}_p^{-\boldsymbol{ZB}}}   \ = \  \bar{\alpha}_{p,1}^{\boldsymbol{ZB}}\bar \lambda_{p,1}^{-\boldsymbol{ZB}} \boldsymbol{{\bar{R}}}_{p,1}^{\boldsymbol{ZB}}  \ &+ \ \bar{\alpha}_{p,2}^{\boldsymbol{ZB}} \bar \lambda_{p,2}^{-\boldsymbol{ZB}} \boldsymbol{{\bar{R}}}_{p,2}^{\boldsymbol{ZB}}  \ + \ \bar{\alpha}_{p,3}^{\boldsymbol{ZB}} \bar \lambda_{p,3}^{-\boldsymbol{ZB}} \boldsymbol{{\bar{R}}}_{p,3}^{\boldsymbol{ZB}} 
 \end{split}
\end{align}
We now define the standard Courant splitting for the eigenvalues as 
 \begin{equation}
 \bar \lambda^{\pm}_{p,i} = \frac{\bar \lambda_{p,i} \pm |\bar \lambda_{p,i}|}{2} 
 \end{equation}
 On using standard upwinding along with above definition, we finally get
\begin{equation}\label{ZBS-FDS-p-part}
 \Delta{\boldsymbol{F}_p^{+\boldsymbol{ZB}}} - \Delta{\boldsymbol{F}_p^{-\boldsymbol{ZB}}}  \ = \ \sum_{i = 1}^{3} \bar{\alpha}_{p,i}^{\boldsymbol{ZB}}\left|\bar \lambda_{p,i}^{\boldsymbol{ZB}}\right| \boldsymbol{{\bar{R}}}_{p,i}^{\boldsymbol{ZB}}
\end{equation} 
To determine right side of (\ref{ZBS-FDS-p-part}) completely, we need to find average values of eigenvalues along with coefficients which are attached with LI eigenvectors. First we consider linearization equation, $\Delta \boldsymbol{F}_{p}^{\boldsymbol{ZB}}  \ = \ \boldsymbol{\bar{A}}_{p}^{\boldsymbol{ZB}}\Delta{\boldsymbol{U}}$, of pressure subsystem of ZBS-FDS scheme. In expanded form it can be written as
\begin{equation}
  \begin{bmatrix}
         0 \\[0.3em]
       \Delta(p) \\[0.3em]
       \Delta(pu)
         \end{bmatrix} \ = \
    \begin{bmatrix}
         \ 0  & 0 & 0  \\[0.3em]
         \frac{1}{2} (\gamma -1) \bar u^2  &  -(\gamma -1)\bar u  &  (\gamma -1) \\[0.3em]
         \ -\frac{\bar u \bar a^2}{\gamma} + \frac{(\gamma - 1)}{2} \bar u^3  &  \frac{\bar a^2}{\gamma} - (\gamma - 1) \bar u^2  &  (\gamma - 1)\bar u 
        \end{bmatrix} 
        \begin{bmatrix}
         \Delta (\rho)  \\[0.3em]
       \Delta  (\rho u) \\[0.3em]
       \Delta  (\rho E)
        \end{bmatrix}
       \end{equation}
From the second equation, we get
\begin{equation}
 \Delta p   \  =  \  \frac{1}{2} (\gamma -1) \bar u^2 \Delta \rho  \  -  \ (\gamma -1)\bar u \Delta  (\rho u)  \  +  \  (\gamma -1) \Delta  (\rho E)
\end{equation}
or
\begin{align}
\begin{split}
\Delta p   \  =  \  \frac{1}{2} (\gamma -1) \bar u^2 \Delta \rho  \  -  \ (\gamma -1)\bar u \Delta  (\rho u)  \  +  \  (\gamma -1) \Delta(\frac{p}{\gamma - 1}) \\ \ + \ \frac{1}{2} (\gamma - 1) \Delta(\rho u^2)
\end{split}
\end{align}
It further reduces to 
\begin{align}\label{average_rel_1}
 \begin{split}
     {\bar{u}}^2 \Delta (\rho) \ - \ 2 \bar{u} \Delta  (\rho u) \ + \ \Delta (\rho u^{2})\ = \ 0
 \end{split}
\end{align}
which gives average value for conserved variable $\bar u$ as given below.
\begin{equation}
 \\  \  \bar{u}    \ = \  \frac{\sqrt{\rho_L} u_L  \ + \ \sqrt{\rho_R} u_R }{\sqrt{\rho_L} \ + \ \sqrt{\rho_R}}
\end{equation} 
Other root is being neglected as it contains negative sign in the denominator. Let us consider the relation
\begin{equation}\label{average_rel_2}
 \Delta(\rho u) \ = \ \bar{\rho} \Delta{u} \ + \ \bar{u} \Delta{\rho}
\end{equation} in expanded form it is written as,
\begin{equation}
 \rho_R u_R - \rho_L u_L   \ = \  \bar{\rho}(u_R - u_L) \ + \  \bar{u}(\rho_R - \rho_L)
\end{equation}
On using $\bar{u}$ in the above relation we get $\bar{\rho} = \sqrt{\rho_L\rho_R}$. 

\begin{figure}[!ht]
 \centerline{%
\subfigure[]{%
\includegraphics[trim=0 5 35 5, clip, width=0.55\textwidth]{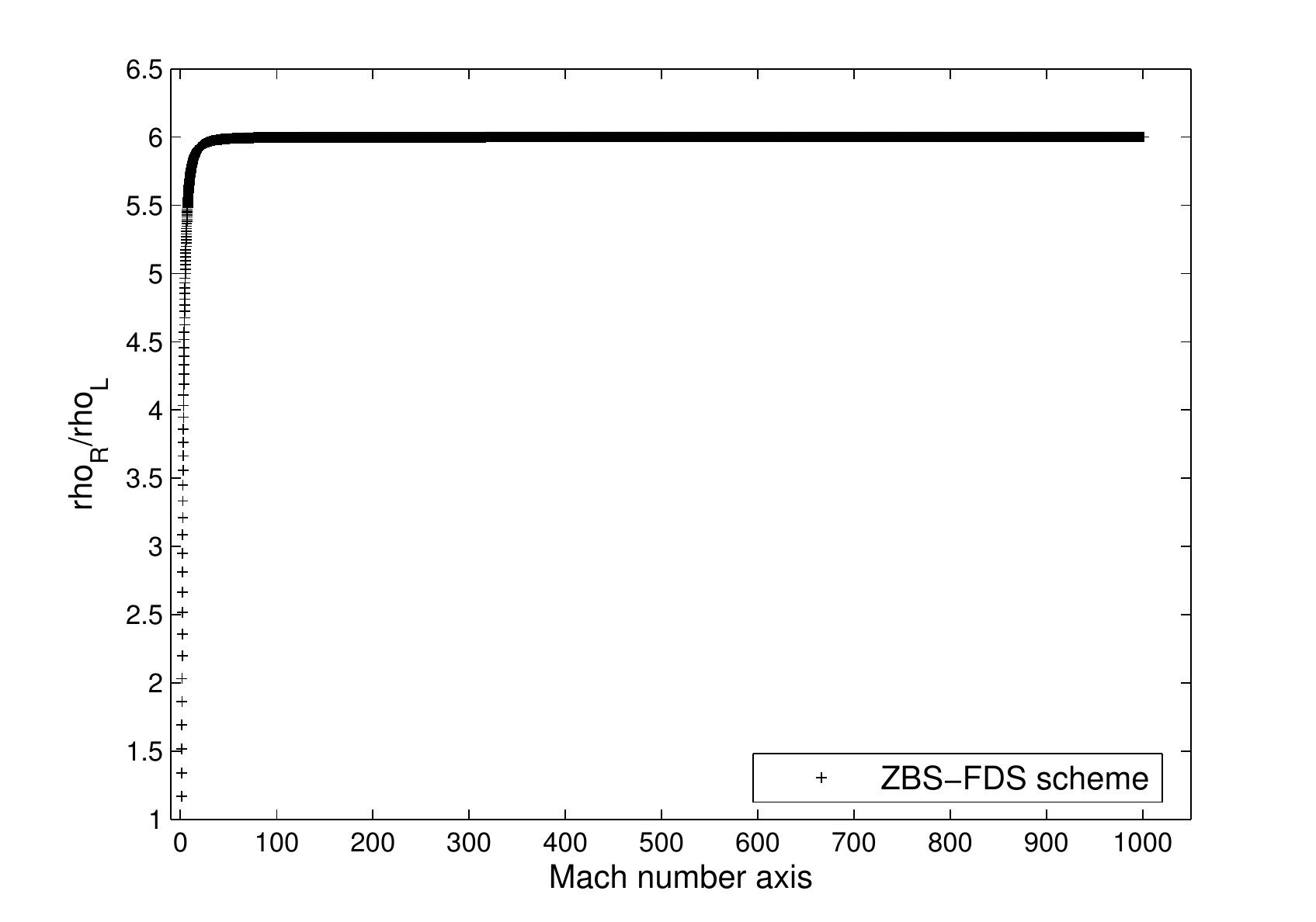}%
\label{ZBS-FDS_d_ratio.eps}}%
\subfigure[]{%
\includegraphics[trim=0 5 35 5, clip, width=0.55\textwidth]{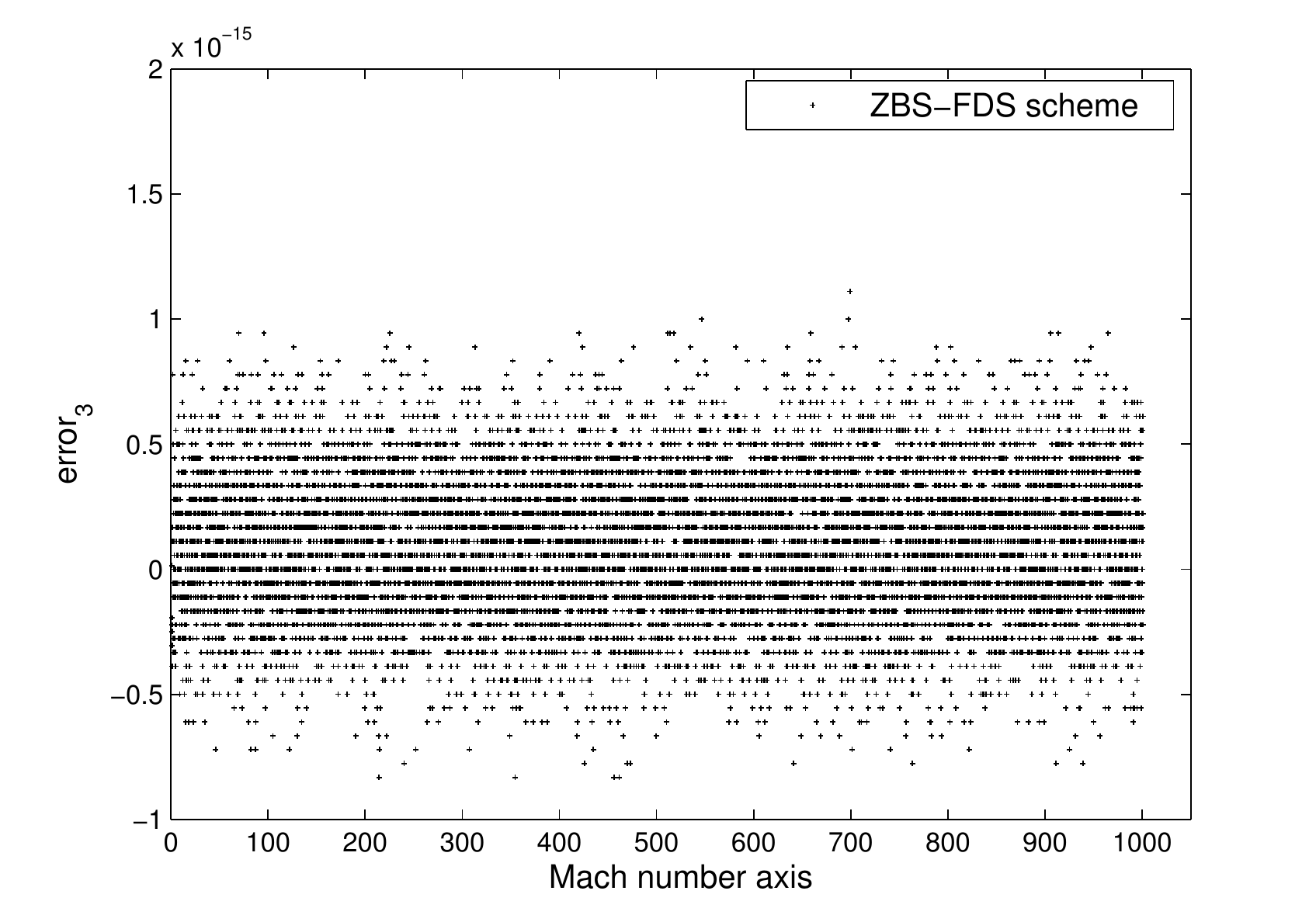}%
\label{ZBS-FDS_error-3.pdf}
}
}%
\caption{(a) represents ratio of densities results and (b) represents error$_3$ results at microscopic level}
\end{figure}

\begin{figure}[!ht]
\begin{center}
\includegraphics[trim=5 5 35 5, clip, width=0.6\textwidth]{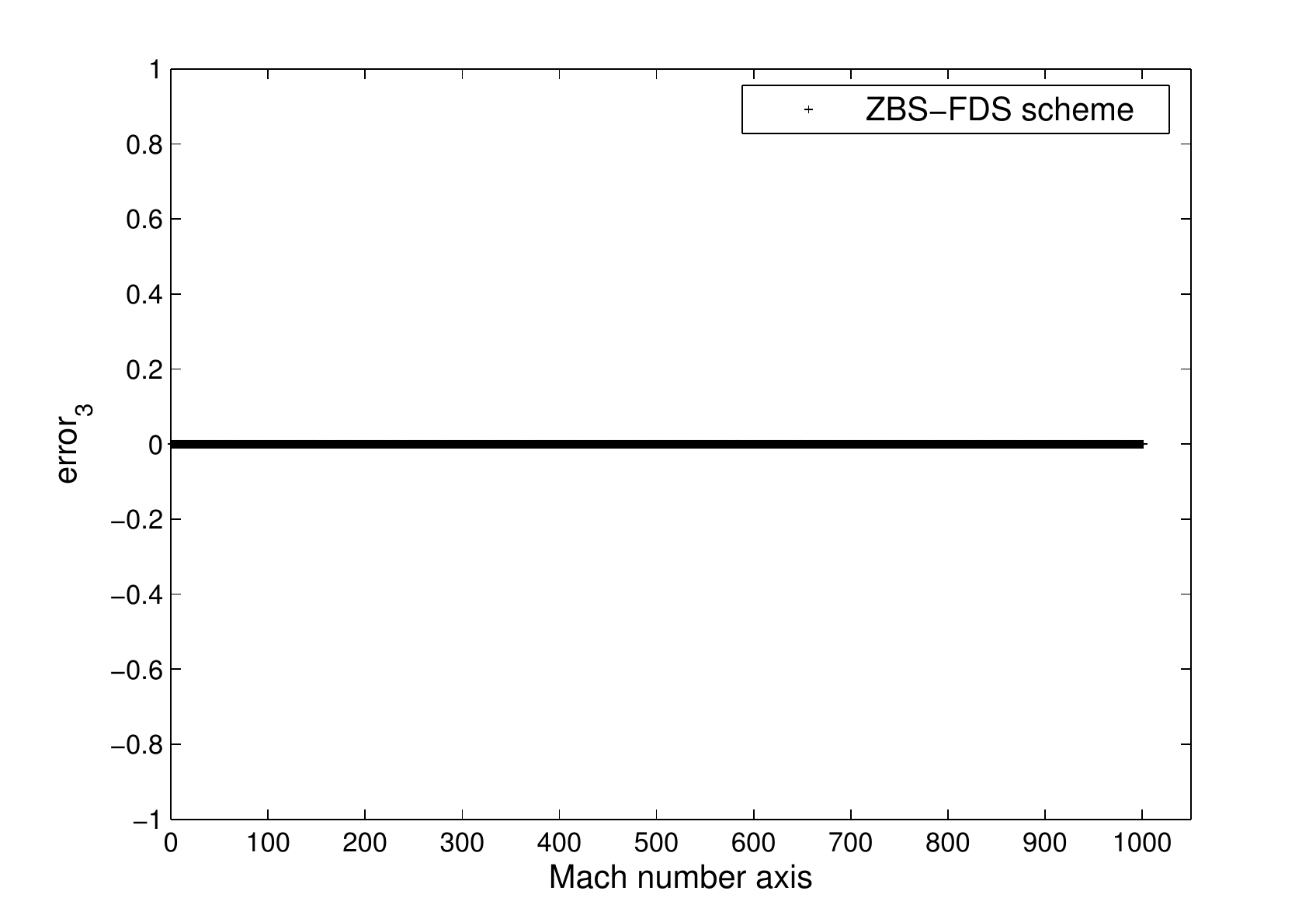}
\caption{represents error$_3$ results at macroscopic level}
\label{ZBS-FDS_error-3_macro.eps}
\end{center}
\end{figure}

Similarly, third equation can be written as 
\begin{eqnarray}
 \Delta(pu) \ & = &  \  \left(-\frac{\bar u \bar a^2}{\gamma} + \frac{(\gamma - 1)}{2} \bar u^3\right)\Delta\rho \  \nonumber \\ & & +  \   \left(\frac{\bar a^2}{\gamma} - (\gamma - 1) \bar u^2\right)\Delta (\rho u)  + \ (\gamma - 1)\bar u \Delta( \rho E)
\end{eqnarray} 
On expanding further, above equation looks like
\begin{align}
 \begin{split}
  \Delta(pu) \ &= \  \left(-\frac{\bar u \bar a^2}{\gamma} + \frac{(\gamma - 1)}{2} \bar u^3\right)\Delta\rho  \ + \   \left(\frac{\bar a^2}{\gamma} - (\gamma - 1) \bar u^2\right)\Delta{(\rho u)} \\ 
  \ &+ \ \bar u\Delta(p) + \frac{1}{2}(\gamma - 1)\bar{u}\Delta{(\rho u^{2})}
 \end{split}
\end{align}
We next use (\ref{average_rel_1}) and (\ref{average_rel_2}) in above equation and after cancellations of some terms we get,
\begin{equation}
 \Delta (pu)   \  =  \  \bar u \Delta p    \  +   \   \bar p \Delta u   
\end{equation}
On using $\bar p =  \dfrac{\bar{a}^{2}\bar{\rho}}{\gamma}$, we have
\begin{equation}\label{average_rel_3}
\Delta(a^2 \rho u)   \  -  \  \bar{u}\Delta(a^2\rho) \  =  \  \bar{a}^2\bar{\rho} \Delta{u}
\end{equation}
Let $\eta$ be any flow variable.  Consider the relation
\begin{equation}
 \Delta (\rho u \eta)  \ -  \ \bar u \Delta (\rho \eta)  \ = \ \bar \rho \bar \eta \Delta u  
\end{equation}
which can be written as 
\begin{equation}
 \Delta (\rho u \eta)  \ -  \ \bar u \Delta (\rho \eta)  \ = \ \sqrt{\rho_L \rho_R}( \frac{\sqrt{\rho_L} \eta_L  \ + \ \sqrt{\rho_R} \eta_R }{\sqrt{\rho_L} \ + \ \sqrt{\rho_R}}) \Delta u
\end{equation}
We put $\eta = a^2$, then above relation becomes an equation and further on comparing it with (\ref{average_rel_3}), we get average value of acoustic speed from  
\begin{equation}
 \bar{a}^2  \ = \ \frac{\sqrt{\rho_L} a^{2}_L  \ + \ \sqrt{\rho_R} a^{2}_R }{\sqrt{\rho_L} \ + \ \sqrt{\rho_R}}
\end{equation}
Similarly, we find coefficients which are attached with LI eigenvector of pressure subsystem for ZBS-FDS scheme as follows.  
\begin{equation}
 \Delta \boldsymbol{U}   \ = \  \sum_{i = 1}^{3} \bar{\alpha}_{p,i}^{\boldsymbol{ZB}}\boldsymbol{{\bar{R}}}_{p,i}^{\boldsymbol{ZB}}
\end{equation}
In expanded form, 
\begin{equation} 
  \begin{bmatrix} 
                 \Delta(\rho) \\[0.3em]  
		\Delta(\rho u) \\[0.3em] 
		\Delta(\rho E) 
		\end{bmatrix}   = 
   \bar{\alpha}_{p,1}^{\boldsymbol{ZB}}\begin{bmatrix}
                  \ 0   \\[0.3em]
                  \ 1     \\[0.3em]
                  \ \bar{u} - \frac{\bar{a}}{\sqrt{\gamma(\gamma - 1)}}
       \end{bmatrix}    + 
  \bar{\alpha}_{p,2}^{\boldsymbol{ZB}}\begin{bmatrix}
                  \ 1    \\[0.3em]
                  \ \bar{u}     \\[0.3em]
                  \ \frac{1}{2}\bar{u}^{2} 
          \end{bmatrix}    + 
   \bar{\alpha}_{p,3}^{\boldsymbol{ZB}}\begin{bmatrix}
                    \ 0  \\[0.3em]
                    \ 1 \\[0.3em]
                    \ \bar{u} + \frac{\bar{a}}{\sqrt{\gamma(\gamma - 1)}}
      \end{bmatrix}  
\end{equation}
On comparing first equation we have,
\begin{equation}\label{eq_1_ZBS-FDS wave_strenghts_p_part}
   \\  \bar{\alpha}_{p,2}^{\boldsymbol{ZB}}  \ = \ \Delta{\rho}   
\end{equation}
From second equation, we get  
\begin{equation}
 \Delta(\rho u)   \  =  \   \bar{\alpha}_{p,1}^{\boldsymbol{ZB}} \ + \ \bar{u} \bar{\alpha}_{p,2}^{\boldsymbol{ZB}}  \ + \ \bar{\alpha}_{p,3}^{\boldsymbol{ZB}}
\end{equation}
On using (\ref{eq_1_ZBS-FDS wave_strenghts_p_part}) in the above equation we get expression as
\begin{equation}\label{eq_2_ZBS-FDS wave_strenghts_p_part}
 \\    \bar{\alpha}_{p,1}^{\boldsymbol{ZB}} \ + \ \bar{\alpha}_{p,3}^{\boldsymbol{ZB}}  \ = \ \bar{\rho}\Delta{u}
\end{equation}
Similarly, from third equation 
\begin{equation}
 \Delta(\rho E)  \ = \  \bar{\alpha}_{p,1}^{\boldsymbol{ZB}}\bigg\{\bar{u} - \frac{\bar{a}}{\sqrt{\gamma(\gamma - 1)}}\bigg\} \ + \ \frac{1}{2}\bar{u}^{2}\bar{\alpha}_{p,2}^{\boldsymbol{ZB}} \ + \ \bar{\alpha}_{p,3}^{\boldsymbol{ZB}}\bigg\{\bar{u} + \frac{\bar{a}}{\sqrt{\gamma(\gamma - 1)}}\bigg\}
\end{equation}
or
\begin{align}
 \begin{split}
\Delta\Big(\frac{p}{(\gamma - 1)} \ + \ \frac{1}{2} \rho u^2\Big)  \ &= \  \bar{u}(\bar{\alpha}_{p,1}^{\boldsymbol{ZB}} \ + \ \bar{\alpha}_{p,3}^{\boldsymbol{ZB}}) + \frac{\bar{a}}{\sqrt{\gamma(\gamma - 1)}}(\bar{\alpha}_{p,3}^{\boldsymbol{ZB}} \ - \ \bar{\alpha}_{p,1}^{\boldsymbol{ZB}}) \\ \ &+ \ \frac{1}{2}\bar{u}^{2}\bar{\alpha}_{p,2}^{\boldsymbol{ZB}}
 \end{split}
\end{align}
on using (\ref{eq_2_ZBS-FDS wave_strenghts_p_part}) in the above equation, we get
\begin{equation}
\\  \frac{1}{(\gamma - 1)}\Delta{p} \ + \ \frac{1}{2}\Delta(\rho u^2)  \ = \ \bar{\rho}\bar{u}\Delta{u} \ + \ \frac{\bar{a}}{\sqrt{\gamma(\gamma - 1)}}(\bar{\alpha}_{p,3}^{\boldsymbol{ZB}} \ - \ \bar{\alpha}_{p,1}^{\boldsymbol{ZB}}) \ + \ \frac{1}{2}\bar{u}^{2}\Delta{\rho}
\end{equation}
Now, $\frac{1}{2}\Delta(\rho u^2) \ = \ \bar{\rho}\bar{u}\Delta{u} \ + \ \frac{1}{2}\bar{u}^{2}\Delta{\rho} $ is an equation for above defined averages values of  $\bar{\rho} \ \textrm{and} \ \bar{u}$. Thus we are left with
\begin{equation}\label{eq_3_ZBS-FDS wave_strenghts_p_part}
 \bar{\alpha}_{p,3}^{\boldsymbol{ZB}} \ - \ \bar{\alpha}_{p,1}^{\boldsymbol{ZB}}  \ = \  \sqrt{\frac{\gamma}{\gamma - 1}}\frac{\Delta{p}}{\bar{a}}
\end{equation}  
On solving (\ref{eq_2_ZBS-FDS wave_strenghts_p_part}) and (\ref{eq_3_ZBS-FDS wave_strenghts_p_part}) simultaneously, we get
\begin{align}
 \begin{split}
  \bar{\alpha}_{p,1}^{\boldsymbol{ZB}} \ &= \   \frac{\bar{\rho}\Delta{u}}{2} \ - \ \sqrt{\frac{\gamma}{\gamma - 1}}\frac{\Delta{p}}{2\bar{a}} \ \    \ \textrm{and} \  \\
  \bar{\alpha}_{p,3}^{\boldsymbol{ZB}} \ &= \   \frac{\bar{\rho}\Delta{u}}{2} \ + \ \sqrt{\frac{\gamma}{\gamma - 1}}\frac{\Delta{p}}{2\bar{a}}
\end{split}
\end{align}
Let us consider convective subsystem of Zha and Bilgen type splitting. On comparing (\ref{conservation_ZB_c_part}) and (\ref{quasi_form_ZB_c_part}) and writing in finite difference analogue, we have
\begin{equation}\label{FDS_ZBS_c_part}
 \Delta{\boldsymbol{F}_{c}^{\boldsymbol{ZB}}} \ = \ \boldsymbol{\bar{A}}_{c}^{\boldsymbol{ZB}}\Delta{\boldsymbol{U}}
\end{equation}
where $\boldsymbol{\bar{A}}_{c}^{\boldsymbol{ZB}}$ is now a function of left and right states, {\em i.e.}, $\boldsymbol{\bar{A}}_{c}^{\boldsymbol{ZB}} = \boldsymbol{\bar{A}}_{c}^{\boldsymbol{ZB}}(\boldsymbol{U}_L,\boldsymbol{U}_R)$. It is worth nothing that (\ref{FDS_ZBS_c_part}) is just a relation and may not become an equation for already defined average values. Further  $\Delta{\boldsymbol{U}}$ is a column vector and from the theory of Jordan Forms we are able to get complete set of LI generalized eigenvectors. Thus we can form generalized basis for $\Delta{\boldsymbol{U}}$, {\em i.e.},
\begin{equation}
 \Delta \boldsymbol{U}   \ = \  \sum_{i = 1}^{3} \bar{\alpha}_{c,i}^{\boldsymbol{ZB}}\boldsymbol{\bar{R}}_{c,i}^{\boldsymbol{ZB}} 
\end{equation}
On using above expression in (\ref{FDS_ZBS_c_part}), we get
\begin{equation}
\Delta{\boldsymbol{F}_{c}^{\boldsymbol{ZB}}} \ = \  \boldsymbol{\bar{A}}_{c}^{\boldsymbol{ZB}} \sum_{i = 1}^{3} \bar{\alpha}_{c,i}^{\boldsymbol{ZB}}\boldsymbol{{\bar{R}}}_{c,i}^{\boldsymbol{ZB}}
\end{equation}
or
\begin{equation}
\Delta{\boldsymbol{F}_c^{\boldsymbol{ZB}}}   \ = \   \bar{\alpha}_{c,1}^{\boldsymbol{ZB}} \boldsymbol{\bar{A}}_{c}^{\boldsymbol{ZB}} \boldsymbol{{\bar{R}}}_{c,1}^{\boldsymbol{ZB}}  \ + \ \bar{\alpha}_{c,2}^{\boldsymbol{ZB}} \boldsymbol{\bar{A}}_{c}^{\boldsymbol{ZB}} \boldsymbol{{\bar{R}}}_{c,2}^{\boldsymbol{ZB}}  \ + \ \bar{\alpha}_{c,3}^{\boldsymbol{ZB}} \boldsymbol{\bar{A}}_{c}^{\boldsymbol{ZB}} \boldsymbol{{\bar{R}}}_{c,3}^{\boldsymbol{ZB}}
\end{equation} 
which is further equal to
\begin{equation} 
\Delta{\boldsymbol{F}_c^{\boldsymbol{ZB}}}   \ = \   \bar{\alpha}_{c,1}^{\boldsymbol{ZB}}\bar \lambda_{c,1}^{\boldsymbol{ZB}} \boldsymbol{{\bar{R}}}_{c,1}^{\boldsymbol{ZB}}  \ + \ \bar{\alpha}_{c,2}^{\boldsymbol{ZB}}\big( \bar \lambda_{c,2}^{\boldsymbol{ZB}} \boldsymbol{{\bar{R}}}_{c,2}^{\boldsymbol{ZB}} + \boldsymbol{{\bar{R}}}_{c,1}^{\boldsymbol{ZB}}\big) \ + \ \bar{\alpha}_{c,3}^{\boldsymbol{ZB}} \bar \lambda_{c,3}^{\boldsymbol{ZB}} \boldsymbol{{\bar{R}}}_{c,3}^{\boldsymbol{ZB}}
\end{equation} 
Now, $\Delta{\boldsymbol{F}_p^{+\boldsymbol{ZB}}}$ must have contribution of positive part of eigenvalues only, {\em i.e.},
\begin{align}\label{ZBS-FDS_positive_c_part}
 \begin{split}
  \Delta{\boldsymbol{F}_c^{+\boldsymbol{ZB}}}   \ &= \  \bar{\alpha}_{c,1}^{\boldsymbol{ZB}}\bar \lambda_{c,1}^{+\boldsymbol{ZB}} \boldsymbol{{\bar{R}}}_{c,1}^{\boldsymbol{ZB}}  \ + \ \bar{\alpha}_{c,2}^{\boldsymbol{ZB}} \bar \lambda_{c,2}^{+\boldsymbol{ZB}} \boldsymbol{{\bar{R}}}_{c,2}^{\boldsymbol{ZB}}  \ + \ \bar{\alpha}_{c,3}^{\boldsymbol{ZB}} \bar \lambda_{c,3}^{+\boldsymbol{ZB}} \boldsymbol{{\bar{R}}}_{c,3}^{\boldsymbol{ZB}} \\
  \ &+ \   \bar{\alpha}_{c,2}^{\boldsymbol{ZB}}\boldsymbol{{\bar{R}}}_{c,1}^{\boldsymbol{ZB}}
 \end{split}
\end{align}
Similarly,
\begin{align}\label{ZBS-FDS_negative_c_part}
 \begin{split}
  \Delta{\boldsymbol{F}_c^{-\boldsymbol{ZB}}}   \ &= \  \bar{\alpha}_{c,1}^{\boldsymbol{ZB}}\bar \lambda_{c,1}^{-\boldsymbol{ZB}} \boldsymbol{{\bar{R}}}_{c,1}^{\boldsymbol{ZB}}  \ + \ \bar{\alpha}_{c,2}^{\boldsymbol{ZB}} \bar \lambda_{c,2}^{-\boldsymbol{ZB}} \boldsymbol{{\bar{R}}}_{c,2}^{\boldsymbol{ZB}}  \ + \ \bar{\alpha}_{c,3}^{\boldsymbol{ZB}} \bar \lambda_{c,3}^{-\boldsymbol{ZB}} \boldsymbol{{\bar{R}}}_{c,3}^{\boldsymbol{ZB}} \\
  \ &+ \   \bar{\alpha}_{c,2}^{\boldsymbol{ZB}}\boldsymbol{{\bar{R}}}_{c,1}^{\boldsymbol{ZB}}
 \end{split}
\end{align}
Again, we define the standard Courant splitting for the eigenvalues as 
 \begin{equation}
 \bar \lambda^{\pm}_{c,i} = \frac{\bar \lambda_{c,i} \pm |\bar \lambda_{c,i}|}{2} 
 \end{equation}
 On using standard upwinding along with above definition, we finally get
\begin{equation}
 \Delta{\boldsymbol{F}_c^{+\boldsymbol{ZB}}} - \Delta{\boldsymbol{F}_c^{-\boldsymbol{ZB}}}  \ = \ \sum_{i = 1}^{3} \bar{\alpha}_{c,i}^{\boldsymbol{ZB}}\left|\bar \lambda_{c,i}^{\boldsymbol{ZB}}\right| \boldsymbol{{\bar{R}}}_{c,i}^{\boldsymbol{ZB}}
\end{equation}
As we can see the resultant of extra contribution, which is coming because of weak hyperbolicity of convective subsystem, turns out be equal to zero. Unlike pressure subsystem here we don't need to find wave strengths because all eigenvalues corresponding to convective subsystem are equal, which leads to
\begin{equation}
 \Delta{\boldsymbol{F}_c^{+\boldsymbol{ZB}}} - \Delta{\boldsymbol{F}_c^{-\boldsymbol{ZB}}}  \ = \ \left|\bar \lambda_{c}^{\boldsymbol{ZB}}\right| \sum_{i = 1}^{3} \bar{\alpha}_{c,i}^{\boldsymbol{ZB}} \boldsymbol{{\bar{R}}}_{c,i}^{\boldsymbol{ZB}}
\end{equation}
or
\begin{equation}
 \Delta{\boldsymbol{F}_c^{+\boldsymbol{ZB}}} - \Delta{\boldsymbol{F}_c^{-\boldsymbol{ZB}}}  \ = \ \left|\bar \lambda_{c}^{\boldsymbol{ZB}}\right| \Delta{\boldsymbol{U}}
\end{equation}
Now,
\begin{equation}
\Delta{\boldsymbol{U}} \ = \  \begin{bmatrix}
         \Delta U_1  \\[0.3em]
       \Delta  U_2 \\[0.3em]
       \Delta U_3
        \end{bmatrix}  \ = \ 
\begin{bmatrix}
         \Delta (\rho)  \\[0.3em]
       \Delta  (\rho u) \\[0.3em]
       \Delta  (\rho E)
        \end{bmatrix}
       \end{equation}
where,
\begin{align}
 \begin{split}
\Delta{U_1}  \  &= \  \rho_R - \rho_L   \\
 \Delta{U_2}  \ &= \ \Delta(\rho u) \ = \ \bar{\rho} \Delta{u}  \ + \ \bar{u} \Delta{\rho} \ \textrm{and} \ \\
\Delta{U_3} \ &= \  \Delta(\rho E)  \ = \  \Delta \Big(\frac{p}{(\gamma - 1)} \ + \ \frac{1}{2} \rho u^{2}\Big)
\\
 \ &= \ \frac{1}{\gamma -1} \Delta{p}  \ + \ \frac{1}{2} \big(\bar{u}^2 \Delta{\rho} \ + \ 2 \bar{\rho} \bar{u} \Delta{u})
\end{split}
\end{align}
Here, we did an experiment to check $\Delta{(\rho E)} \ =  \  \frac{1}{\gamma -1} \Delta{p}  \ + \ \frac{1}{2} \big(\bar{u}^2 \Delta{\rho} \ + \ 2 \bar{\rho} \bar{u} \Delta{u})$ holds or not. As we know from the theory of gas dynamics \cite{Zucrow_Gas_dynamics}, ratio of densities, {\em i.e.}, $
(\frac{\rho_r}{\rho_l})$  attains a constant value of $6$ as Mach number $M \rightarrow \infty$.  Next, we define 
\begin{equation}
 \textrm{error}_3  \ = \  \Delta{(\rho E)}  \ - \  \frac{1}{\gamma -1} \Delta{p}  \ - \ \frac{1}{2} \big(\bar{u}^2 \Delta{\rho} \ + \ 2 \bar{\rho} \bar{u} \Delta{u}) 
\end{equation}
and we consider a test case with variable Mach number taken from \cite{steady_shock_problem}, which is given below.  
\begin{equation*}
 \begin{bmatrix}
        p_l \\[0.3em]
         \rho_l \\[0.3em]
        u_l
         \end{bmatrix} \ = \
    \begin{bmatrix}
         \frac{1}{\gamma M^2}  \\[0.3em]
         1.0  \\[0.3em]
         1.0 
        \end{bmatrix}  \ \textrm{and} \
      \begin{bmatrix}
        p_r \\[0.3em]
         \rho_r \\[0.3em]
        u_r
         \end{bmatrix} \ = \
    \begin{bmatrix}
         p_l \frac{2 \gamma M^{2} \ - \ (\gamma - 1)}{\gamma + 1}  \\[0.8em]
         \dfrac{\frac{\gamma + 1}{\gamma - 1} \frac{p_r}{p_l} + 1}{\frac{\gamma + 1}{\gamma  - 1}  +  \frac{p_r}{p_l}}   \\[1.5em]
         \sqrt {\gamma \frac{(2+(\gamma - 1)M^{2})p_r}{(2\gamma M^{2} + (1 - \gamma ))\rho_r}}
        \end{bmatrix}  
\end{equation*}
As per expectations density ratio approaches to limit $6$ as shown in Figure \ref{ZBS-FDS_d_ratio.eps} and error$_3$, which is  given in Figure \ref{ZBS-FDS_error-3.pdf} is not exactly zero and it fluctuates between limits $-10^{-15}$ to $10^{-15}$, which is anyhow very small. The possible reason of this could be the generation of round-off error and if we take macroscopic scale, error looks almost equal to zero as shown in Figure (\ref{ZBS-FDS_error-3_macro.eps}).
 \\
The final expressions in the finite volume framework, with ZBS-FDS scheme for Euler equations are as follows.   
 \begin{equation}
 \boldsymbol{U}^{n+1}_{j} \ = \ \boldsymbol{U}^{n}_{j} - \frac{\Delta t}{\Delta x} 
 \left[ \boldsymbol{F}^{n}_{j+\frac{1}{2}} \ - \ \boldsymbol{F}^{n}_{j-\frac{1}{2}} \right] 
 \end{equation} 
 where the cell-interface fluxes, $\boldsymbol{F}_{I} \ = \ \boldsymbol{F}_{j\pm\frac{1}{2}}$, are defined by 
 \begin{equation}
 \boldsymbol{F}_{I} = \frac{1}{2} \left[\boldsymbol{F}_{L} \ + \ \boldsymbol{F}_{R} \right] 
 - \frac{1}{2} \left[\left(\Delta{\boldsymbol{F}_c^{+\boldsymbol{ZB}}} - \Delta{\boldsymbol{F}_c^{-\boldsymbol{ZB}}}\right) \ + \ \left(\Delta{\boldsymbol{F}_p^{+\boldsymbol{ZB}}} - \Delta{\boldsymbol{F}_p^{-\boldsymbol{ZB}}}\right) \right] 
 \end{equation}  
where 
 \begin{equation}
\Delta{\boldsymbol{F}_c^{+\boldsymbol{ZB}}} - \Delta{\boldsymbol{F}_c^{-\boldsymbol{ZB}}}  \ = \ \left|\bar \lambda_{c}^{\boldsymbol{ZB}}\right|\begin{bmatrix}
         \rho_R - \rho_ L  \\[0.8em]
          \bar{\rho} \Delta{u}  \ + \ \bar{u} \Delta{\rho}  \\[0.8em]
       \frac{1}{\gamma -1} \Delta{p}  \ + \ \frac{1}{2} \big(\bar{u}^2 \Delta{\rho} \ + \ 2 \rho \bar{u} \Delta{u})
        \end{bmatrix}
\end{equation} 
and
\begin{equation}
\Delta{\boldsymbol{F}_p^{+\boldsymbol{ZB}}} - \Delta{\boldsymbol{F}_p^{-\boldsymbol{ZB}}}  \ = \ \sum_{i = 1}^{3} \bar{\alpha}_{p,i}^{\boldsymbol{ZB}}\left|\bar \lambda_{p,i}^{\boldsymbol{ZB}}\right| \boldsymbol{{\bar{R}}}_{p,i}^{\boldsymbol{ZB}}
\end{equation}
\subsection{Formulation of TVS-FDS scheme}
Let us consider pressure subsystem corresponding to Toro and V\'azquez type splitting and on comparing  (\ref{conservation_TV_p_part}) and (\ref{quasi_form_TV_p_part}), and on using the finite difference analogue we have,
\begin{equation}
 \Delta{\boldsymbol{F}_{p}^{\boldsymbol{TV}}} \ = \ \boldsymbol{\bar{A}}_{p}^{\boldsymbol{TV}}\Delta{\boldsymbol{U}}
\end{equation}
Similar procedure like what we did for pressure subsystem of Zha and Bilgen type splitting  is followed here and finally we get 
\begin{equation}\label{TVS-FDS-p-part}
 \Delta{\boldsymbol{F}_p^{+\boldsymbol{TV}}} - \Delta{\boldsymbol{F}_p^{-\boldsymbol{TV}}}  \ = \ \sum_{i = 1}^{3} \bar{\alpha}_{p,i}^{\boldsymbol{TV}}\left|\bar \lambda_{p,i}^{\boldsymbol{TV}}\right| \boldsymbol{{\bar{R}}}_{p,i}^{\boldsymbol{TV}}
\end{equation} 
To determine right side of (\ref{TVS-FDS-p-part}) completely, we need to find average values of eigenvalues along with wave strengths.  Consider the linearization equation of pressure subsystem for TVS-FDS scheme, $\Delta \boldsymbol{F}_{p}^{\boldsymbol{TV}}  \ = \ \boldsymbol{\bar{A}}_{p}^{\boldsymbol{TV}}\Delta{\boldsymbol{U}}$, {\em i.e.},
\begin{equation}
  \begin{bmatrix}
         0 \\[0.3em]
       \Delta  p \\[0.3em]
       (\frac{\gamma}{\gamma-1}) \Delta p u
         \end{bmatrix} \ = \
    \begin{bmatrix}
         \ 0  & 0 & 0  \\[0.3em]
         \frac{1}{2} (\gamma -1) \bar u^2  &  -(\gamma -1)\bar u  &  (\gamma -1) \\[0.3em]
         \ -\frac{\bar u \bar a^2}{(\gamma - 1)} + \frac{1}{2}\gamma \bar u^3  &  \frac{\bar a^2}{(\gamma - 1)} - \gamma \bar u^2  &  \gamma \bar u 
        \end{bmatrix} 
        \begin{bmatrix}
         \Delta (\rho)  \\[0.3em]
       \Delta  (\rho u) \\[0.3em]
       \Delta  (\rho E)
        \end{bmatrix}
       \end{equation}
From the second equation, we get
\begin{equation}
 \Delta p   \  =  \  \frac{1}{2} (\gamma -1) \bar u^2 \Delta \rho  \  -  \ (\gamma -1)\bar u \Delta  (\rho u)  \  +  \  (\gamma -1) \Delta  (\rho E)
\end{equation}
or
\begin{align}
\begin{split}
 \Delta p   \  &=  \  \frac{1}{2} (\gamma -1) \bar u^2 \Delta \rho  \  -  \ (\gamma -1)\bar u \Delta  (\rho u) \\ \  &+  \  (\gamma -1) \Delta(\frac{p}{\gamma - 1}) \ + \ \frac{1}{2} (\gamma - 1) \Delta(\rho u^2) 
\end{split}
\end{align}
which further reduces to 
\begin{equation}
 \   {\bar{u}}^2 \Delta (\rho) \ - \ 2 \bar{u} \Delta  (\rho u) \ + \ \Delta (\rho u^{2})\ = \ 0
\end{equation}
\begin{equation}
 \Rightarrow  \\  \bar{u}    \ = \  \frac{\sqrt{\rho_L} u_L  \ + \ \sqrt{\rho_R} u_R }{\sqrt{\rho_L} \ + \ \sqrt{\rho_R}}
\end{equation}
Average density can be found by substituting $\bar{u}$ in the relation
\begin{equation}
 \rho_R u_R - \rho_L u_L   \ = \  \bar{\rho}(u_R - u_L) \ + \  \bar{u}(\rho_R - \rho_L)
\end{equation}
and after some simple calculation, we get $\bar{\rho} = \sqrt{\rho_L\rho_R}$. 
Similarly, third equation can be written as 
\begin{align}
 \begin{split}
\left(\frac{\gamma}{\gamma-1}\right) \Delta( p u ) \ &= \  \left(-\frac{\bar u \bar a^2}{(\gamma - 1)}  \ + \ \frac{1}{2}\gamma \bar u^3\right) \Delta \rho \\ \ &+ \   \left(\frac{\bar a^2}{(\gamma - 1)} \  -  \ \gamma \bar u^2\right) \Delta (\rho u)  \ + \ \gamma \bar u \Delta( \rho E)
 \end{split}
\end{align}
On expanding we have,
\begin{align}
 \begin{split}
\left(\frac{\gamma}{\gamma-1}\right) \Delta( p u ) \ &= \  \left(-\frac{\bar u \bar a^2}{(\gamma - 1)}  \ + \ \frac{1}{2}\gamma \bar u^3\right) \Delta \rho  \ + \   \left(\frac{\bar a^2}{(\gamma - 1)} \  -  \ \gamma \bar u^2\right) \Delta (\rho u) \\  \ &+ \ \gamma \bar u \Delta\Big(\frac{p}{\gamma - 1} + \frac{1}{2} \rho u^{2} \Big)
 \end{split}
\end{align}
After further simplifications, above relation reduces to
\begin{equation}\label{recall_1}
\gamma \Delta(pu)   \  =  \    \bar{a}^2 \bar{\rho}\Delta{u} \ + \ \gamma \bar{u}\Delta{p}   
\end{equation}
On using $p =  \dfrac{a^2\rho}{\gamma}$, we have
\begin{equation}
\Delta(a^2 \rho u)   \  =  \  \bar{u}\Delta(a^2\rho) \  +  \  \bar{a}^2\bar{\rho} \Delta{u}
\end{equation}
as explained earlier, $\bar{a}^2$ comes out as
\begin{equation}
\bar{a}^2  \ = \ \frac{\sqrt{\rho_L} a^{2}_L  \ + \ \sqrt{\rho_R} a^{2}_R }{\sqrt{\rho_L} \ + \ \sqrt{\rho_R}}
\end{equation}
Similarly, wave strengths for pressure subsystem of TVS-FDS scheme can be calculated from the relation
\begin{equation}
 \Delta \boldsymbol{U}   \ = \  \sum_{i = 1}^{3} \bar{\alpha}_{p,i}^{\boldsymbol{TV}}\boldsymbol{{\bar{R}}}_{p,i}^{\boldsymbol{TV}}
\end{equation}
{\em i.e.}, 
\begin{equation} 
  \begin{bmatrix} 
                 \Delta(\rho) \\[0.3em]  
		\Delta(\rho u) \\[0.3em] 
		\Delta(\rho E) 
		\end{bmatrix}   =  
    \bar{\alpha}_{p,1}^{\boldsymbol{TV}}\begin{bmatrix}
                  \ 0    \\[0.3em]
                  \ 1     \\[0.3em]
                  \ \bar u + \frac{1}{2}(\frac{\bar u - \bar{\beta}}{\gamma - 1}) 
       \end{bmatrix}  + 
  \bar{\alpha}_{p,2}^{\boldsymbol{TV}} \begin{bmatrix}
                  \ 1    \\[0.3em]
                  \ \bar u     \\[0.3em]
                  \ \frac{1}{2} \bar u^2 
          \end{bmatrix}   + 
   \bar{\alpha}_{p,3}^{\boldsymbol{TV}} \begin{bmatrix}
                    \ 0    \\[0.3em]
                  \ 1     \\[0.3em]
                  \ \bar u + \frac{1}{2}(\frac{\bar u + \bar {\beta}}{\gamma - 1}) 
      \end{bmatrix}  
\end{equation}
From the first equation, we get 
\begin{equation}
 \bar{\alpha}_{p,2}^{\boldsymbol{TV}} \ = \  \Delta \rho
\end{equation}
Similarly, the second equation gives 
\begin{equation} 
 \Delta(\rho u)  \ = \  \bar{\alpha}_{p,1}^{\boldsymbol{TV}}   \ + \ \bar{u} \bar{\alpha}_{p,1}^{\boldsymbol{TV}} \ + \ \bar{\alpha}_{p,3}^{\boldsymbol{TV}} 
\end{equation}
\begin{equation} \label{momentum_p}
  \Rightarrow \\ \bar{\alpha}_{p,1}^{\boldsymbol{TV}}  \ + \ \bar{\alpha}_{p,3}^{\boldsymbol{TV}} \ = \   \bar{\rho}\Delta{u} 
\end{equation}
Third equation implies 
\begin{equation} 
\Delta (\rho E)  \  =  \   \{ u + \frac{1}{2}(\frac{u - \bar{\beta}}{\gamma - 1})\}\bar{\alpha}_{p,1}^{\boldsymbol{TV}} \ + \   \frac{1}{2}u^2\bar{\alpha}_{p,2}^{\boldsymbol{TV}} \ + \  \{ u + \frac{1}{2}(\frac{u + \bar{\beta}}{\gamma - 1})\}\bar{\alpha}_{p,3}^{\boldsymbol{TV}}
\end{equation}
On rearrangement of terms and after some algebra, the above equation reduces to 
\begin{equation} \label{internal_p}
 \bar{\alpha}_{p,3}^{\boldsymbol{TV}}  \ - \ \bar{\alpha}_{p,1}^{\boldsymbol{TV}}   \ = \ \frac{1}{\bar{\beta}} \{ 2 \Delta p  \ - \ \bar u \bar \rho \Delta u \}
\end{equation}
Finally on comparing (\ref{momentum_p})\ and\  (\ref{internal_p}), we get both $\bar{\alpha}_{p,1}^{TV}$ and $\bar{\alpha}_{p,3}^{TV}$ as
 \begin{eqnarray}
  \bar{\alpha}_{p,1}^{\boldsymbol{TV}}  \ = \ \frac{1}{2} \bar{\rho} \Delta{u} \ + \ \frac{1}{2 \bar{\beta}} \bar{\rho} \bar{u} \Delta u  \ - \ \frac{\Delta p}{\bar \beta}   \\  
  \bar{\alpha}_{p,3}^{\boldsymbol{TV}} \ = \ \frac{1}{2} \bar \rho \Delta u    \ - \ \frac{1}{2 \bar \beta} \bar \rho \bar u \Delta u  \ + \ \frac{\Delta p}{\bar \beta} 
 \end{eqnarray}  
 Therefore, the wave strengths are finally given by 
 \begin{align}
 \begin{split}
  \bar{\alpha}_{p,1}^{\boldsymbol{TV}} \ &= \ \frac{1}{2} \bar \rho \Delta u    \ + \ \frac{1}{2 \bar \beta} \bar \rho \bar u \Delta u  \ - \frac{\Delta p}{\bar \beta} \ , \ 
\bar{\alpha}_{p,2}^{\boldsymbol{TV}} \ = \  \Delta \rho \\   \ \ \textrm{and} \ \
 \bar{\alpha}_{p,3}^{\boldsymbol{TV}}  \ &= \ \frac{1}{2} \bar \rho \Delta u    \ - \ \frac{1}{2 \bar \beta} \bar \rho \bar u \Delta u  \ + \ \frac{\Delta p}{\bar \beta} 
\end{split}
\end{align}
We know that the convective subsystem is weakly hyperbolic but we can still form a basis of generalized eigenvectors, {\em i.e.},
 \begin{equation}
 \Delta \boldsymbol{U}   \ = \  \sum_{i = 1}^{3} \bar{\alpha}_{c,i}^{\boldsymbol{TV}}\boldsymbol{{\bar{R}}}_{c,i}^{\boldsymbol{TV}} 
\end{equation} and on comparing (\ref{conservation_TV_c_part}) and (\ref{quasi_form_TV_c_part}), and after writing in a finite difference form, we have the relation
\begin{equation}
 \\    \Delta{\boldsymbol{F}_{c}^{\boldsymbol{TV}}} \ = \  \boldsymbol{\bar{A}}_{c}^{\boldsymbol{TV}} \sum_{i = 1}^{3} \bar{\alpha}_{c,i}^{\boldsymbol{TV}}\boldsymbol{{\bar{R}}}_{c,i}^{\boldsymbol{TV}}
\end{equation}
 or
 \begin{equation} \label{full_convection_flux_for_C_part}
\Delta{\boldsymbol{F}_c^{\boldsymbol{TV}}}   \ = \   \bar{\alpha}_{c,1}^{\boldsymbol{TV}} \boldsymbol{\bar{A}}_{c}^{\boldsymbol{TV}} \boldsymbol{{\bar{R}}}_{c,1}^{\boldsymbol{TV}}  \ + \ \bar{\alpha}_{c,2}^{\boldsymbol{TV}} \boldsymbol{\bar{A}}_{c}^{\boldsymbol{TV}} \boldsymbol{{\bar{R}}}_{c,2}^{\boldsymbol{TV}}  \ + \ \bar{\alpha}_{c,3}^{\boldsymbol{TV}} \boldsymbol{\bar{A}}_{c}^{\boldsymbol{TV}} \boldsymbol{{\bar{R}}}_{c,3}^{\boldsymbol{TV}}
\end{equation} 
Since $\bar A$ is non-diagonalizable, this means 
\begin{equation*}
 \boldsymbol{\bar{A}}_{c}^{\boldsymbol{TV}} \boldsymbol{{\bar{R}}}_{c,i}^{\boldsymbol{TV}}  \ \neq \  \bar \lambda_{c,i}^{\boldsymbol{TV}} \boldsymbol{{\bar{R}}}_{c,i}^{\boldsymbol{TV}} 
\end{equation*} for some i's.
We know $\boldsymbol{{\bar{R}}}_{c,3}^{\boldsymbol{TV}}$  is a generalized eigenvector and corresponding to $\boldsymbol{{\bar{R}}}_{c,2}^{\boldsymbol{TV}}$, a Jordan chain of order two will be formed, {\em i.e.}, 
\begin{equation}
 \boldsymbol{\bar{A}}_{c}^{\boldsymbol{TV}} \boldsymbol{{\bar{R}}}_{c,2}^{\boldsymbol{TV}}  \ = \  \bar \lambda_{c,2}^{\boldsymbol{TV}} \boldsymbol{{\bar{R}}}_{c,2}^{\boldsymbol{TV}}\
 ~~\mbox{and}~~
 \boldsymbol{\bar{A}}_{c}^{\boldsymbol{TV}} \boldsymbol{{\bar{R}}}_{c,3}^{\boldsymbol{TV}}  \ = \  \bar \lambda_{c,3}^{\boldsymbol{TV}} \boldsymbol{{\bar{R}}}_{c,3}^{\boldsymbol{TV}} \ + \ \boldsymbol{{\bar{R}}}_{c,2}^{\boldsymbol{TV}}
\end{equation}
On using above relations in (\ref{full_convection_flux_for_C_part}), we get
\begin{equation} 
\Delta{\boldsymbol{F}_c^{\boldsymbol{TV}}}   \ = \   \bar{\alpha}_{c,1}^{\boldsymbol{TV}}\bar \lambda_{c,1}^{\boldsymbol{TV}} \boldsymbol{{\bar{R}}}_{c,1}^{\boldsymbol{TV}}  \ + \ \bar{\alpha}_{c,2}^{\boldsymbol{TV}} \bar \lambda_{c,2}^{\boldsymbol{TV}} \boldsymbol{{\bar{R}}}_{c,2}^{\boldsymbol{TV}}  \ + \ \bar{\alpha}_{c,3}^{\boldsymbol{TV}} \bar \lambda_{c,3}^{\boldsymbol{TV}} \boldsymbol{{\bar{R}}}_{c,3}^{\boldsymbol{TV}}  \ + \ \bar \alpha_{c,3}^{\boldsymbol{TV}}\boldsymbol{{\bar{R}}}_{c,2}^{\boldsymbol{TV}}
\end{equation} 
After using standard Courant splitting for the eigenvalues,
$\Delta{\boldsymbol{F}_c^{+\boldsymbol{TV}}}$ and $\Delta{\boldsymbol{F}_c^{-\boldsymbol{TV}}}$ are given by 
\begin{align}\label{flux_c_plus}
 \begin{split}
  \Delta{\boldsymbol{F}_c^{+\boldsymbol{TV}}}   \ = \  \bar{\alpha}_{c,1}^{\boldsymbol{TV}}\bar \lambda_{c,1}^{+\boldsymbol{TV}} \boldsymbol{{\bar{R}}}_{c,1}^{\boldsymbol{TV}}  \ &+ \ \bar{\alpha}_{c,2}^{\boldsymbol{TV}} \bar \lambda_{c,2}^{+\boldsymbol{TV}} \boldsymbol{{\bar{R}}}_{c,2}^{\boldsymbol{TV}}  \ + \ \bar{\alpha}_{c,3}^{\boldsymbol{TV}} \bar \lambda_{c,3}^{+\boldsymbol{TV}} \boldsymbol{{\bar{R}}}_{c,3}^{\boldsymbol{TV}} \\
\ &+ \ \bar \alpha_{c,3}^{\boldsymbol{TV}}\boldsymbol{{\bar{R}}}_{c,2}^{\boldsymbol{TV}}
 \end{split}
\end{align}
and
\begin{align}\label{flux_c_minus}
 \begin{split}
\Delta{\boldsymbol{F}_c^{-\boldsymbol{TV}}}   \ = \  \bar{\alpha}_{c,1}^{\boldsymbol{TV}}\bar \lambda_{c,1}^{-\boldsymbol{TV}} \boldsymbol{{\bar{R}}}_{c,1}^{\boldsymbol{TV}}  \ &+ \ \bar{\alpha}_{c,2}^{\boldsymbol{TV}} \bar \lambda_{c,2}^{-\boldsymbol{TV}} \boldsymbol{{\bar{R}}}_{c,2}^{\boldsymbol{TV}}  \ + \ \bar{\alpha}_{c,3}^{\boldsymbol{TV}} \bar \lambda_{c,3}^{-\boldsymbol{TV}} \boldsymbol{{\bar{R}}}_{c,3}^{\boldsymbol{TV}} \\
\ &+ \ \bar \alpha_{c,3}^{\boldsymbol{TV}}\boldsymbol{{\bar{R}}}_{c,2}^{\boldsymbol{TV}}
 \end{split}
\end{align}
Therefore, we have
\begin{align}
\begin{split}
\Delta{\boldsymbol{F}_c^{+\boldsymbol{TV}}} - \Delta{\boldsymbol{F}_c^{-\boldsymbol{TV}}}  \ &= \  \sum_{i = 1}^{3} \bar{\alpha}_{c,i}^{\boldsymbol{TV}}\bar \lambda_{c,i}^{+\boldsymbol{TV}} \boldsymbol{{\bar{R}}}_{c,i}^{\boldsymbol{TV}}  \ + \ \bar \alpha_{c,3}^{\boldsymbol{TV}}\boldsymbol{{\bar{R}}}_{c,2}^{\boldsymbol{TV}} \\
\ &- \ \sum_{i = 1}^{3} \bar{\alpha}_{c,i}^{\boldsymbol{TV}}\bar \lambda_{c,i}^{-\boldsymbol{TV}} \boldsymbol{{\bar{R}}}_{c,i}^{\boldsymbol{TV}}  \ - \ \bar \alpha_{c,3}^{\boldsymbol{TV}}\boldsymbol{{\bar{R}}}_{c,2}^{\boldsymbol{TV}}
\end{split}
\end{align}   
\begin{equation}
\\ \Rightarrow  \Delta{\boldsymbol{F}_c^{+\boldsymbol{TV}}} - \Delta{\boldsymbol{F}_c^{-\boldsymbol{TV}}}  \ = \ \sum_{i = 1}^{3} \bar{\alpha}_{c,i}^{\boldsymbol{TV}}\left(\bar \lambda_{c,i}^{+\boldsymbol{TV}} \ - \ \bar \lambda_{c,i}^{-\boldsymbol{TV}}\right) \boldsymbol{{\bar{R}}}_{c,i}^{\boldsymbol{TV}}
\end{equation} 
or
\begin{equation}
 \Delta{\boldsymbol{F}_c^{+\boldsymbol{TV}}} - \Delta{\boldsymbol{F}_c^{-\boldsymbol{TV}}}  \ = \ \sum_{i = 1}^{3} \bar{\alpha}_{c,i}^{\boldsymbol{TV}}\left|\bar \lambda_{c,i}^{\boldsymbol{TV}}\right| \boldsymbol{{\bar{R}}}_{c,i}^{\boldsymbol{TV}}
\end{equation} 
Like previous case, resultant of extra contribution comes out equal to zero.
Above relation can be further written as,
\begin{align}
 \begin{split}
  \Delta{\boldsymbol{F}_c^{+\boldsymbol{TV}}} - \Delta{\boldsymbol{F}_c^{-\boldsymbol{TV}}}  \ &= \  \left|\bar \lambda_{c}^{\boldsymbol{TV}}\right|\big(\bar{\alpha}_{c,2}^{\boldsymbol{TV}} \boldsymbol{{\bar{R}}}_{c,2}^{\boldsymbol{TV}} \ + \ \bar{\alpha}_{c,3}^{\boldsymbol{TV}}\boldsymbol{{\bar{R}}}_{c,3}^{\boldsymbol{TV}}\big)   \\
 \textrm{where}
 \left|\bar \lambda_{c}^{\boldsymbol{TV}}\right| \  &= \  \left|\bar \lambda_{c,2}^{\boldsymbol{TV}}\right| = \left|\bar \lambda_{c,3}^{\boldsymbol{TV}}\right|. 
 \end{split}
\end{align}
 This can be further written as
 \begin{equation}\label{TVS-FDS_c_part}
  \Delta{\boldsymbol{F}_c^{+\boldsymbol{TV}}} - \Delta{\boldsymbol{F}_c^{-\boldsymbol{TV}}}  \ = \  \left|\bar \lambda_{c}^{\boldsymbol{TV}}\right|\big[\Delta \boldsymbol{U} \ - \ \bar{\alpha}_{c,1}^{\boldsymbol{TV}}\boldsymbol{{\bar{R}}}_{c,1}^{\boldsymbol{TV}}\big]
 \end{equation}
In order to determine (\ref{TVS-FDS_c_part}) fully, we need to find wave strengths which can be calculated from the relation
\begin{equation}
 \Delta \boldsymbol{U}   \ = \  \sum_{i = 1}^{3} \bar{\alpha}_{c,i}^{\boldsymbol{TV}}\boldsymbol{{\bar{R}}}_{c,i}^{\boldsymbol{TV}}
\end{equation}
or
\begin{equation} 
  \begin{bmatrix} 
                 \Delta(\rho) \\[0.3em]  
		\Delta(\rho u) \\[0.3em] 
		\Delta(\rho E) 
		\end{bmatrix}  \  = \  
   \bar{\alpha}_{c,1}^{\boldsymbol{TV}} \begin{bmatrix}
                  \ 0    \\[0.3em]
                  \ 0     \\[0.3em]
                  \ 1 
       \end{bmatrix}   \ + \ 
  \bar{\alpha}_{c,2}^{\boldsymbol{TV}} \begin{bmatrix}
                  \ 1    \\[0.3em]
                  \ \bar{u}     \\[0.3em]
                  \ \frac{1}{2} \bar{u}^{2} 
          \end{bmatrix}   \ + \ 
    \bar{\alpha}_{c,3}^{\boldsymbol{TV}}\begin{bmatrix}
                    \  x_{1} \\[0.3em]
                    \ 1 + \bar{u}x_{1} \\[0.3em]
                    \ \bar u + \frac{1}{2}\bar{u}^{2}x_{1}
      \end{bmatrix}
\end{equation}
From the first equation, we get
\begin{equation}\label{eq_1_TVFS_wave_strenghts}
\Delta(\rho) \ = \ \bar{\alpha}_{c,2}^{\boldsymbol{TV}} \ + \ x_{1}\bar{\alpha}_{c,3}^{\boldsymbol{TV}}
\end{equation}
From second equation, we have 
\begin{equation}
 \Delta(\rho u)   \  =  \   \bar{u}\bar{\alpha}_{c,2}^{\boldsymbol{TV}}  \  +  \  (1+\bar{u}x_{1})\bar{\alpha}_{c,3}^{\boldsymbol{TV}}
\end{equation}
\begin{equation}
 \\     \Rightarrow   \bar{\rho} \Delta u    \ + \  \bar{u} \Delta \rho    \ = \  \bar{u}\big(\bar{\alpha}_{c,2}^{\boldsymbol{TV}} \ + \ x_{1}\bar{\alpha}_{c,3}^{\boldsymbol{TV}}\big)   \  +  \  \bar{\alpha}_{c,3}^{\boldsymbol{TV}}
\end{equation}
On using (\ref{eq_1_TVFS_wave_strenghts}) in the above equation and after cancellation of some terms, we get
\begin{equation}\label{eq_2_TVFS_wave_strenghts}
\\    \bar{\alpha}_{c,3}^{\boldsymbol{TV}}  \ = \   \bar{\rho}\Delta{u}
\end{equation}
On using (\ref{eq_2_TVFS_wave_strenghts}) in (\ref{eq_1_TVFS_wave_strenghts}) we get 
\begin{equation}
 \\    \bar{\alpha}_{c,2}^{\boldsymbol{TV}} \ = \ \Delta{\rho} \ - \ x_{1}\bar{\rho}\Delta{u}
\end{equation}
Similarly, from third equation we have 
\begin{equation}
 \Delta(\rho E)  \ = \  \bar{\alpha}_{c,1}^{\boldsymbol{TV}}   \ + \  \frac{1}{2}{\bar{u}}^2 \bar{\alpha}_{c,2}^{\boldsymbol{TV}}  \ + \  \big(\bar{u} + x_{1}\frac{1}{2}{\bar{u}}^2\big)\bar{\alpha}_{c,3}^{\boldsymbol{TV}} 
\end{equation}
\begin{equation}
 \Rightarrow \Delta(\frac{p}{\gamma - 1} \ + \ \frac{1}{2} \rho u^2)  \ = \ \bar{\alpha}_{c,1}^{\boldsymbol{TV}}   \ + \  \frac{1}{2}{\bar{u}}^2 \Big(\Delta{\rho} \ - \ x_{1}\bar{\rho}\Delta{u}\Big) \ + \  \big(\bar{u} + x_{1}\frac{1}{2}{\bar{u}}^2\big) \bar{\rho}\Delta{u}
\end{equation}
After little algebra we get,
\begin{equation}
 \\  \bar{\alpha}_{c,1}^{\boldsymbol{TV}}   \  =  \   \frac{1}{(\gamma - 1)}\Delta{p}
\end{equation}  
Therefore, the wave strengths are finally given by 
\begin{equation} 
\bar{\alpha}_{c,1}^{\boldsymbol{TV}}   \ = \  \frac{1}{(\gamma - 1)}\Delta{p} \ , \ \bar{\alpha}_{c,2}^{\boldsymbol{TV}} \ = \  \Delta{\rho} \ - \ x_{1}\bar{\rho}\Delta{u}   \ \ \textrm{and} \ \ \bar{\alpha}_{c,3}^{\boldsymbol{TV}} \ = \   \bar{\rho}\Delta{u}
\end{equation}
On using above calculated wave strength $\bar{\alpha}_{c,1}^{\boldsymbol{TV}}$ along with eigenvector $\boldsymbol{{\bar{R}}}_{c,1}^{\boldsymbol{TV}} \ = \ \boldsymbol{e}_3$  in (\ref{TVS-FDS_c_part}), we get
\begin{equation}
\Delta{\boldsymbol{F}_c^{+\boldsymbol{TV}}} - \Delta{\boldsymbol{F}_c^{-\boldsymbol{TV}}}  \ = \ \left|\bar \lambda_{c}^{\boldsymbol{TV}}\right|\begin{bmatrix}
         \rho_R - \rho_ L  \\[0.8em]
          \bar{\rho} \Delta{u}  \ + \ \bar{u} \Delta{\rho}  \\[0.9em]
        \frac{1}{2} \big(\bar{u}^2 \Delta{\rho} \ + \ 2 \bar{\rho} \bar{u} \Delta{u})
        \end{bmatrix}
\end{equation} 
Thus, the final expressions for TVS-FDS scheme is written as,
 \begin{equation}
 \boldsymbol{U}^{n+1}_{j} \ = \ \boldsymbol{U}^{n}_{j} - \frac{\Delta t}{\Delta x} 
 \left[ \boldsymbol{F}^{n}_{j+\frac{1}{2}} \ - \ \boldsymbol{F}^{n}_{j-\frac{1}{2}} \right] 
 \end{equation} 
 where the cell-interface fluxes, $\boldsymbol{F}_{I} \ = \ \boldsymbol{F}_{j\pm\frac{1}{2}}$, are defined by 
 \begin{align}
  \begin{split}
\boldsymbol{F}_{I} \ = \ \frac{1}{2} \left[\boldsymbol{F}_{L} \ + \ \boldsymbol{F}_{R} \right]  
 \ - \ \frac{1}{2} \left[ \left(\Delta{\boldsymbol{F}_c^{+\boldsymbol{TV}}} - \Delta{\boldsymbol{F}_c^{-\boldsymbol{TV}}}\right) + \left(\Delta{\boldsymbol{F}_p^{+\boldsymbol{TV}}} - \Delta{\boldsymbol{F}_p^{-\boldsymbol{TV}}}\right) \right] 
  \end{split}
\end{align}  
\section{Results and discussion}
We first consider here smooth solution problem with periodic boundary conditions to check experimental order of convergence for both constructed upwind schemes. After that, both ZBS-FDS and TVS-FDS schemes are tested on various one-dimensional Riemann problems of gas dynamics. For most of numerical examples, computational domain  lies between $0$ and $1$, {\em i.e.}, $0\leq x \leq 1.0$ except for Sod's shock tube and shock-entropy test problems. For each problem computational domain is divided into 100 equally spaced cells, except for the shock-entropy interaction test case and for the blast problem test case in which computational domain is divided into $800$ and $3000$ equally spaced cells respectively.
\subsection{Experimental order of convergence (EOC)}
Order of accuracy of both constructed upwind schemes can be determined using the EOC analysis, as given below.
\begin{equation}\label{EOC_formula_1}
     E \ = \ C(\Delta{x})^{s}
\end{equation}
where, $E$ is an error between the exact solution and the numerical solution using some appropriate norm. In particular, we are taking three norms namely, $L_{1}, L_{2}$ and $L_{\infty}$.  Here, $C$ is a constant, $\Delta{x}$ is grid spacing and $s$ is the order of accuracy, which need to be calculated.  On taking logarithms on both sides of (\ref{EOC_formula_1}), we get
\begin{equation}\label{EOC_formula_2}
 log \ {E} \ = \ log \ {C} \ + \ s \ log \ {\Delta{x}}
\end{equation}
which is the equation of a straight line with slope s.  For a given norm, let us initially take $\Delta{x} \ = \ h_{1}$ and on using this in (\ref{EOC_formula_2}), we get 
\begin{equation}\label{EOC_formula_3}
 log \ {E}_{norm,h_{1}} \ = \ log \ {C} \ + \ s \ log \ h_1
\end{equation}
Next we take $\Delta{x} \ = \ h_{2}$, preferably $h_{2} \ = \ \frac{h_{1}}{2}$, with same norm in (\ref{EOC_formula_3}), we get 
\begin{equation}\label{EOC_formula_4}
 log \ {E}_{norm,h_{2}} \ = \ log \ {C} \ + \ s \ log \ h_2
\end{equation}
and on subtracting (\ref{EOC_formula_4}) from (\ref{EOC_formula_3}), formula for finding experimental order of convergence ``s'' comes out as follows.
\begin{equation}
 s \ = \   \frac{\left(log \ {E}_{norm,h_{1}} \ - \ log \ {E}_{norm,h_{2}}\right)}{\left(log \ h_1 \ - \ log \ h_2\right)}
\end{equation}
To check the performance of both schemes in term of errors associated with each norm for different grid sizes, we choose a test case from \cite{Arun_et_al} with initial smooth conditions, for one dimensional Euler system, which are given below.   
\begin{equation*}
 \rho{(x,t)} \ = \ 1.0 + 0.2 \sin(\pi(x-ut)), \  \  u{(x,t)} \ = \ 0.1,  \   \   p{(x,t)} \ = \ 0.5.
\end{equation*}
 For the present case, final solutions remain smooth and periodic boundary conditions are being employed to get meaningful solutions.  Computational domain is chosen as $[0, 2]$, {\em i.e.}, $0\leq x \leq 2$ and all solutions are obtained at final time $t = 0.5$.  $L_{1}$ error norm for both schemes are given in Table \ref{norm-1_smooth_prob} and as per expectations, there is reduction in error on refinement of mess size.  Similarly, $L_{2}$ error norm and $L_{\infty}$ error norm are given in Table \ref{norm-2_smooth_prob} and Table \ref{norm-3_smooth_prob}.  It is clear from all three tables that performance of both schemes is similar, if no discontinuity is present in the solution. 
\begin{table}[!ht]
\caption{$L_{1}$ error norm for smooth solution problem corresponding to both schemes}\label{norm-1_smooth_prob}
\centerline{%
\begin{tabular}{|c|c|c|c|c|c|}
\hline
\textrm{grid points}{}~~ & ZBS-FDS scheme & EOC & TVS-FDS scheme & EOC \\
\hline
$40$ & 0.004476  & - & 0.004476 & - \\
\hline
$80$ & 0.002529  & 0.82 & 0.002529 & 0.82 \\
\hline
$160$ & 0.001258  &  1.007 & 0.001258 & 1.007 \\
\hline
$320$ & 0.000624  & 1.011 & 0.000624 & 1.011 \\
\hline
$640$ & 0.000308  & 1.018 & 0.000308 & 1.018\\
\hline
\end{tabular}
}%
\end{table}
\begin{table}[!ht]
\caption{$L_{2}$ error norm for smooth solution problem corresponding to both schemes}\label{norm-2_smooth_prob}
\centerline{%
\begin{tabular}{|c|c|c|c|c|c|}
\hline
\textrm{grid points}{}~~ & ZBS-FDS scheme & EOC & TVS-FDS scheme & EOC \\
\hline
$40$ & 0.005238  & - & 0.005238 & - \\
\hline
$80$ & 0.003227  & 0.6988 & 0.003227 & 0.6988 \\
\hline
$160$ & 0.001707  &  0.9187 & 0.001707 & 0.9187 \\
\hline
$320$ & 0.000885  & 0.9477 & 0.000885 & 0.9477 \\
\hline
$640$ & 0.000452  & 0.9693 & 0.000452 & 0.9693\\
\hline
\end{tabular}
}%
\end{table}
\begin{table}[!ht]
\caption{$L_{\infty}$ error norm for smooth solution problem corresponding to both schemes}\label{norm-3_smooth_prob}
\centerline{%
\begin{tabular}{|c|c|c|c|c|c|}
\hline
\textrm{grid points}{}~~ & ZBS-FDS scheme & EOC & TVS-FDS scheme & EOC \\
\hline
$40$ & 0.019067  & - & 0.019067 & - \\
\hline
$80$ & 0.013968  & 0.44 & 0.013968 & 0.44 \\
\hline
$160$ & 0.007847  &  0.83 & 0.007847 & 0.83 \\
\hline
$320$ & 0.003994  & 0.974 & 0.003994 & 0.974 \\
\hline
$640$ & 0.002006  & 0.9935 & 0.002006 & 0.9935\\
\hline
\end{tabular}
}%
\end{table}
\subsection{Sod's shock tube problem and Lax problem}
\begin{figure}[!ht]
\centerline{%
\subfigure[]{%
\includegraphics[trim=0 5 35 5, clip, width=0.55\textwidth]{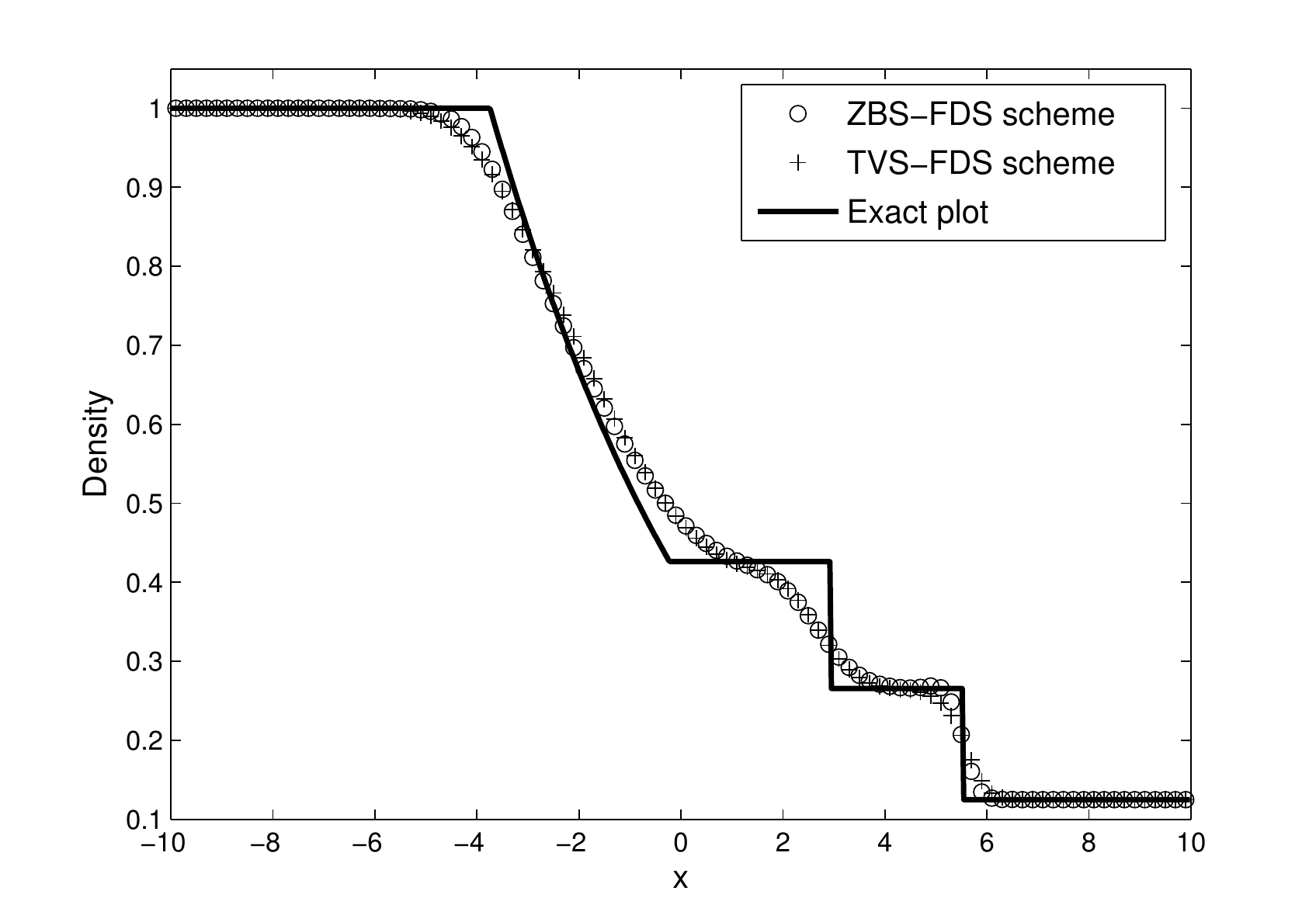}%
\label{Laney1}
}%
\subfigure[]{%
\includegraphics[trim=0 5 35 5, clip, width=0.55\textwidth]{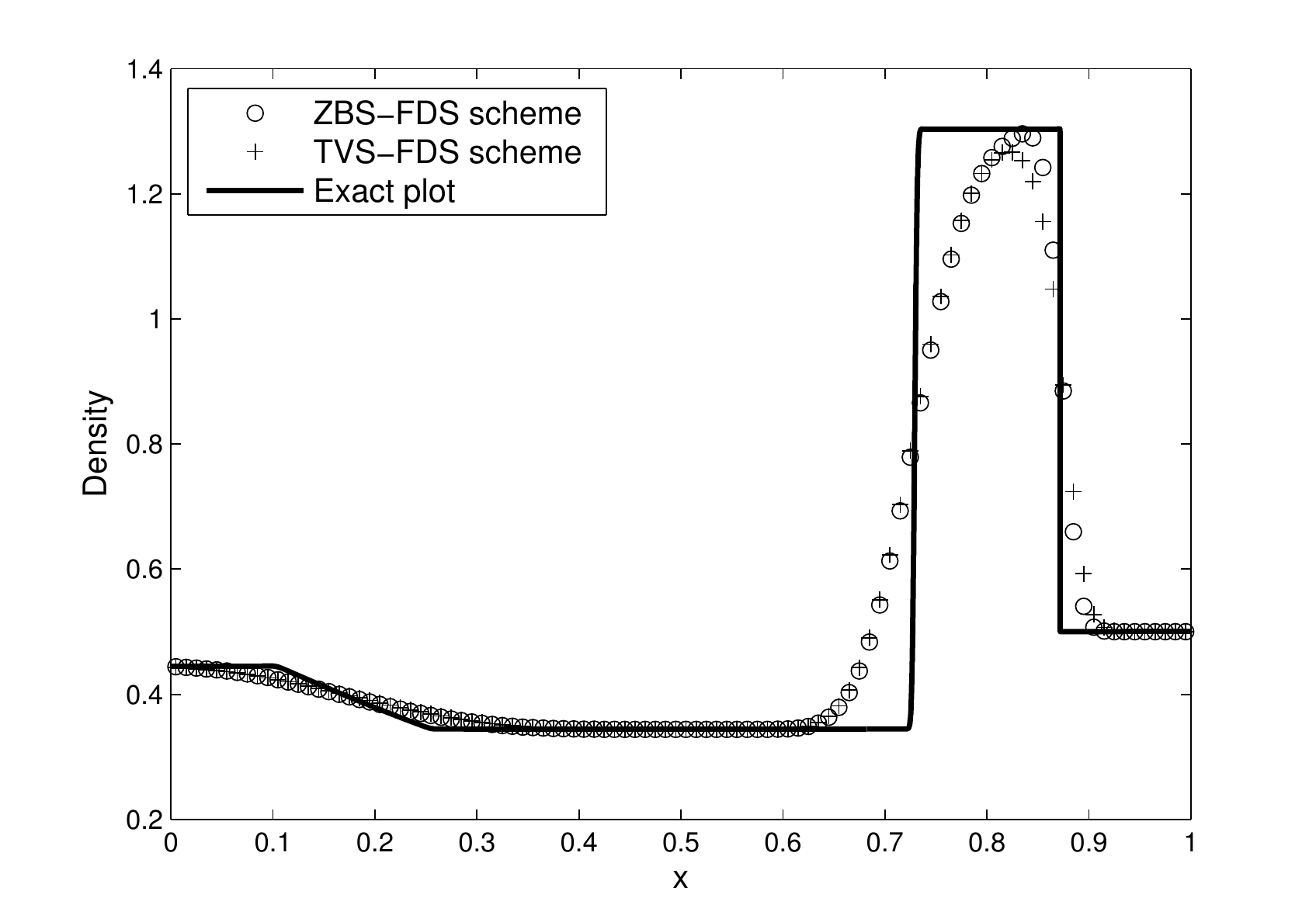}%
\label{wes2_d}
}%
}
\caption{(a) represents results of density variable for Sod's shock tube problem  and (b) represents density plots for Lax problem.}
\end{figure}

First we consider Sod's shock tube problem in which, an initial discontinuity in the middle evolves to, going from right to left as we observe, a shock, a contact discontinuity and an expansion fan.  The initial conditions \cite{Laney_CFD} are $(\rho_L, u_L, p_L)  = (1.0, 0.0, 100000.0)$, $(\rho_R, u_R, p_R) = (0.125, 0.0, 10000.0)$ with initial discontinuity at $x_{o} = 0.0$ and computational domain lies between $-10$ to $10$. All numerical results are obtained at final time $t$ = $0.01$. For this test case, both ZBS-FDS scheme and TVS-FDS scheme exhibit almost similar results except near normal shock region, where ZBS-FDS scheme performs slightly better than TVS-FDS scheme. Results corresponding to density variable are presented in Figure \ref{Laney1}.  We also present error analysis of both schemes corresponding to $L_{1}$-norm and $L_{2}$-norm, and results are given in Table \ref{norm-1_Laney1}, Table \ref{norm-2_Laney1} respectively. Error analysis indicates that ZBS-FDS scheme is a little more accurate than TVS-FDS scheme. Next we consider Lax test case for which initial conditions are given as $(\rho_L, u_L, p_L) = (0.445, 0.698, 3.528)$, $(\rho_R, u_R, p_R) = (0.5, 0.0, 0.571)$ with $x_{o} = 0.5$ and all numerical results are obtained at final time $t=0.15$. The contact and shock discontinuities here are stronger than those in Sod's shock tube problem.  Results of both ZBS-FDS and TVS-FDS schemes are given in Figure \ref{wes2_d}.  For this problem too, ZBS-FDS scheme performs slightly better than TVS-FDS scheme.
 
\begin{table}[!ht]
\caption{$L_{1}$ error norm for the Sod's shock tube problem for both schemes}\label{norm-1_Laney1}
\centerline{%
\begin{tabular}{|c|c|c|c|}
\hline
\textrm{grid points}{}~~ & ZBS-FDS scheme & TVS-FDS scheme  \\
\hline
$40$ & 0.502947  &  0.582406 \\
\hline
$80$ & 0.352076  &  0.397561  \\
\hline
$160$ & 0.235865  & 0.268590  \\
\hline
$320$ & 0.156230  & 0.176909  \\
\hline
$640$ & 0.101461  & 0.114140 \\
\hline
\end{tabular}
}%
\end{table}
\begin{table}[!ht]
\caption{$L_{2}$ error norm for the Sod's shock tube problem for both schemes}\label{norm-2_Laney1}
\centerline{%
\begin{tabular}{|c|c|c|c|}
\hline
\textrm{grid points}{}~~ & ZBS-FDS scheme & TVS-FDS scheme  \\
\hline
$40$ & 0.177260  & 0.196831 \\
\hline
$80$ & 0.134255  & 0.144438  \\
\hline
$160$ & 0.097736  & 0.105535  \\
\hline
$320$ & 0.073471  & 0.078268  \\
\hline
$640$ & 0.055553  & 0.058590 \\
\hline
\end{tabular}
}%
\end{table}
\subsection{Sonic point problem and strong shock problem}
\begin{figure}[!ht]
\centerline{%
\subfigure[]{%
\includegraphics[trim=0 5 35 5, clip, width=0.55\textwidth]{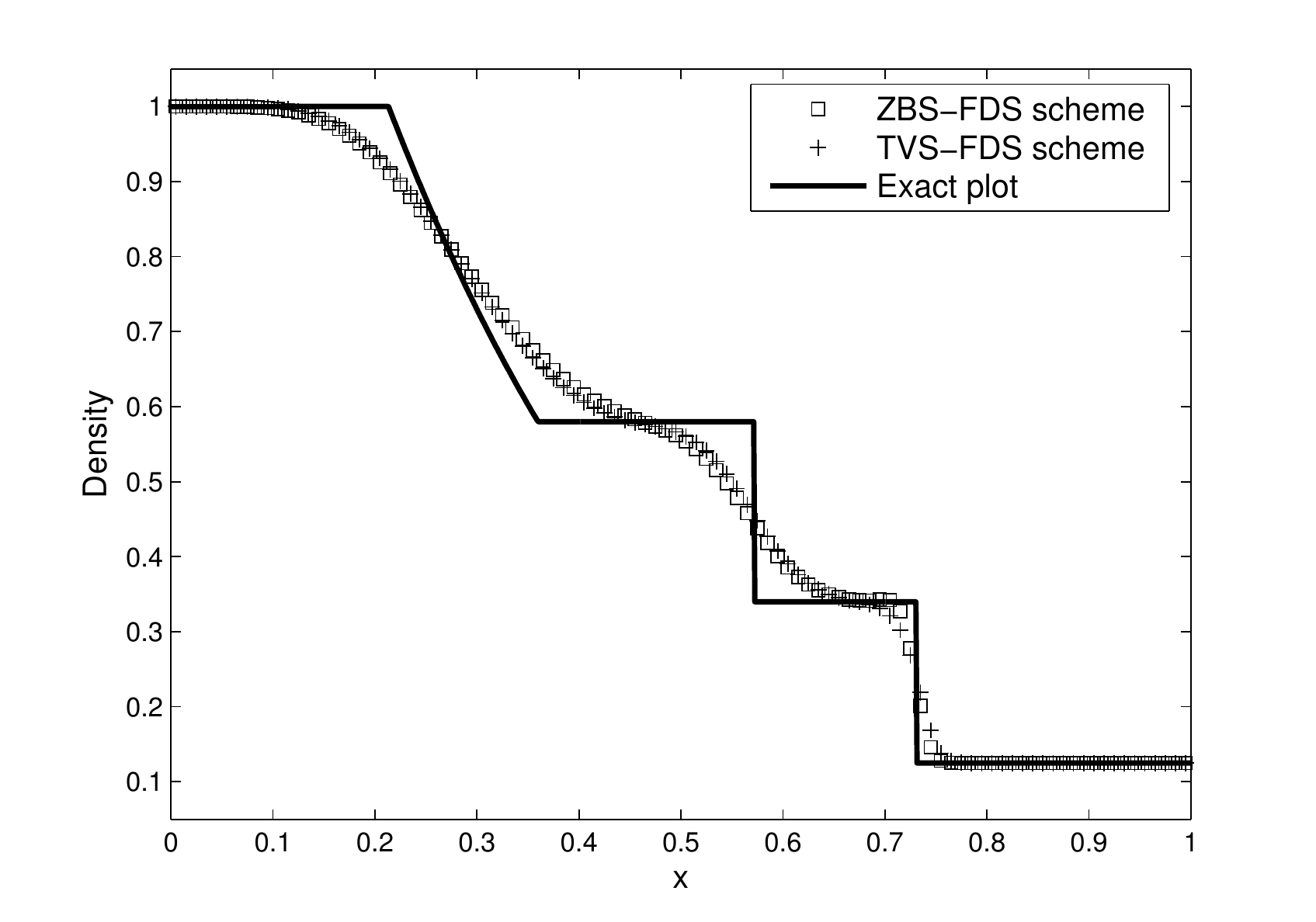}%
\label{Toro1_d.eps}
}%
\subfigure[]{%
\includegraphics[trim=0 5 35 5, clip, width=0.55\textwidth]{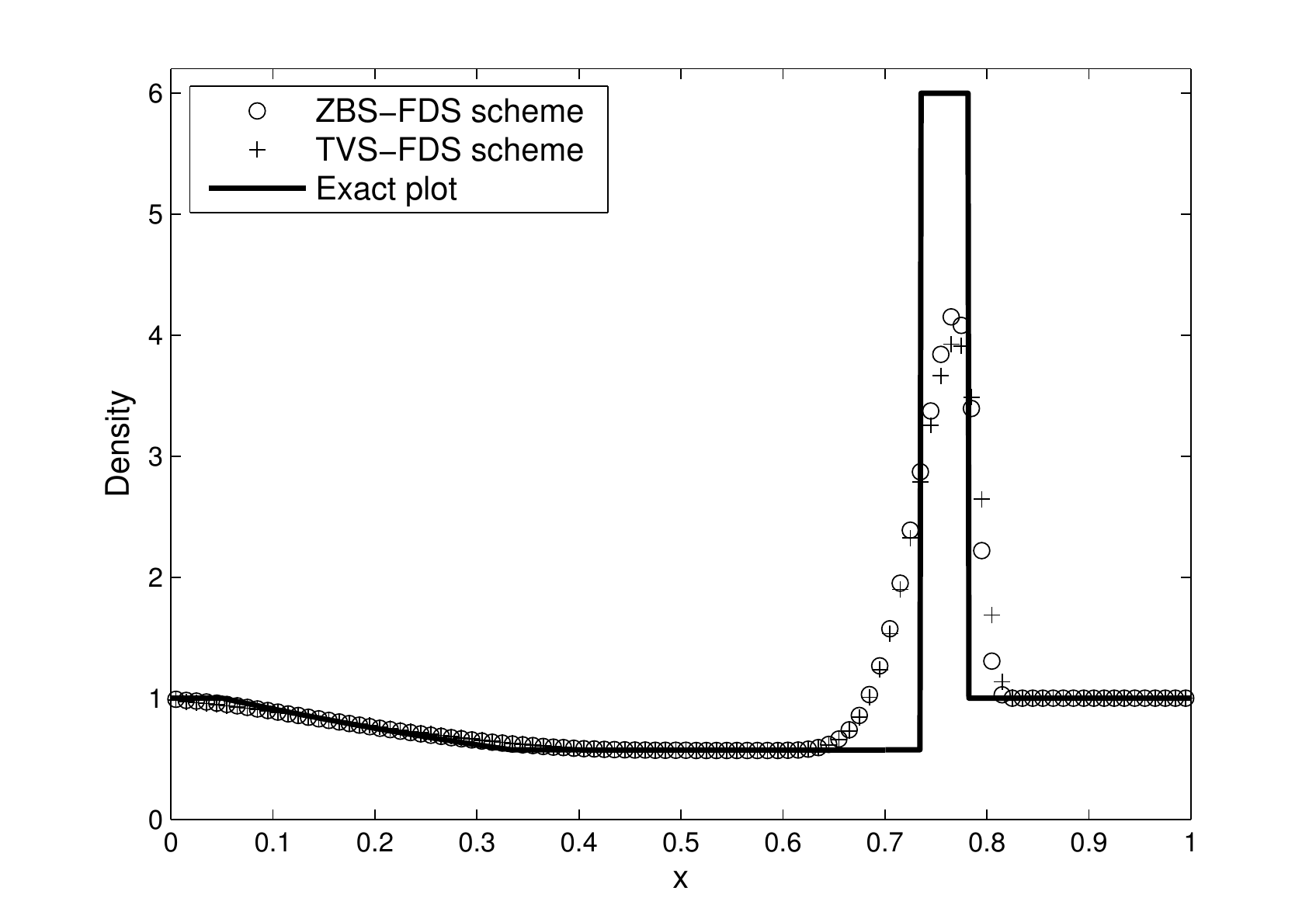}%
\label{toro3_d.eps}
}
}%
\caption{(a) represents density plots for sonic point problem and (b) represents density plots for strong shock problem.}
\end{figure}

Next, we present numerical results of both schemes for a modified version  of Sod's problem. For this problem, solution has a right shock wave, a right travelling contact discontinuity and a left sonic rarefaction wave. This test case is useful in assessing  the entropy condition  satisfaction property of numerical methods. Initial conditions for this problem are given as $(\rho_L, u_L, p_L) = (1.0, 0.75, 1.0)$, $(\rho_R, u_R, p_R) = (0.125, 0.0, 0.1)$ with initial discontinuity at $x_{o} = 0.3$ and all numerical solutions are obtained at final time $t = 0.2$.  For this test case, low diffusive schemes like Roe scheme may violate the entropy condition and give unphysical rarefaction shocks in the expansion region at sonic points.  To avoid this drawback, additional numerical diffusion is typically required, which is usually introduced as an entropy fix and one such famous fix is given by Harten \cite{Harten_entropy_fix}.  Because of sufficient inbuilt numerical diffusion, both ZBS-FDS scheme and TVS-FDS scheme are seen to satisfy the entropy condition, as can be seen in the results shown in Figure \ref{Toro1_d.eps} with no rarefaction shock or non-smoothness being present in the solution.  For this problem, error analysis of $L_1$-norm and $L_2$-norm show TVS-FDS scheme is a slightly more accurate as given in Table \ref{norm-1_Toro1} and Table \ref{norm-2_Toro1}.  Second test case is taken from \cite{Toro_Book_on_Riemann_Solvers} and is designed to assess the robustness and accuracy of numerical methods.  Its solution consists of a strong shock wave with Mach number $198$, a contact discontinuity and a left rarefaction wave.  Initial conditions are given as $(\rho_L, u_L, p_L) = (1.0, 0.0, 1000.0)$, $(\rho_R, u_R, p_R) = (1.0, 0.0, 0.01)$ with $x_{o}=0.5$ and all solutions are obtained at time $t=0.012$.  Both schemes work well and results are given in Figure \ref{toro3_d.eps}. Error analysis of $L_1$-norm and $L_2$-norm show ZBS-FDS scheme is a little more accurate and results are given in Table \ref{norm-1_strong_shock} and Table \ref{norm-2_strong_shock}.
\begin{table}[!ht]
\caption{$L_{1}$ error norm for the sonic point problem for both schemes}\label{norm-1_Toro1}
\centerline{%
\begin{tabular}{|c|c|c|c|}
\hline
\textrm{grid points}{}~~ & ZBS-FDS scheme & TVS-FDS scheme  \\
\hline
$40$ & 0.038718  &  0.036894 \\
\hline
$80$ & 0.028065  &  0.026387  \\
\hline
$160$ & 0.019058  & 0.017863  \\
\hline
$320$ & 0.012698  & 0.011879  \\
\hline
$640$ & 0.008396 & 0.007822 \\
\hline
\end{tabular}
}%
\end{table} 
\begin{table}[!ht]
\caption{$L_{2}$ error norm for sonic point problem for both schemes}\label{norm-2_Toro1}
\centerline{%
\begin{tabular}{|c|c|c|c|}
\hline
\textrm{grid points}{}~~ & ZBS-FDS scheme & TVS-FDS scheme  \\
\hline
$40$ & 0.053061  & 0.051063 \\
\hline
$80$ & 0.043452  & 0.040727  \\
\hline
$160$ & 0.033069  & 0.031264  \\
\hline
$320$ & 0.025245  & 0.024187  \\
\hline
$640$ & 0.019600  & 0.018961 \\
\hline
\end{tabular}
}%
\end{table} 
\begin{table}[!ht]
\caption{$L_{1}$ error norm for the strong shock problem for both schemes}\label{norm-1_strong_shock}
\centerline{%
\begin{tabular}{|c|c|c|c|}
\hline
\textrm{grid points}{}~~ & ZBS-FDS scheme & TVS-FDS scheme  \\
\hline
$40$ & 0.317106  &  0.334709 \\
\hline
$80$ & 0.241142  &  0.258266  \\
\hline
$160$ & 0.180898  & 0.192025  \\
\hline
$320$ & 0.131432  & 0.138044  \\
\hline
$640$ & 0.088449  & 0.092496 \\
\hline
\end{tabular}
}%
\end{table}
\begin{table}[!ht]
\caption{$L_{2}$ error norm for strong shock problem for both schemes}\label{norm-2_strong_shock}
\centerline{%
\begin{tabular}{|c|c|c|c|}
\hline
\textrm{grid points}{}~~ & ZBS-FDS scheme & TVS-FDS scheme  \\
\hline
$40$ & 0.824979  & 0.856110 \\
\hline
$80$ & 0.665983  & 0.699312  \\
\hline
$160$ & 0.558651  & 0.574795  \\
\hline
$320$ & 0.473423  & 0.479052  \\
\hline
$640$ & 0.366668  & 0.372136 \\
\hline
\end{tabular}
}%
\end{table}
\subsection{Stationary contact discontinuity}  
A contact discontinuity occurs when a family of characteristics are parallel to each other in the $x-t$ domain.  Since fluid velocity is the same on both sides, contact discontinuities move with fluid.  The initial conditions as given in \cite{Toro_Book_on_Riemann_Solvers} are $(\rho_L, u_L, p_L) = (1.4, 0.0, 1.0)$ and $(\rho_R, u_R, p_R) = (1.0, 0.0, 1.0)$. The initial discontinuity is present at $x_{o} = 0.5$.  Both ZBS-FDS and TVS-FDS schemes capture the steady contact discontinuity exactly, as shown in Figure \ref{contact_d.eps}.          
\subsection{Strong shock problem with slowly moving contact discontinuity}
\begin{figure}[!ht]
\centerline{%
\subfigure[]{%
\includegraphics[trim=0 5 35 5, clip, width=0.55\textwidth]{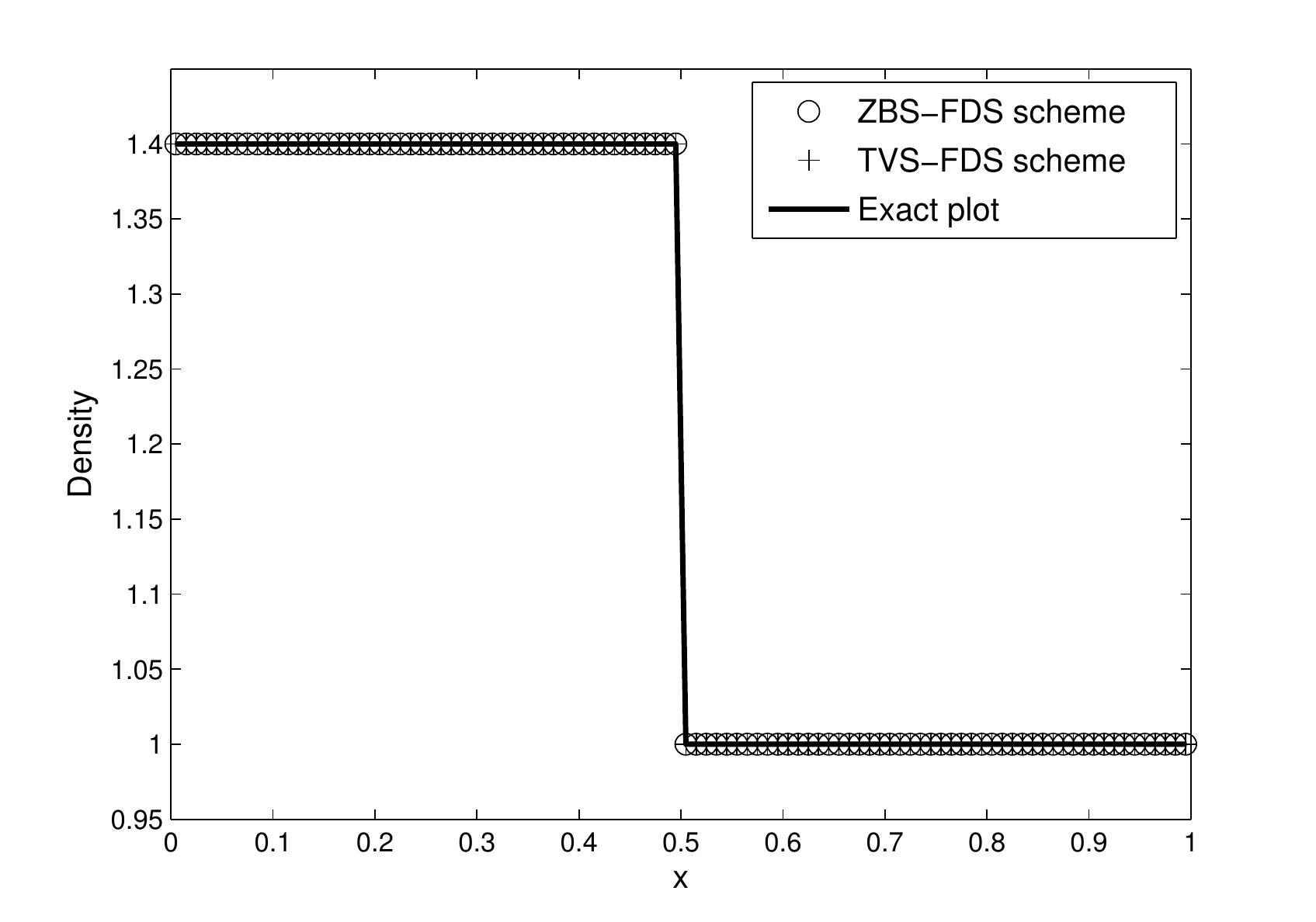}%
\label{contact_d.eps}
}%
\subfigure[]{%
\includegraphics[trim=0 5 35 5, clip, width=0.55\textwidth]{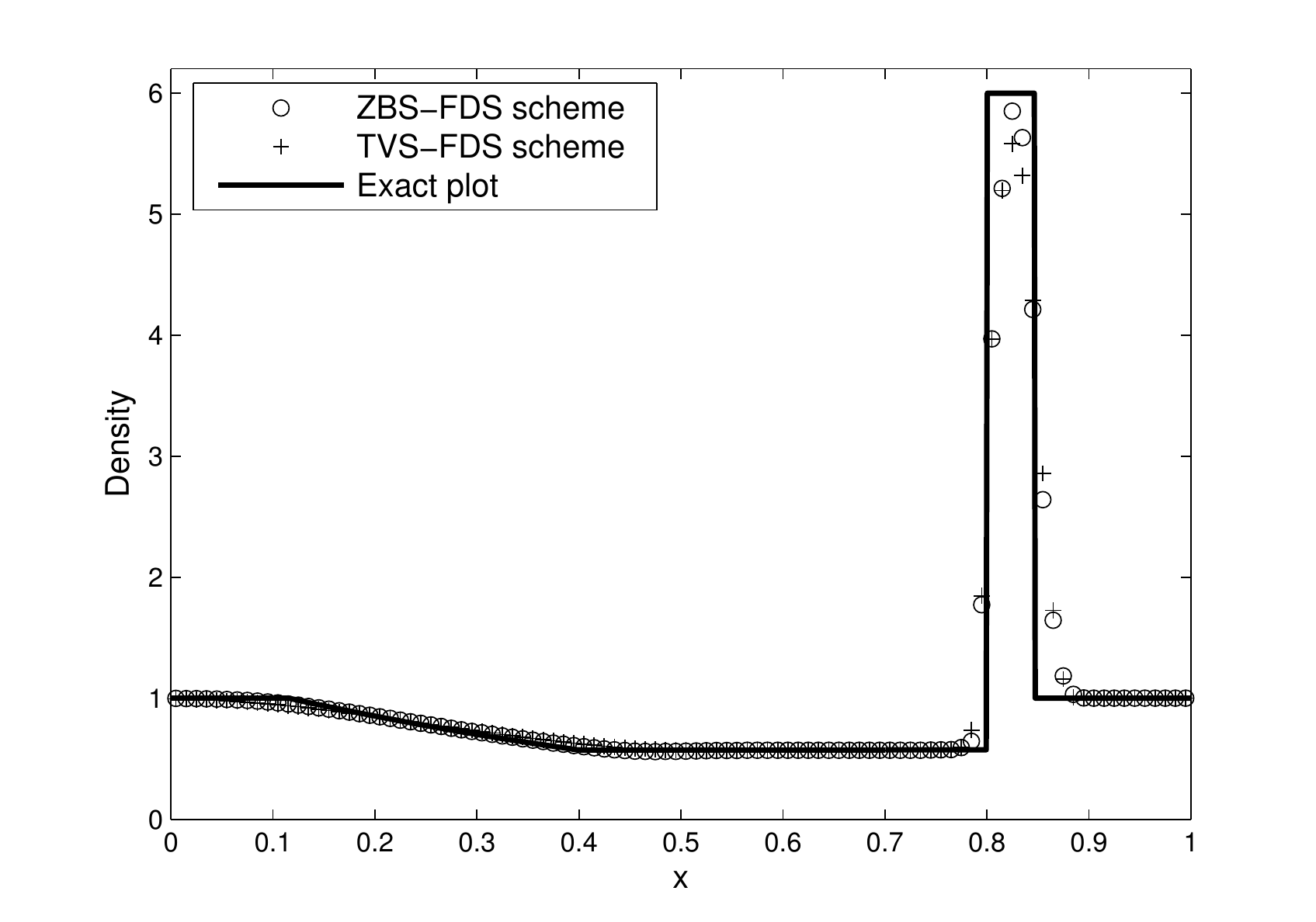}%
\label{toro5_d.eps}
}
}%
\caption{(a) represents density plot for stationary contact discontinuity problem and (b) represents density plots for strong shock problem with slowly moving contact discontinuity.}
\end{figure}

 This test case is also devised to test the robustness of numerical methods but the main reason for devising this test case is to assess the ability of numerical methods to resolve slowly-moving contact discontinuities. The exact solution of this test consists of a left rarefaction wave, a right-travelling shock wave and a slowly moving contact discontinuity.  Initial conditions are given as $(\rho_L, u_L, p_L) \  =  \ (1.0, -19.59745, 1000.0)$, $(\rho_R, u_R, p_R) = (1.0, -19.59745, 0.01)$ with $x_{o}=0.8$ and all numerical solutions are obtained at time $t=0.012$.  In case of ZBS-FDS scheme, numerical solution goes towards the top portion of slowly moving contact wave whereas, for TVS-FDS scheme it is a little below as given in Figure \ref{toro5_d.eps}.  

\subsection{Slowly moving shock} 
\begin{figure}[!ht]
\centerline{%
\subfigure[]{%
\includegraphics[trim=0 5 35 5, clip, width=0.55\textwidth]{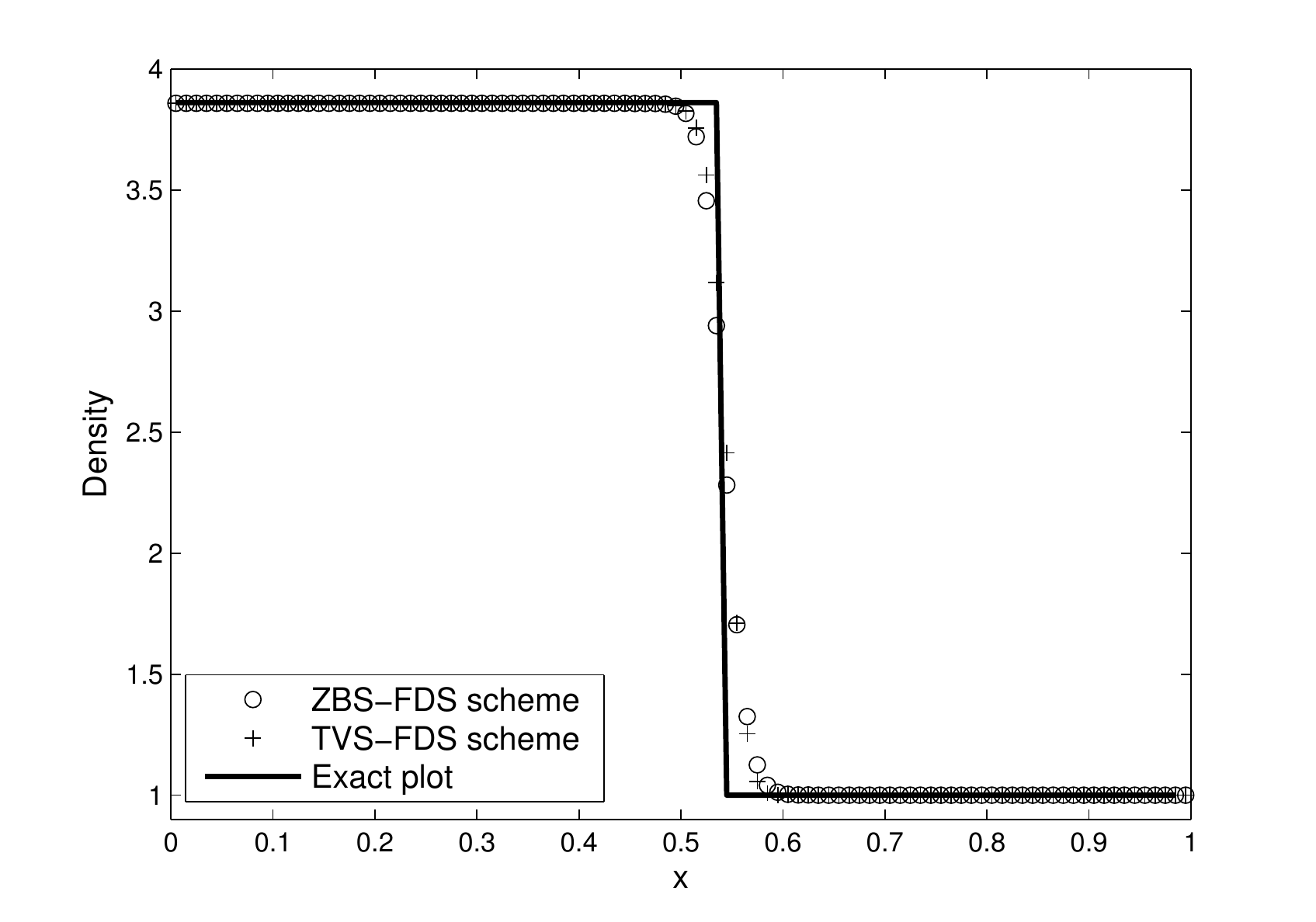}%
\label{MS_d.eps}
}%
\subfigure[]{%
\includegraphics[trim=0 5 35 5, clip, width=0.55\textwidth]{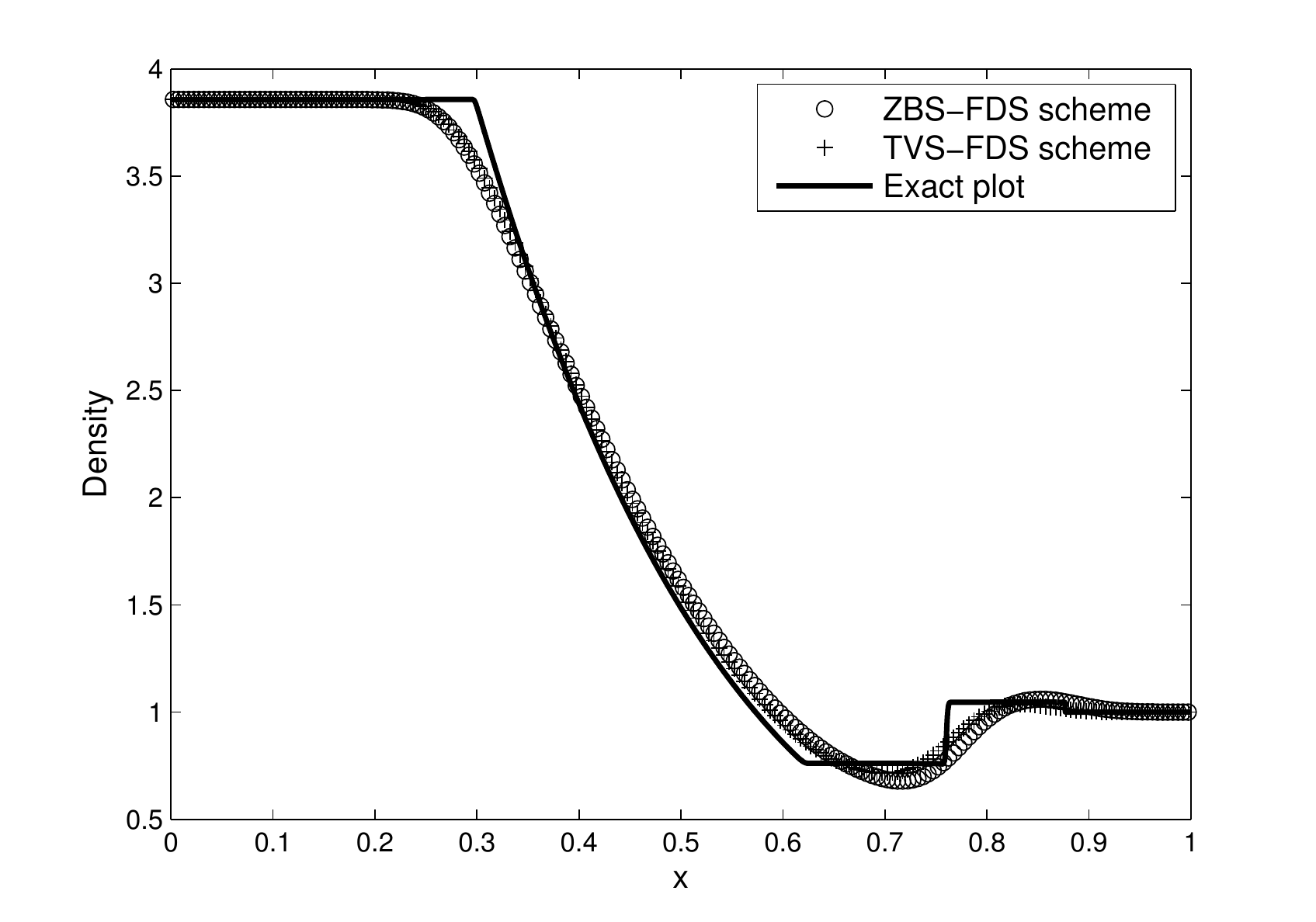}%
\label{wes3_d.eps}
}
}%
\caption{(a) represents density plots for slowly moving shock problem and (b) represents density plots for Mach $3$ problem.}
\end{figure}

Sometimes numerical methods tend to produce oscillations near the shock regions, which are completely unphysical. The oscillations associated with slowly moving shock problems are usually linked with lack of sufficient numerical diffusion in the scheme. We took a test case from \cite{Stiriba_Donat_postshock_oscillations} with  initial conditions as 
$(\rho_L, m_L, E_L) = (3.86, -3.1266, 27.0913)$ and $(\rho_R, m_R, E_R) = (1.0, -3.44, 8.4168)$, where $m = \rho u$ is momentum and $E$ is total energy. Final solutions are obtained at  $t = 4$ units and results are given in Figure \ref{MS_d.eps}. 
    
\subsection{Mach 3 problem}
The initial conditions for this problem are, $(\rho_L, u_L, p_L) = (3.857, 0.92, 10.333)$ and $(\rho_R, u_R, p_R) = (1.0, 3.55, 1.0)$ with $x_{o} = 0.4$ and all solutions are obtained at  $t = 0.1$ units. This problem consists of a supersonic flow with Mach number $3$ in expansion region and it produces a strong expansion fan. Low diffusive upwind schemes such as Roe's approximate solver fail for this problem and require an entropy fix.  According to Wesseling \cite{wesseling}, even after use of Harten's entropy fix, Roe scheme still gives sonic glitch.   Both ZBS-FDS and TVS-FDS schemes perform well without needing any entropy fix and results are given in Figure \ref{wes3_d.eps}.    
\subsection{Interacting blast wave problem}\label{blast_prob}
\begin{figure}[!ht]
\centerline{%
\subfigure[]{%
\includegraphics[trim=0 5 35 5, clip, width=0.55\textwidth]{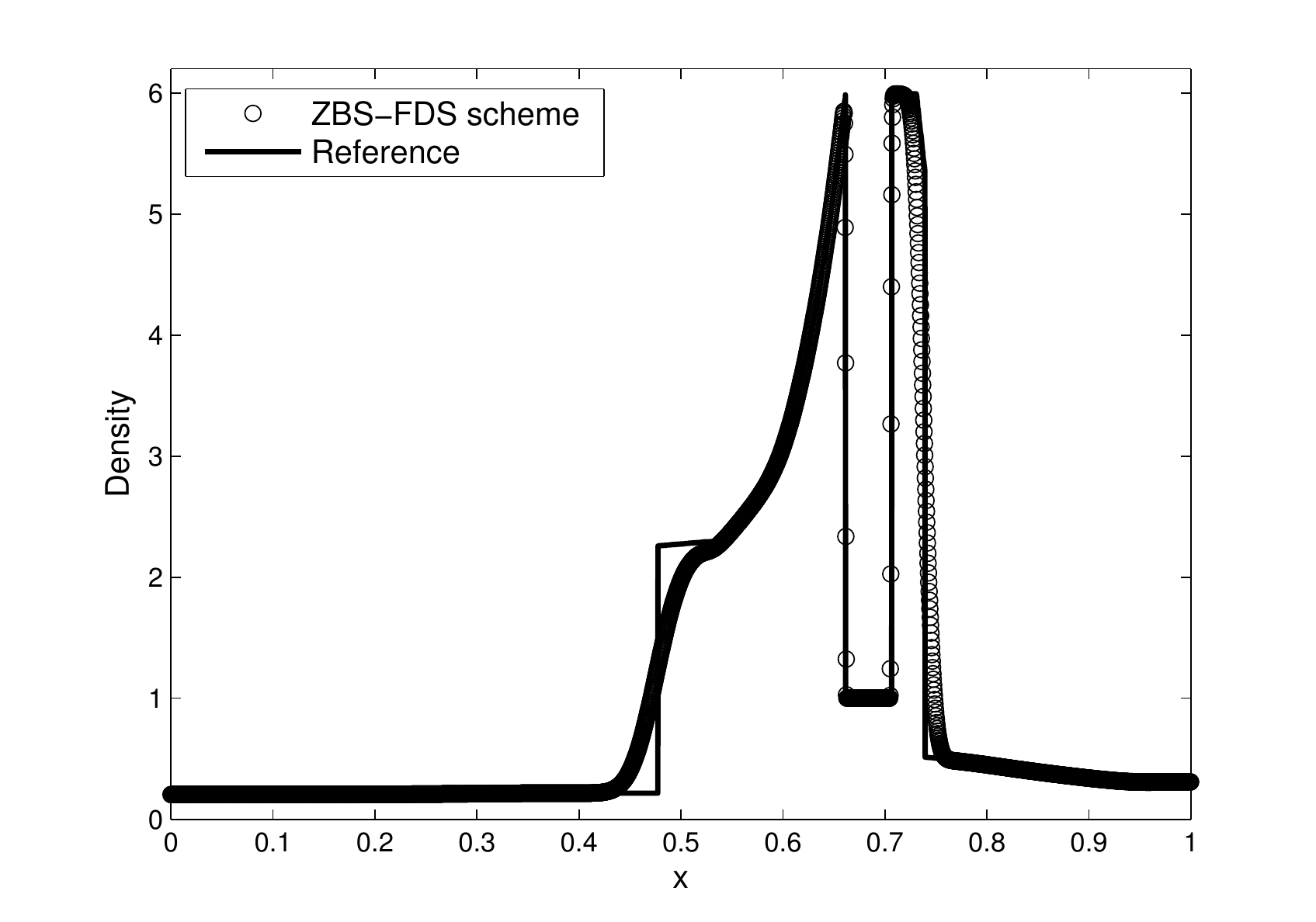}%
\label{blast1_d.eps}
}%
\subfigure[]{%
\includegraphics[trim=0 5 35 5, clip, width=0.55\textwidth]{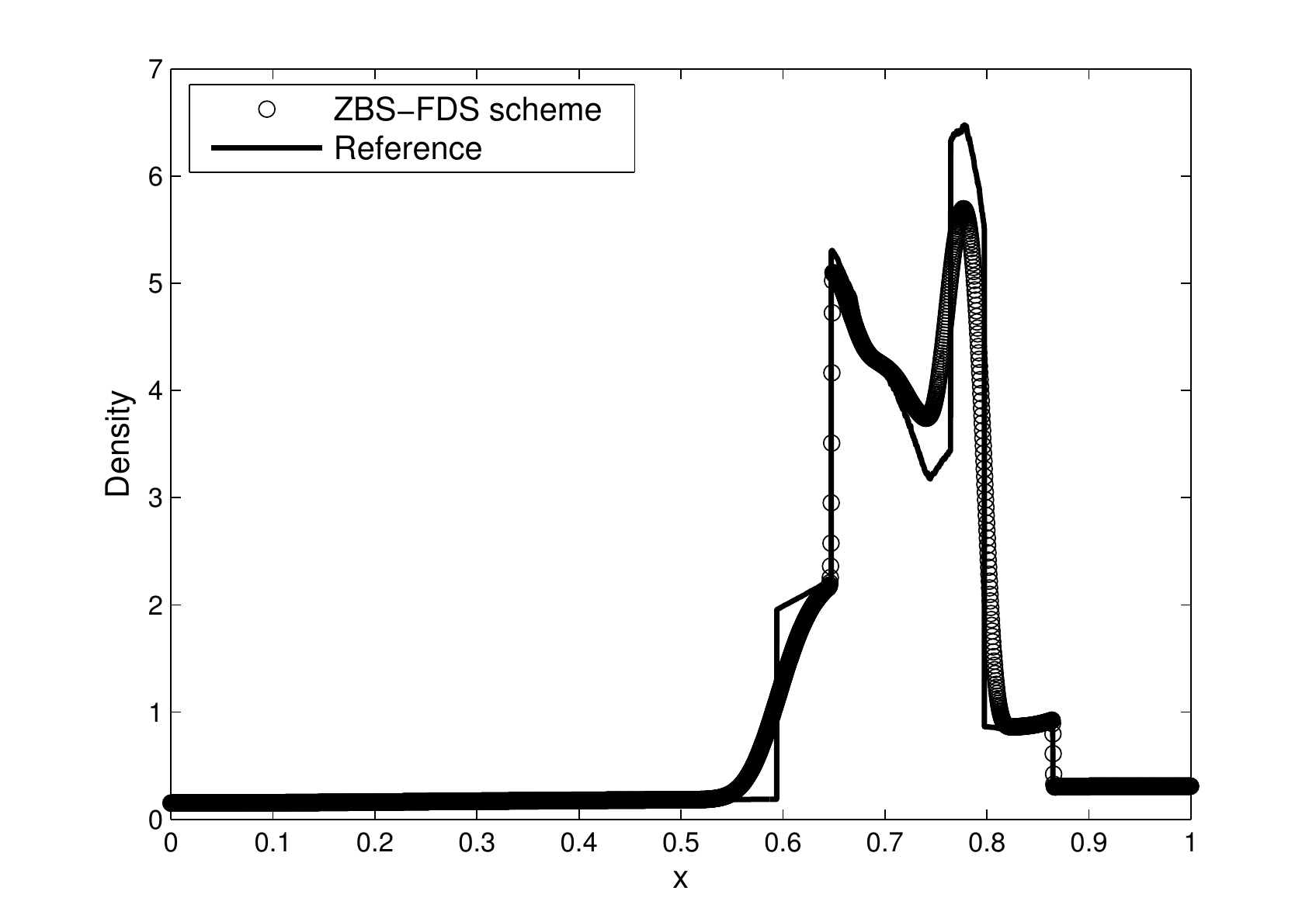}%
\label{blast2_d.eps}
}
}%
\caption{(a) Density plot, for blast wave problem at time t = 0.026 units and (b) density plot, for same problem at time t = 0.038 units.}
\end{figure}

This is one of the most severe test cases used to assess the numerical algorithm for its performance and is taken from Woodward and Colella \cite{Woodward_colella_JCP_1984}.  Computational domain is divided into $3000$ equally spaced finite volumes. Initial conditions for density and velocity are constants and given by $\rho =  1.0$, $u  =  0$. For pressure variable, two discontinuities are present at position $x  =  0.1 \ \textrm{and} \ 0.9$. Initially, $p_{L}$ =  $1000$ if $x  \in  [0.0 , 0.1]$ , $p_{M}$  =  $0.01$ if $x 
\in [0.1,0.9]$ and $p_{R}$ = $100$ if $x  \in [0.9 , 1.0]$. Solution of this problem consists of multiple shocks, contact discontinuities and expansions waves.  Results for ZBS-FDS scheme are given at two different time levels, as shown in Figure \ref{blast1_d.eps} and \ref{blast2_d.eps}.  For this test case, TVS-FDS scheme {\em blew up} in our simulations.  
\subsection{Shock-entropy wave interaction}
\begin{figure}[!ht]
 \centerline{%
\subfigure[]{%
\includegraphics[trim=0 5 35 5, clip, width=0.55\textwidth]{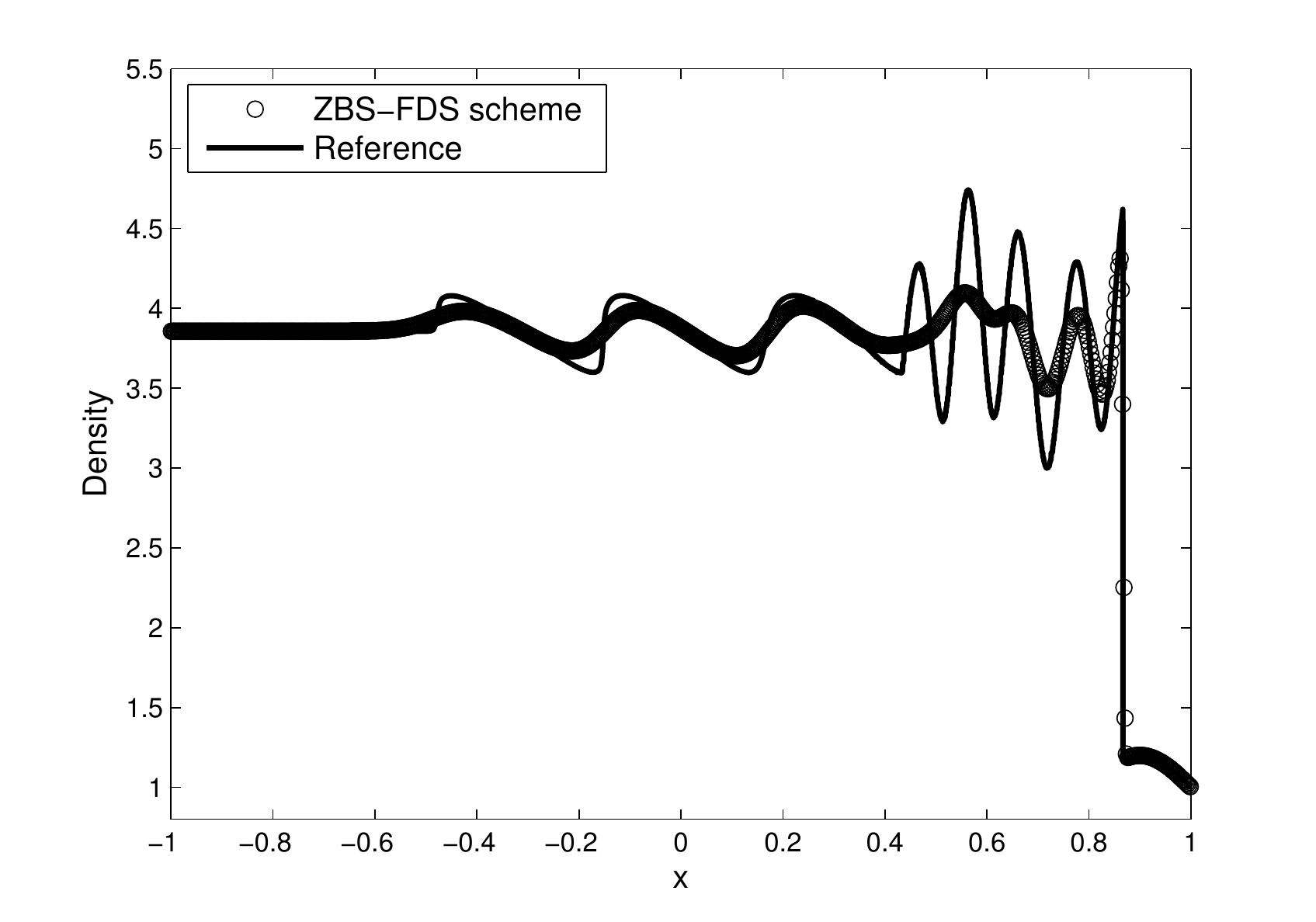}%
\label{ZBS-FDS_shock_entropy_d.eps}
}%
\subfigure[]{%
\includegraphics[trim=0 5 35 5, clip, width=0.55\textwidth]{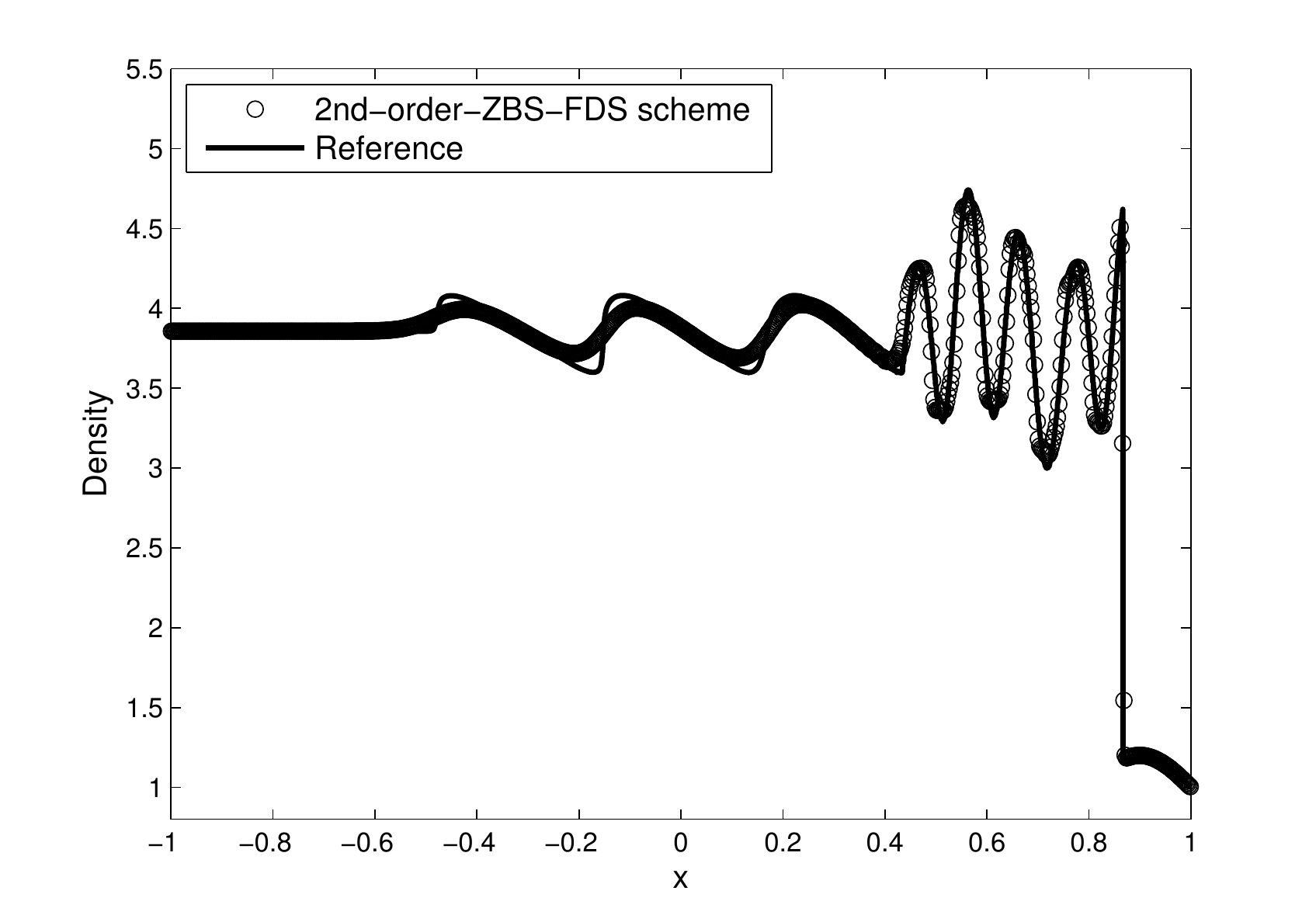}%
\label{ZBS-FDS_shock_entropy_2nd_order_d.eps}
}
}%
\caption{(a) 1st-order results for ZBS-FDS scheme, for shock-entropy wave interaction problem and (b) 2nd-order results for ZBS-FDS scheme for same problem.}
\end{figure}

\begin{figure}[!ht]
\begin{center}
\includegraphics[trim=5 5 35 5, clip, width=0.7\textwidth]{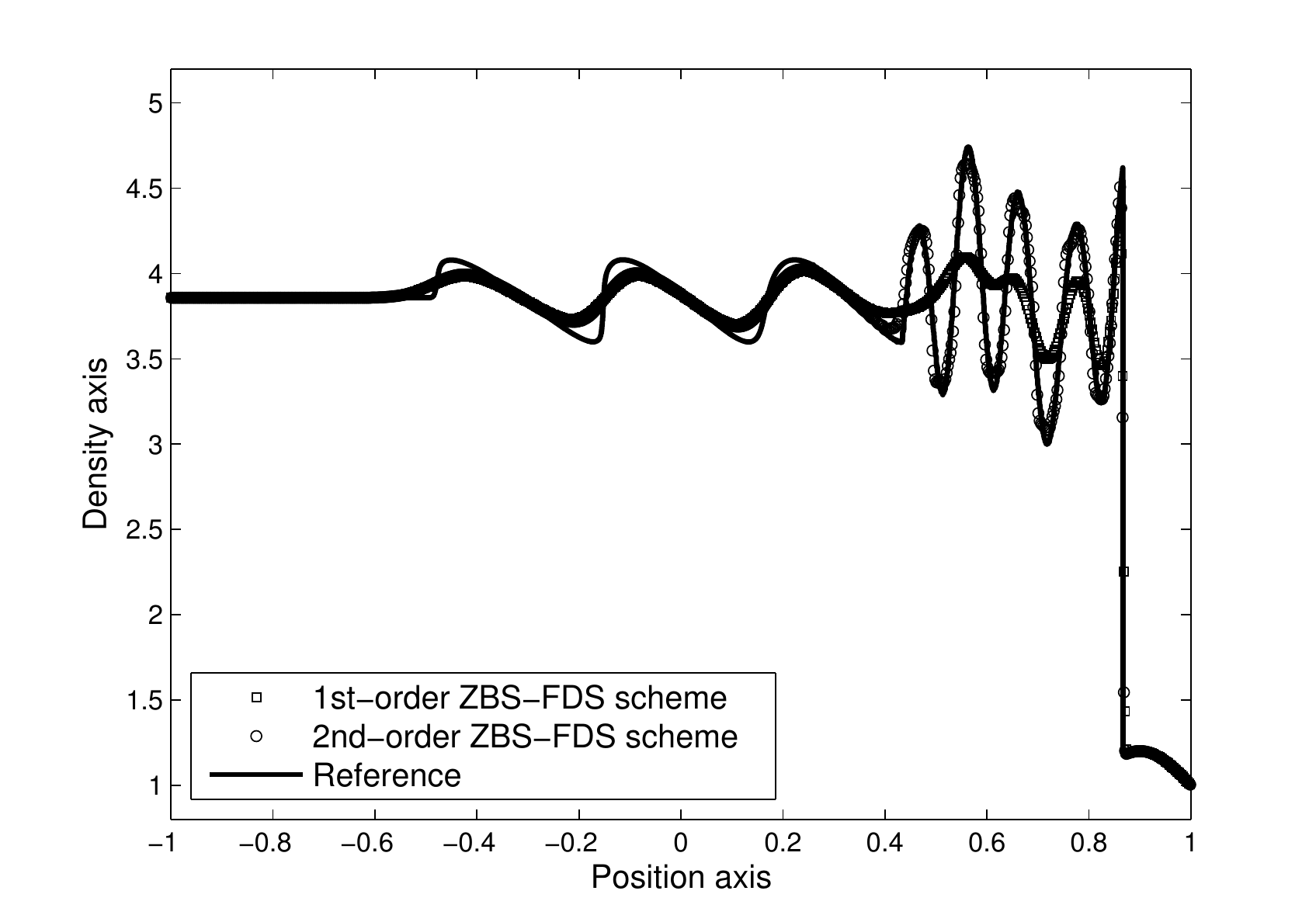}
\caption{Comparison of 1st-order and 2nd-order numerical results of ZBS-FDS scheme for shock-entropy wave interaction problem.}
\label{ZBS-FDS_shock_entropy_comp_d.eps}
\end{center}
\end{figure}

\begin{figure}[!ht]
 \centerline{%
\subfigure[]{%
\includegraphics[trim=0 5 35 5, clip, width=0.55\textwidth]{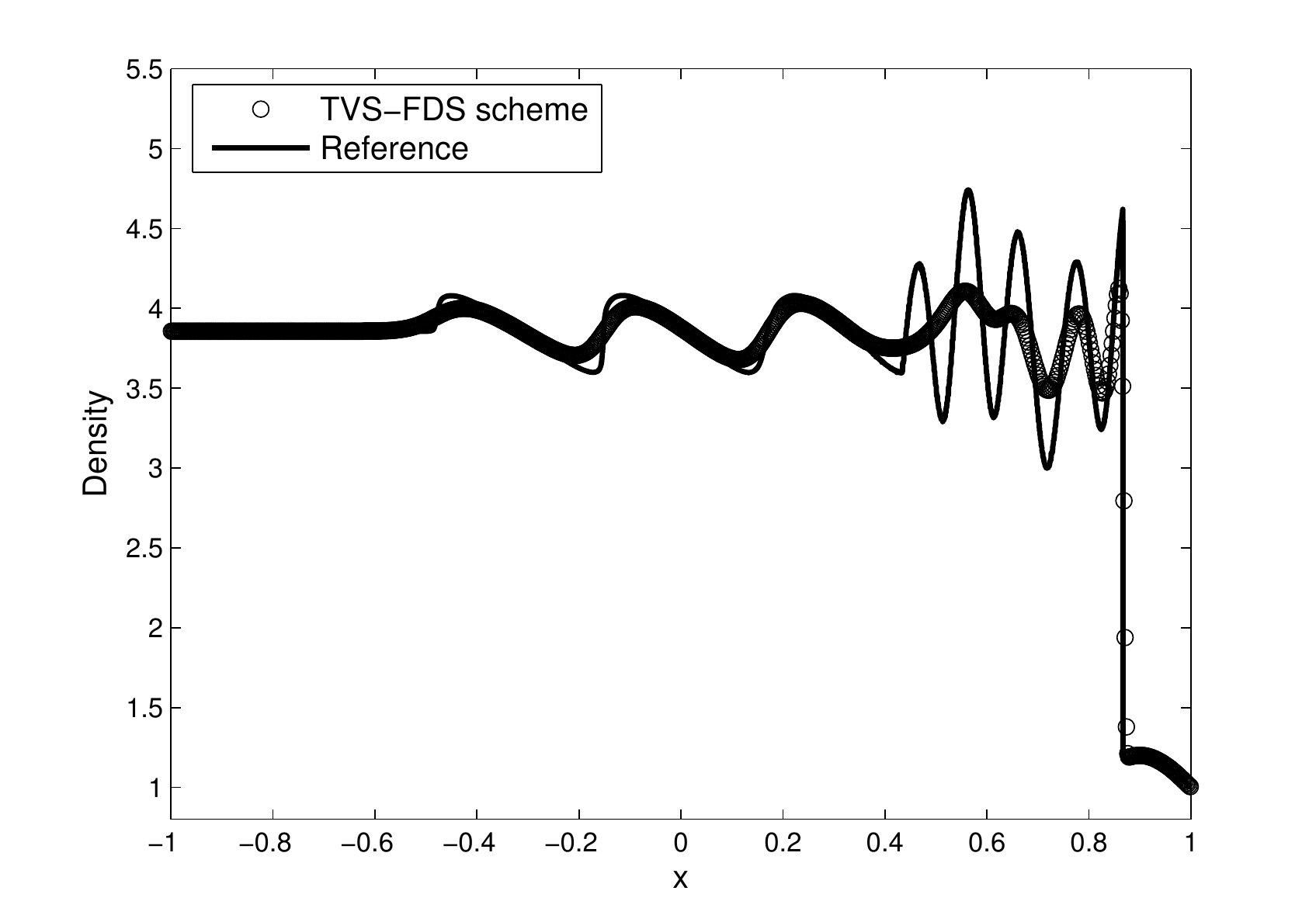}%
\label{TVS-FDS_shock_entropy_d.eps}
}%
\subfigure[]{%
\includegraphics[trim=0 5 35 5, clip, width=0.55\textwidth]{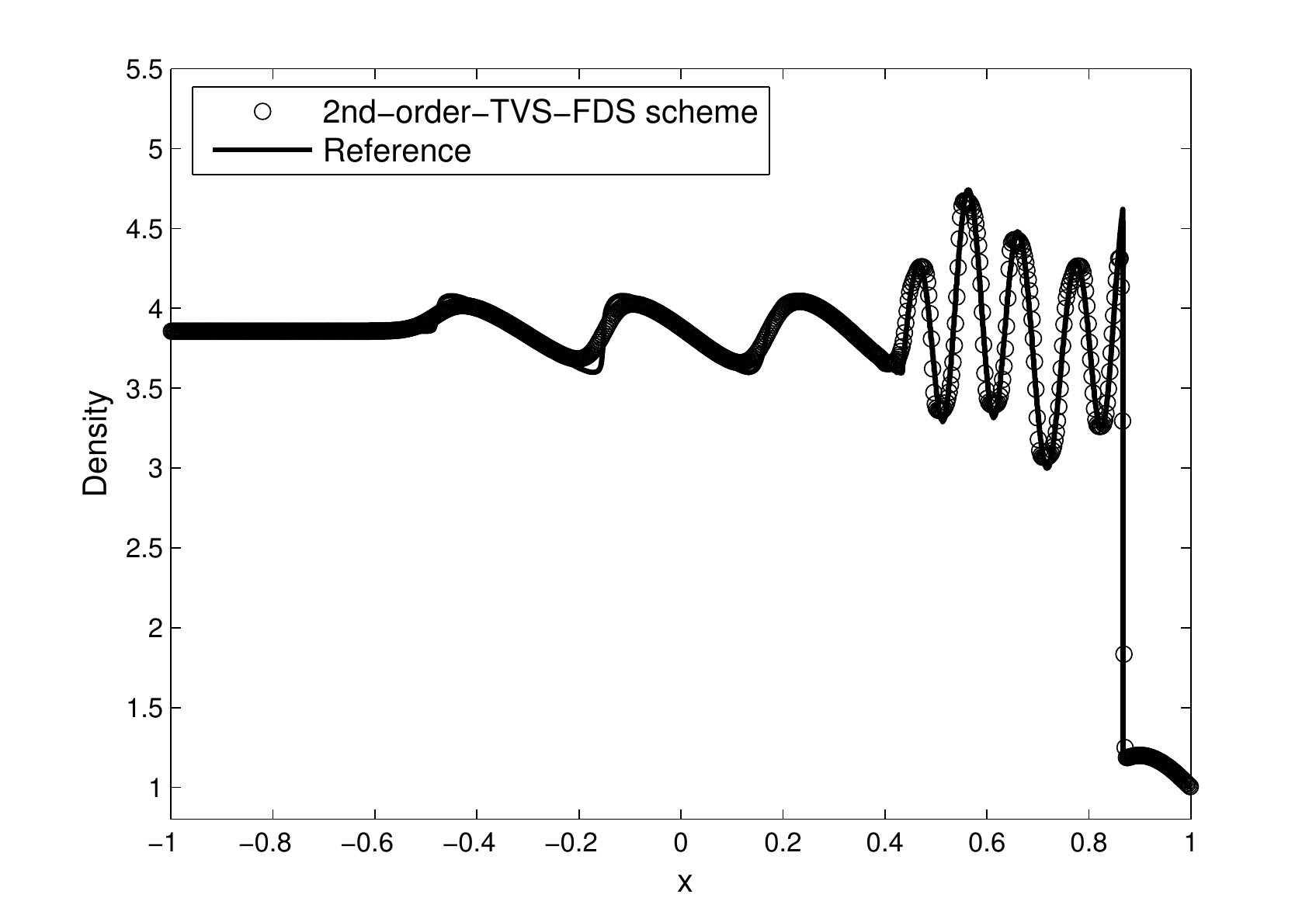}%
\label{shock_entropy_2nd_TVS_d.eps}
}
}%
\caption{(a) Density plot for 1st-order TVS-FDS scheme, for shock-entropy wave interaction problem and (b) density plot for 2nd-order TVS-FDS scheme for same problem.}
\end{figure}

The shock-entropy wave interaction test case considered here for testing the present schemes is taken from \cite{Balsara and Shu}, with computational domain $x \in [-1,1]$ being divided into $800$ equally spaced finite volumes and all solutions are obtained at final time $t = 0.47$.  The initial conditions are given below.  
\begin{align}
\begin{split}
   (\rho_{L}, u_{L}, p_{L}) \ &= \ \big[3.857143, 2.629369, 10.3333\big] \ \ \textrm{if}   \  \  x \ < \ -0.8   \\
   (\rho_{R}, u_{R}, p_{R}) \ &= \ \big[1 + 0.2\sin(5\pi x), 0, 1\big] \ \ \textrm{if}    \  \   x \ > \ -0.8.
\end{split}
\end{align}
In this problem, a Mach 3 shock wave interacts with density disturbances created by perturbing the initial density. This initial disturbance gives rise to the continuous interaction of smooth flow with the discontinuities. Similar kind of interaction can be observed in compressible turbulence. Therefore, this is a  suitable test case to test the scheme for its ability to resolve complex interactions, which can be used in turbulent computations.  First order results for ZBS-FDS scheme are presented in Figure \ref{ZBS-FDS_shock_entropy_d.eps}. To achieve second order accuracy, we used Venkatakrishnan's limiter which is a modified version of van Albada limiter \cite{Venkatakrishnan_limiter} and deals with piecewise linear reconstruction of primitive variables. As an example, let us consider a piecewise linear reconstruction for density variable, {\em i.e.}, to obtain 
\begin{equation}
 \rho_{i+1/2}^{L} \ = \  \rho_{i} \ + \ \frac{1}{2} \frac{\left(\Delta_{+}^{2} + \epsilon^{2}\right)\Delta_{-} \ + \ \left(\Delta_{-}^{2} + \epsilon^{2}\right)\Delta_{+}}{\Delta_{+}^{2} + \Delta_{+}^{2} + 2\epsilon^{2}}
\end{equation}
where,
 \begin{align}
 \begin{split}
\Delta_{+}  \ &= \ \rho_{i+1} - \rho_{i}    \\
\Delta_{-}  \ &= \ \rho_{i} - \rho_{i-1}   
\end{split}
\end{align}
and
\begin{equation}
\epsilon^{2}  \ = \   (K\Delta{x})^{3} 
\end{equation}
Similarly, other primitive variables can be reconstructed. Here, $K$ is a constant and $\Delta{x}$ is a grid spacing. Large values of $K$ indicate no limiting and in the present case, we take $K = 0.1$.  For this problem, second order results for ZBS-FDS scheme are computed and are given in Figure \ref{ZBS-FDS_shock_entropy_2nd_order_d.eps}.  Comparison of both first order results and second order results for ZBS-FDS scheme are given in Figure \ref{ZBS-FDS_shock_entropy_comp_d.eps}.  Results for TVS-FDS scheme are presented in Figure \ref{TVS-FDS_shock_entropy_d.eps}, \ref{shock_entropy_2nd_TVS_d.eps}.  Both the schemes produce results of nearly similar accuracy for this test case. In both cases, second order accurate results are substantially better, compared to the first order accurate results.

\section{Two-dimensional Euler system} 
The 2-D Euler equations form a system of four coupled non-linear hyperbolic PDEs with independent space variables $x,y$ and independent time variable $t$. In the differential, as well as conservative, form the 2-D Euler system can be written as 
\begin{equation}\label{2-D_differential_form}
 \dfrac{\partial \boldsymbol{U}}{\partial{t}}  \ + \ \dfrac{\partial \boldsymbol{F}_{1}}{\partial{x}} \ + \  \dfrac{\partial \boldsymbol{F}_{2}}{\partial{y}}  \ = \ \boldsymbol{0}
\end{equation}
where, $\boldsymbol{U}$ is vector of conserved variables and $\boldsymbol{F}_{1}$, $\boldsymbol{F}_{2}$ are flux vectors which are given as follows.
\begin{equation} 
 \boldsymbol{U} = \begin{bmatrix} 
                 \rho \\[0.4em]  
		 \rho u \\[0.4em] 
		 \rho v \\[0.4em]
		 \rho E 
		\end{bmatrix} \ \textrm{,} \ 
 \boldsymbol{F}_{1}  = \begin{bmatrix} 
                  \rho u \\[0.4em] 
	          \rho u^{2} + p  \\[0.4em]
	          \rho uv      \\[0.4em]
	          \rho u E + p u 
				\end{bmatrix}  \  \textrm{and} \
\boldsymbol{F}_{2}  = \begin{bmatrix} 
                    \rho v \\[0.4em] 
	            \rho uv   \\[0.4em]
	            \rho v^{2} + p  \\[0.4em]
	            \rho v E + p v
				\end{bmatrix}
\end{equation} 
Using the divergence form, equation (\ref{2-D_differential_form}) can be written as
\begin{equation}\label{FV_2d_Euler}
 \dfrac{\partial  \boldsymbol{U}}{\partial{t}}  \ + \   {\nabla}\boldsymbol{.F}  \ = \ \boldsymbol{0}
\end{equation}
where 
\begin{equation} 
 \boldsymbol{F} = \begin{bmatrix} 
                 \rho u_{\bot} \\[0.4em]  
		 \rho u u_{\bot} + p n_{x} \\[0.4em] 
		 \rho v u_{\bot} + p n_{y} \\[0.4em]
		 \rho E u_{\bot} + p u_{\bot}
		\end{bmatrix} 
\end{equation} 
is the flux vector and the vector $u_{\bot}$ is defined as the scalar product of the velocity vector $\boldsymbol{u}$ and the unit normal vector $\boldsymbol{n}$, {\em i.e.}, 
\begin{equation}
 u_{\bot} \ = \ \boldsymbol{u.} \boldsymbol{n} \ = \ n_{x} u +  n_{y} v
\end{equation}
where $n_x$ and $n_y$ represent the direction cosines of the unit normal $\hat{n}$ to the cell-interface and are given by 
\begin{equation}
 n_x \ = \ \dfrac{\Delta{y}}{\Delta{s}} \ , \   n_y \ = \ -\dfrac{\Delta{x}}{\Delta{s}} 
\end{equation}

On integrating (\ref{FV_2d_Euler}) over domain $\varOmega$ with boundary $\partial {\varOmega}$ and on further using Green's theorem, we get 
\begin{equation}
 \dfrac{\partial}{\partial t} \int_{\varOmega} \boldsymbol{U} d\varOmega  \ + \  \oint_{\partial {\varOmega}} \boldsymbol{F} ds  \ = \ \boldsymbol{0}
\end{equation}
On calculating average value of $\boldsymbol{U}$ over $\varOmega$, the first integral can be re-written and the above equation becomes 
\begin{equation}
 \dfrac{\partial \boldsymbol{\bar{U}}}{\partial t}  \ = \   - \dfrac{1}{A} \oint_{\partial {\varOmega}} \boldsymbol{F} ds 
\end{equation}
where $A$ is area enclosed by $\varOmega$.
For typical two dimensional finite volume, for general quadrilaterals, integral as given by second term of above semi-integral from can be approximated by line integrals.  After a little algebra, the discretized finite volume form for 2-D Euler system is given by 
\begin{equation}
 \dfrac{\partial \boldsymbol{\bar{U}}_{m}}{\partial t}  \ = \  - \frac{1}{A_m} \  \sum_{k} \boldsymbol{F}^{k} \Delta{s}^{k}
\end{equation}
\begin{figure}[!ht]
\begin{center}
\includegraphics[trim=0 0 0 0, clip, width=0.55\textwidth]{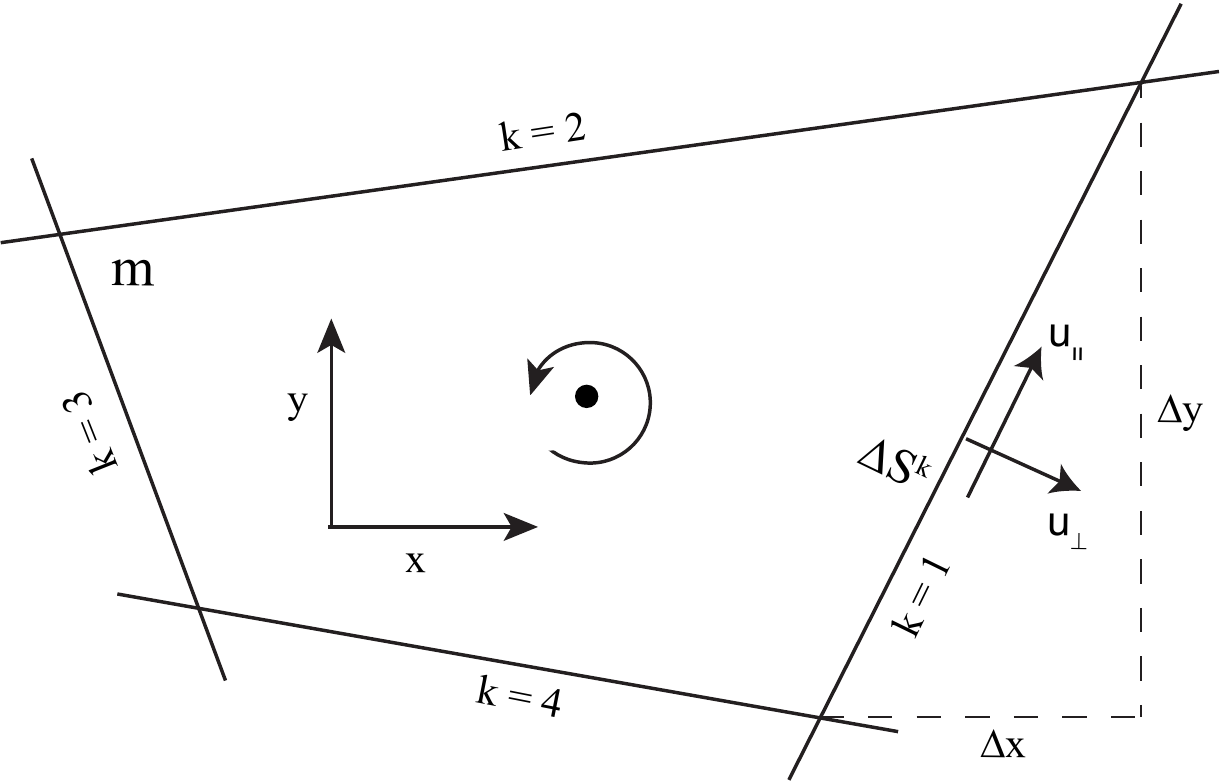}
\caption{Schematic representation of general quadrilateral.}
\label{gen_quard-crop}
\end{center}
\end{figure}
where subscript $m$ denotes the cell number, $k$ is the cell-interface index of the $m^{th}$ cell.  Similarly, $\boldsymbol{F}^{k}$ and $\Delta{s}^{k}$ are the normal flux and the perimeter of the $k^{th}$ face. Further, $A_{m}$ is the area of $m^{th}$ cell.\\

In 1-D, both ZBS-FDS and TVS-FDS schemes performed reasonably well for most of the important test cases except for the blast wave problem (\ref{blast_prob}), where TVS-FDS {\em blew up} quickly. Otherwise, it is bit difficult to judge which of two is more accurate as ZBS-FDS scheme produces slightly better results for Sod's shock tube problem, Lax problem and for both strong shock problems, whereas TVS-FDS scores over ZBS-FDS in case of sonic point problem, slowly moving shock problem and for Mach 3 problem.  In this work, we opt for Zha and Bilgen type splitting in formulating FDS concept based numerical scheme for 2-D Euler system, while emphasizing that this choice is purely based on our convenience and in future other possibility can be explored.  

\subsection{Analysis of Zha and Bilgen type splitting in 2-D}
The flux vector $\boldsymbol{F}$ in the 2-D case is split into a convection part and a pressure part, based on Zha and Bilgen splitting, as follows.  
\begin{equation}
 \boldsymbol{F}  \ = \   \boldsymbol{F}_{c}^{\boldsymbol{ZB}}  \ + \  \boldsymbol{F}_{p}^{\boldsymbol{ZB}}
\end{equation}
where 
\begin{equation} 
\boldsymbol{F}_{c}^{\boldsymbol{ZB}}  = \begin{bmatrix}
                 \rho u_{\bot}   \\[0.4em]  
		 \rho u u_{\bot}  \\[0.4em] 
		 \rho v u_{\bot}   \\[0.4em]
		 \rho E u_{\bot}
		\end{bmatrix}   \  \textrm{and} \
\boldsymbol{F}_{p}^{\boldsymbol{ZB}}  = \begin{bmatrix} 
                  0  \\[0.4em]  
		  p n_{x} \\[0.4em] 
		  p n_{y} \\[0.4em]
		  p u_{\bot}
		\end{bmatrix} 
\end{equation} 
Let $\boldsymbol{A}_{c}^{\boldsymbol{ZB}}$ denote convection flux Jacobian matrix which is given below.
\begin{equation}
\boldsymbol{A}_{c}^{\boldsymbol{ZB}}  \ = \ \begin{bmatrix}
         \ 0  &&  n_x  &&   n_y  &&  0   \\[0.4em]
         \ -u u_\bot   &&  u_\bot + u n_x  &&  u n_y  &&  0   \\[0.4em]
         \  -v u_\bot  &&   v n_x   &&  u_\bot + v n_y  && 0  \\[0.4em]
         \  -E u_\bot  &&  E n_x   &&  E n_y  &&  u_\bot
        \end{bmatrix}
\end{equation}
Eigenvalues corresponding to matrix $\boldsymbol{A}_{c}^{\boldsymbol{ZB}}$ are real and equal with set of eigenvalues as $u_\bot, u_\bot, u_\bot, u_\bot$. Analysis of $\boldsymbol{A}_{c}^{\boldsymbol{ZB}}$ shows that it has a defective set of LI eigenvectors, {\em i.e.},
\begin{equation}
 \boldsymbol{R}_{c,1}^{\boldsymbol{ZB}} \ = \ \begin{bmatrix}
                  \ n_x    \\[0.4em]
                  \ u_\bot  \\[0.4em]
                  \ 0        \\[0.4em]
                  \ 0 
                 \end{bmatrix} \
\ \textrm{,} \
\boldsymbol{R}_{c,2}^{\boldsymbol{ZB}} \ = \ \begin{bmatrix}
                  \ n_y    \\[0.4em]
                  \ 0  \\[0.4em]
                  \ u_\bot        \\[0.4em]
                  \ 0 
                 \end{bmatrix} \ \textrm{and} \
\boldsymbol{R}_{c,3}^{\boldsymbol{ZB}} \ = \ \begin{bmatrix}
                  \ 0    \\[0.4em]
                  \ 0  \\[0.4em]
                  \ 0        \\[0.4em]
                  \ 1 
                 \end{bmatrix} 
\end{equation}
On evaluating rank of matrices $(\boldsymbol{A}_{c}^{\boldsymbol{ZB}}- u_\bot \boldsymbol{I_4})$, $(\boldsymbol{A}_{c}^{\boldsymbol{ZB}}- u_\bot \boldsymbol{I_4})^{2} ...$, we find that there will be one Jordan block of order two as $rank(\boldsymbol{A}_{c}^{\boldsymbol{ZB}}- u_\bot \boldsymbol{I_4})^2 \ = \ rank(\boldsymbol{A}_{c}^{\boldsymbol{ZB}}- u_\bot \boldsymbol{I_4})^3$. Let $R(\boldsymbol{A})$ denote the space spanned by the columns of matrix $\boldsymbol{A}_{c}^{\boldsymbol{ZB}} -u_\bot \boldsymbol{I_4}$. Then, as explained in 1-D case, we have
\begin{equation}
 R(\boldsymbol{A}) \ = \ x_{1}\boldsymbol{A}_1 \ + \ x_{2}\boldsymbol{A}_2 \ + \ x_{3}\boldsymbol{A}_3 \ + \ x_4 \boldsymbol{A}_4
\end{equation}
where, $\boldsymbol{A}_1,\boldsymbol{A}_2,\boldsymbol{A}_3 \ and \ \boldsymbol{A}_4$ are column vectors of $\boldsymbol{A_{c}^{ZB}} -u_\bot \boldsymbol{I_4}$. Now
\begin{equation}
 R(\boldsymbol{A}) \ = \ x_{1}\begin{bmatrix}
                  \ -u_\bot    \\[0.4em]
                  \ -u u_\bot    \\[0.4em]
                  \ -v u_\bot    \\[0.4em]
                  \ -E u_\bot  \
                 \end{bmatrix} \ + \ 
               x_{2}\begin{bmatrix}
                  \ n_x    \\[0.4em]
                  \ u n_x    \\[0.4em]
                  \ v n_x   \\[0.4em]
                  \ E n_x  \
                 \end{bmatrix}  \ + \ 
                 x_{3}\begin{bmatrix}
                  \ n_y    \\[0.4em]
                  \ u n_y   \\[0.4em]
                  \ v n_y   \\[0.4em]
                  \ E n_y   \
                 \end{bmatrix} \ + \ 
                 x_{4}\begin{bmatrix}
                  \ 0    \\[0.4em]
                  \ 0   \\[0.4em]
                  \ 0   \\[0.4em]
                  \ 0   \
                 \end{bmatrix}
\end{equation} 
or
\begin{equation}
 R(\boldsymbol{A})  \ = \ (-u_\bot x_{1} + n_x x_2 + n_y x_3)
                 \begin{bmatrix}
                  \ 1    \\[0.4em]
                  \ u    \\[0.4em]
                  \ v    \\[0.4em]
                  \ E
                  \end{bmatrix}
\end{equation} 
The column vector $(1,u,v,E)^{t}$, which is a range space of $R(\boldsymbol{A})$, becomes generalized eigenvector of $\boldsymbol{A}_{c}^{\boldsymbol{ZB}}$ as null space $N(\boldsymbol{AX})$ is just a scalar coefficient.
Let us take  $\boldsymbol{A}_{c}^{\boldsymbol{ZB}}\boldsymbol{X}_1  $ which is equal to 
\begin{equation} 
   \begin{bmatrix} \ 
       \ 0  &&  n_x  &&   n_y  &&  0   \\[0.4em]
         \ -u u_\bot   &&  u_\bot + u n_x  &&  u n_y  &&  0   \\[0.4em]
         \  -v u_\bot  &&   v n_x   &&  u_\bot + v n_y  && 0  \\[0.4em]
         \  -E u_\bot  &&  E n_x   &&  E n_y  &&  u_\bot
        \end{bmatrix}
       \begin{bmatrix}
                  \ 1    \\[0.4em]
                  \ u    \\[0.4em]
                  \ v    \\[0.4em]
                  \ E 
                 \end{bmatrix} 
\end{equation} 
which, on solving, is equal to 
\begin{equation} 
                u_\bot \begin{bmatrix}
                  \ 1    \\[0.3em]
                  \ u    \\[0.3em]
                  \ v   \\[0.3em]
                  \ E  
                 \end{bmatrix} 
\end{equation}
Thus, $\boldsymbol{A}_{c}^{\boldsymbol{ZB}}\boldsymbol{X}_1  \ = \ u_\bot \boldsymbol{X}_1 $ holds.
Now, this generalized eigenvector is expected to form a Jordan chain of order two corresponding to matrix $\boldsymbol{A}_{c}^{\boldsymbol{ZB}}$, {\em i.e.},
\begin{align}
 \begin{split}
 \boldsymbol{A}_{c}^{\boldsymbol{ZB}}\boldsymbol{X}_1 \ &= \ u_\bot \boldsymbol{X}_1  \\
  \boldsymbol{A}_{c}^{\boldsymbol{ZB}}\boldsymbol{X}_2 \ &= \ u_\bot \boldsymbol{X}_{2} \ + \ \boldsymbol{X}_1
 \end{split}
\end{align}
Other generalized eigenvector $\boldsymbol{X}_2$ can be found from second relation and, in expanded from, it is written as
\begin{equation}
 \begin{bmatrix}
         \ 0  &  n_x  &   n_y  &  0   \\[0.4em]
         \ -u u_\bot   &  u_\bot + u n_x  &  u n_y  &  0   \\[0.4em]
         \  -v u_\bot  &   v n_x   &  u_\bot + v n_y  & 0  \\[0.4em]
         \  -E u_\bot  &  E n_x   &  E n_y  &  u_\bot
        \end{bmatrix} \begin{bmatrix}
                  \ x_{1}    \\[0.4em]
                  \ x_{2}    \\[0.4em]
                  \ x_{3}    \\[0.4em]
                  \ x_4     
                 \end{bmatrix}  \ = \ u_\bot \begin{bmatrix}
                  \ x_{1}    \\[0.4em]
                  \ x_{2}    \\[0.4em]
                  \ x_{3}    \\[0.4em]
                  \ x_4    
                 \end{bmatrix}        \ + \   
                  \begin{bmatrix}
                  \ 1    \\[0.4em]
                  \ u    \\[0.4em]
                  \ v    \\[0.4em]
                  \ E 
                 \end{bmatrix} 
\end{equation}
here, each $x_{i} \in {\rm I\!R}$, where i runs from 1 to 4, and on solving all four simultaneous equations, we get
\begin{equation}
 n_x x_2 \ + \   n_y x_3   \  = \  1  \ + \  u_\bot x_1 
\end{equation}
which together with $(x_1, x_2, x_3, x_4)^t$ defines a generalized eigenvector.   
Now, if we take
\begin{equation}
 \boldsymbol{P} \ = \ \begin{bmatrix}
         \ 1   &&  x_1  &&   n_x  &&  n_y   \\[0.4em]
         \ u   &&  x_2  &&   u_\bot  &&  0   \\[0.4em]
         \ v   &&  x_3  &&   0 && u_\bot  \\[0.4em]
         \ E   &&  x_4   &&  0  &&  0
        \end{bmatrix}
\end{equation}
then 
  \begin{equation}
        \boldsymbol{P}^{-1}\boldsymbol{A}_{c}^{\boldsymbol{ZB}}\boldsymbol{P} \ = \ \ \begin{bmatrix}
       \ u_\bot  &&  1  &&  0 && 0  \\[0.4em]
       \ 0  &&  u_\bot  && 0 && 0 \\[0.4em]
       \ 0  &&  0 &&  u_\bot && 0 \\[0.4em]
       \ 0  &&  0  &&  0  && u_\bot
        \end{bmatrix} \textrm{holds.}
  \end{equation}
\\
Let $\boldsymbol{A}_{p}^{\boldsymbol{ZB}}$ denote the Jacobian matrix corresponding to pressure flux function $\boldsymbol{F}_{p}^{\boldsymbol{ZB}}$. After a little algebra $\boldsymbol{A}_{p}^{\boldsymbol{ZB}}$ comes out equal to 
\begin{equation}
                       \ (\gamma - 1)\ \begin{bmatrix}
                              \ 0  &&  0  &&  0 && 0  \\[0.4em]
                               \  \varTheta^{2} n_x  && -n_x u  && -n_x v && n_x \\[0.4em]
                              \  \varTheta^{2} n_y  &&  -n_y u &&  -n_y v  && n_y \\[0.4em]
\ \big(\varTheta^{2} - \varPhi^{2}\big) u_\bot  && \varPhi^{2}n_x - u_\bot u  && \varPhi^{2} n_y - u_\bot v  && u_\bot
        \end{bmatrix}
  \end{equation}
  where we define 
\begin{align}
\begin{split}
    \varTheta^{2} \ &= \ \dfrac{u^2 + v^2}{2}  \  \ \textrm{and} \  \  \\
    \varPhi^{2} \ &= \ \dfrac{a^2}{\gamma(\gamma -1)}
\end{split}
\end{align}
The eigenvalues of the flux Jacobian matrix $\boldsymbol{A}_{p}^{\boldsymbol{ZB}}$ are:
\begin{equation}
 \lambda_{p,1}^{\boldsymbol{ZB}} \ = \ -\sqrt{\frac{\gamma - 1}{\gamma}} a, \   \lambda_{p,2}^{\boldsymbol{ZB}} \ = \ 0, \ \lambda_{p,3}^{\boldsymbol{ZB}} \ = \ 0, \  \lambda_{p,4}^{\boldsymbol{ZB}} \ = \ \sqrt{\frac{\gamma - 1}{\gamma}} a
\end{equation}
Since all eigenvalues are real and distinct, therefore $\boldsymbol{A}_{p}^{\boldsymbol{ZB}}$ must have full set of LI eigenvectors and are given by:
\begin{equation}
 \boldsymbol{R}_{p,1}^{\boldsymbol{ZB}} \ = \ 
                  \begin{bmatrix}
                  \ 0    \\[0.4em]
                  \ n_x    \\[0.4em]
                  \ n_y    \\[0.4em]
                  \ u_\bot - \dfrac{a}{\sqrt{\gamma(\gamma - 1)}}  \
                 \end{bmatrix} \  \ , \
    \boldsymbol{R}_{p,2}^{\boldsymbol{ZB}} \ = \    
                   \begin{bmatrix}
                  \ u_\parallel    \\[0.4em]
                  \ u u_\parallel + \varTheta^{2} n_y    \\[0.4em]
                  \ v u_\parallel - \varTheta^{2} n_x   \\[0.4em]
                  \ 0  \
                 \end{bmatrix}  
 \end{equation} 
 \begin{equation}
 \boldsymbol{R}_{p,3}^{\boldsymbol{ZB}} \ = \ 
                  \begin{bmatrix}
                  \ 1    \\[0.4em]
                  \ n_x u_\bot    \\[0.4em]
                  \ n_y u_\bot   \\[0.4em]
                  \ u_\bot^2 - \varTheta^{2}  \
                 \end{bmatrix} \  \ , \
    \boldsymbol{R}_{p,4}^{\boldsymbol{ZB}} \ = \    
                   \begin{bmatrix}
                  \ 0    \\[0.4em]
                  \ n_x    \\[0.4em]
                  \ n_y   \\[0.4em]
                  \ u_\bot + \dfrac{a}{\sqrt{\gamma(\gamma - 1)}} \
                 \end{bmatrix}  
 \end{equation} 
 Now, both convection and pressure fluxes at a cell-interface are calculated by using the following.   
  \begin{equation}
\boldsymbol{F}_{c,I}^{\boldsymbol{ZB}} \ = \ \frac{1}{2} \big[\boldsymbol{F}_{c,L}^{\boldsymbol{ZB}} + \boldsymbol{F}_{c,R}^{\boldsymbol{ZB}}\big] 
          - \frac{1}{2}\left|\bar{u}_\bot \right| \Delta{\boldsymbol{U}}
\end{equation} 
where 
\begin{equation}
\Delta{\boldsymbol{U}} \ = \
 \begin{bmatrix}
         \rho_R - \rho_ L  \\[0.8em]
          \bar{\rho} \Delta{u}  \ + \ \bar{u} \Delta{\rho}  \\[0.8em]
          \bar{\rho} \Delta{v}  \ + \ \bar{v} \Delta{\rho}  \\[0.8em]
       \frac{1}{\gamma -1} \Delta{p}  \ + \ \frac{1}{2} \big(\bar{u}^2 + \bar{v}^2\big)\Delta{\rho} \ + \    \bar{\rho} (\bar{u}  \Delta{v} + \bar{v}\Delta{v})
        \end{bmatrix}
\end{equation}
and
\begin{equation}
\boldsymbol{F}_{p,I}^{\boldsymbol{ZB}}  \ = \ \frac{1}{2} \big[\boldsymbol{F}_{p,L}^{\boldsymbol{ZB}} + \boldsymbol{F}_{p,R}^{\boldsymbol{ZB}}\big] - \frac{1}{2} \sum_{i = 1}^{4} \bar{\alpha}_{p,i}^{\boldsymbol{ZB}}\left|\bar \lambda_{p,i}^{\boldsymbol{ZB}}\right| \boldsymbol{{\bar{R}}}_{p,i}^{\boldsymbol{ZB}}
\end{equation}
respectively.
Like in 1-D case, the average quantities are defined as 
\begin{align}
 \begin{split}
  \bar{\rho} \ &= \ \sqrt{\rho_L \rho_R}  \ , \
  \bar{u}   \ = \   \frac{\sqrt{\rho_L} u_L  \ + \ \sqrt{\rho_R} u_R }{\sqrt{\rho_L} \ + \ \sqrt{\rho_R}} \\
  \bar{v}   \ &= \   \frac{\sqrt{\rho_L} v_L  \ + \ \sqrt{\rho_R} v_R }{\sqrt{\rho_L} \ + \ \sqrt{\rho_R}} \ , \
  \bar{a}^2   \ = \   \frac{\sqrt{\rho_L} a_L^2  \ + \ \sqrt{\rho_R} a_R^2 }{\sqrt{\rho_L} \ + \ \sqrt{\rho_R}} \\
  \bar{u}_\bot \ &= \  \bar{u} n_x + \bar{v} n_y \ \ \textrm{and} \ \  \bar{u}_\parallel \ = \  -\bar{u} n_y + \bar{v} n_x  
 \end{split}
\end{align}
where $u_\parallel$ denotes velocity component parallel to cell-interface and is given by 
\begin{equation}
 u_\parallel \ = \ -u n_y \ + \ v n_x
\end{equation}

Similarly, wave strengths $\bar{\alpha}_{p,i}^{\boldsymbol{ZB}}$, where $i$ runs from $1 \ \textrm{to} \ 4$, are given as 
\begin{align}
 \begin{split}
  \bar{\alpha}_{p,1}^{\boldsymbol{ZB}} \ &= \  \dfrac{\bar{\rho} \Delta{u_\bot}}{2} - \sqrt{\dfrac{\gamma}{\gamma -1}} \dfrac{\Delta{p}}{2\bar{a}} \\
  \bar{\alpha}_{p,2}^{\boldsymbol{ZB}} \ &= \ \dfrac{\bar{u}_\parallel \Delta{\rho} + \bar{\rho}\Delta{u_\parallel}}{\bar{\varTheta}^{2}-\bar{u}_\bot^2} \\
  \bar{\alpha}_{p,3}^{\boldsymbol{ZB}} \ &= \ \Delta{\rho} \ - \  \dfrac{\bar{u}_\parallel^{2} \Delta{\rho} + \bar{\rho}\bar{u}_\parallel\Delta{u_\parallel}}{\bar{\varTheta}^{2}-\bar{u}_\bot^2} \\
  \bar{\alpha}_{p,4}^{\boldsymbol{ZB}} \ &= \  \dfrac{\bar{\rho} \Delta{u_\bot}}{2} + \sqrt{\dfrac{\gamma}{\gamma -1}} \dfrac{\Delta{p}}{2\bar{a}}
 \end{split}
\end{align}
where 
\begin{align}
 \begin{split}
  \bar{\varTheta}^{2} \ &= \ \dfrac{\bar{u}^2 + \bar{v}^2}{2}  \\
  \Delta{\rho} \ &= \ \rho_{R} - \rho_{L}  \\
  \Delta{u_\bot} \ &= \  u_{\bot \ R} - u_{\bot \ L}  \\
\Delta{u_\parallel} \ &= \  u_{\parallel \ R} - u_{\parallel \ L}  \\
\Delta{p}  \ &= \  p_R - p_L 
 \end{split}
\end{align}  
\subsection{Numerical examples}
In this subsection, the ZBS-FDS scheme is tested on various well-established benchmark test problems.  Special attention is given to  problems with complex interactions of strong shocks which leads to various shock instabilities.  Many well known upwind schemes are known to produce unphysical features in such cases \cite{Quirk}.
\subsubsection{Oblique shock reflection}
\begin{figure}[!ht]
 \begin{center}
 \begin{minipage}{0.48\linewidth}
 \centering
(a)
 {\includegraphics[trim=0.0cm 0.0cm 0cm 0.0cm, clip, width=\textwidth]{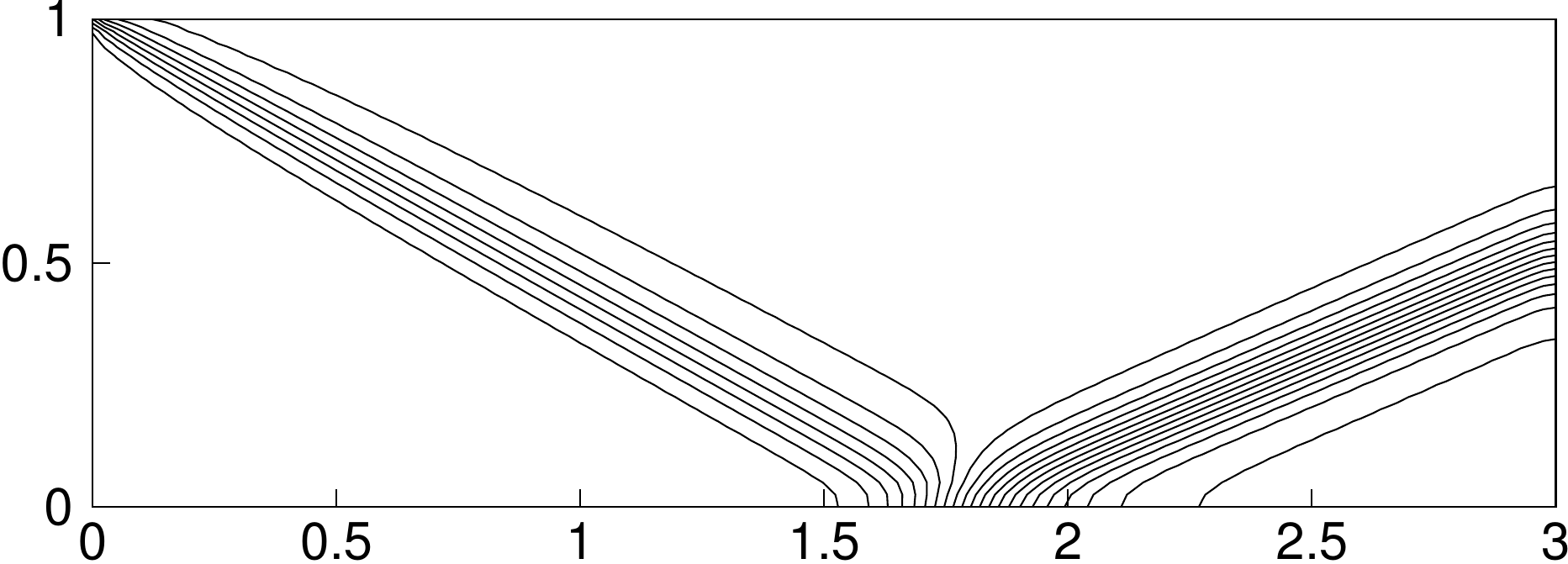}}
 \end{minipage}
 \begin{minipage}{0.48\linewidth}
  \centering
(a)
  {\includegraphics[trim=0.0cm 0.0cm 0cm 0.0cm, clip, width=\textwidth]{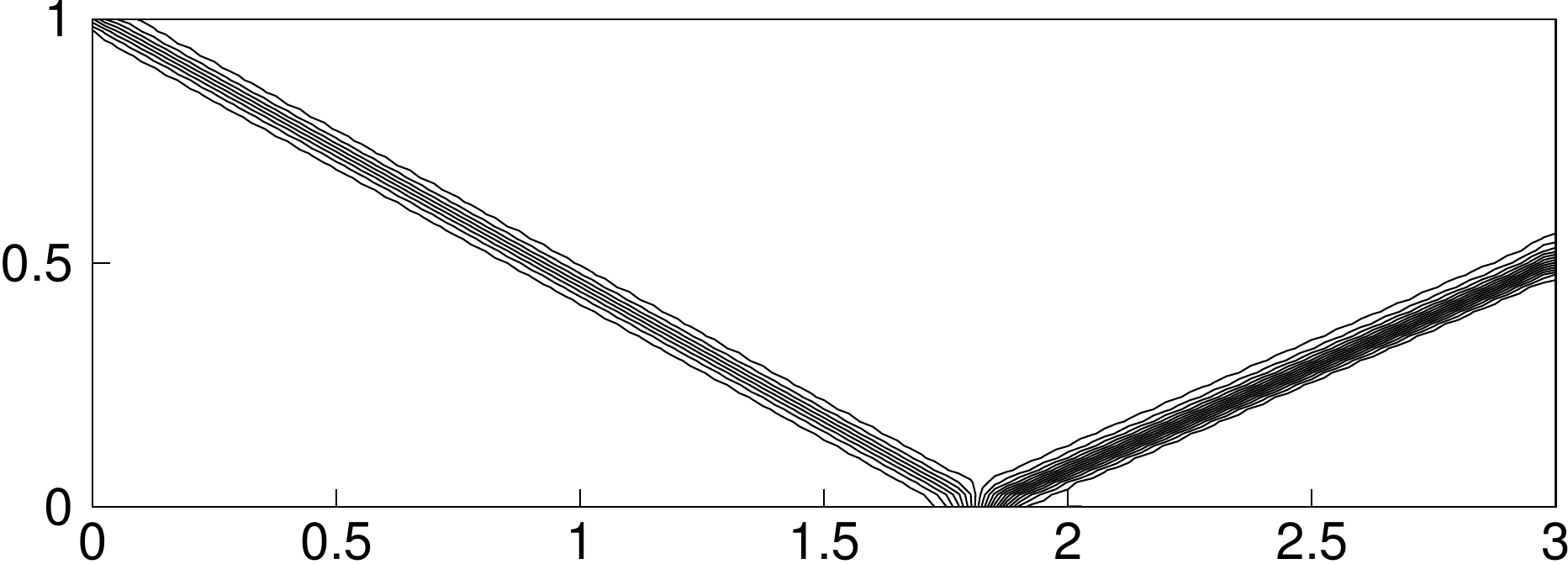}}
\end{minipage}
 \begin{minipage}{0.48\linewidth}
 \centering
(b)
 {\includegraphics[trim=0.0cm 0.0cm 0cm 0.0cm, clip, width=\textwidth]{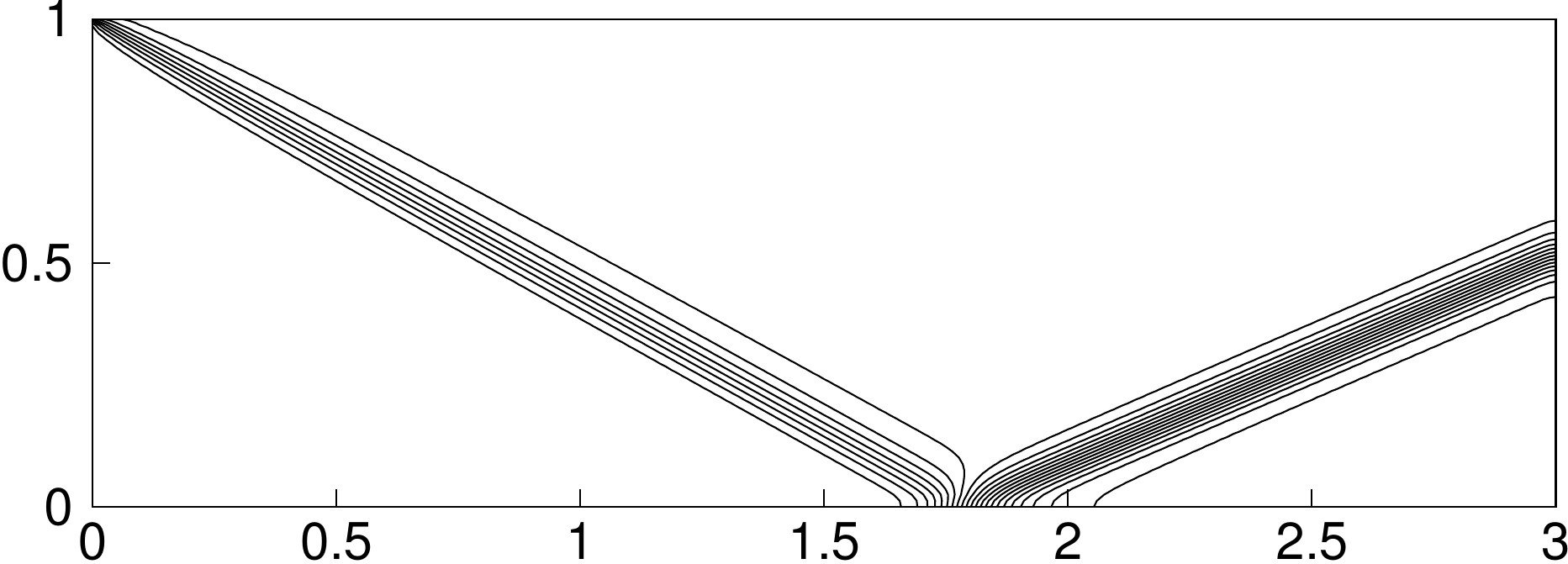}}
 \end{minipage}
 \begin{minipage}{0.48\linewidth}
  \centering
(b)
  {\includegraphics[trim=0.0cm 0.0cm 0cm 0.0cm, clip, width=\textwidth]{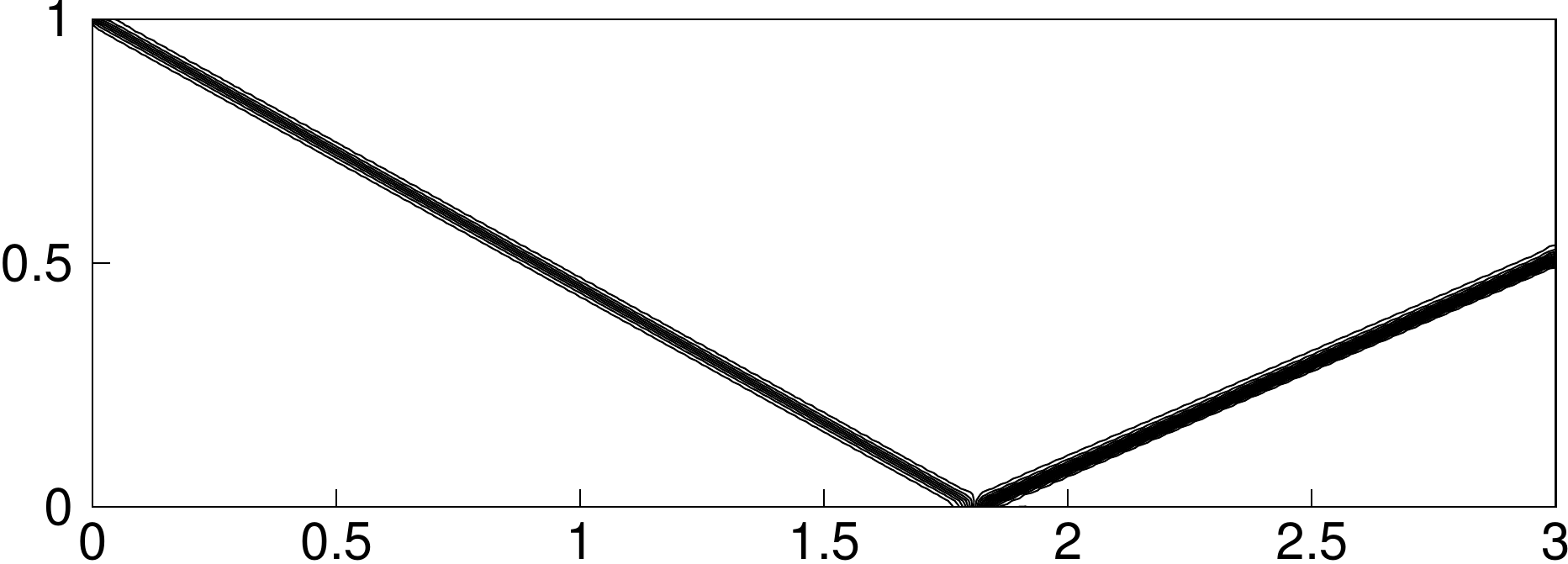}}
\end{minipage}
  \caption{First order results of ZBS-FDS scheme are presented on left, where second order results are given in right side for shock reflection problem; pressure contours (0.7: 0.1: 2.9) on the grids: (a) $120\times40$ and (b) $240\times80$}
 \label{shock_ref_ZBS_FDS}
 \end{center}
 \end{figure}

In this test case \cite{shock_ref_1982}, an oblique shock wave is introduced at the top left corner by means of initial conditions and post shock boundary conditions, at the left and top side of the domain, respectively. The computational domain considered for this test case is $[0,3] \times [0,1]$. The initial conditions for this test problem are as given below.  
\begin{equation*}
\left( \rho, u, v, p \right)_{0,y,t}=\left( 1.0, 2.9, 0, 1/1.4 \right) 
\end{equation*}
\begin{equation*}
\left( \rho, u, v, p \right)_{x,1,t}=\left( 1.69997, 2.61934, -0.50633, 1.52819 \right)
\end{equation*}
The incident shock angle measured from the top side of the domain is $29^0$ and the free stream Mach number $M=2.9$.  Wall boundary conditions are prescribed at the bottom boundary and supersonic outflow boundary conditions are used at the right side of the computational domain.  Both first order and second order results on various grids are presented in Figure \ref{shock_ref_ZBS_FDS}. 
\subsubsection{Supersonic flow across a compression ramp in a wind tunnel}  
\begin{figure}[!htbp]
 \begin{center}
 \begin{minipage}{0.48\linewidth}
 \centering
(a)
 {\includegraphics[trim=0.0cm 0.0cm 0cm 0.0cm, clip, width=\textwidth]{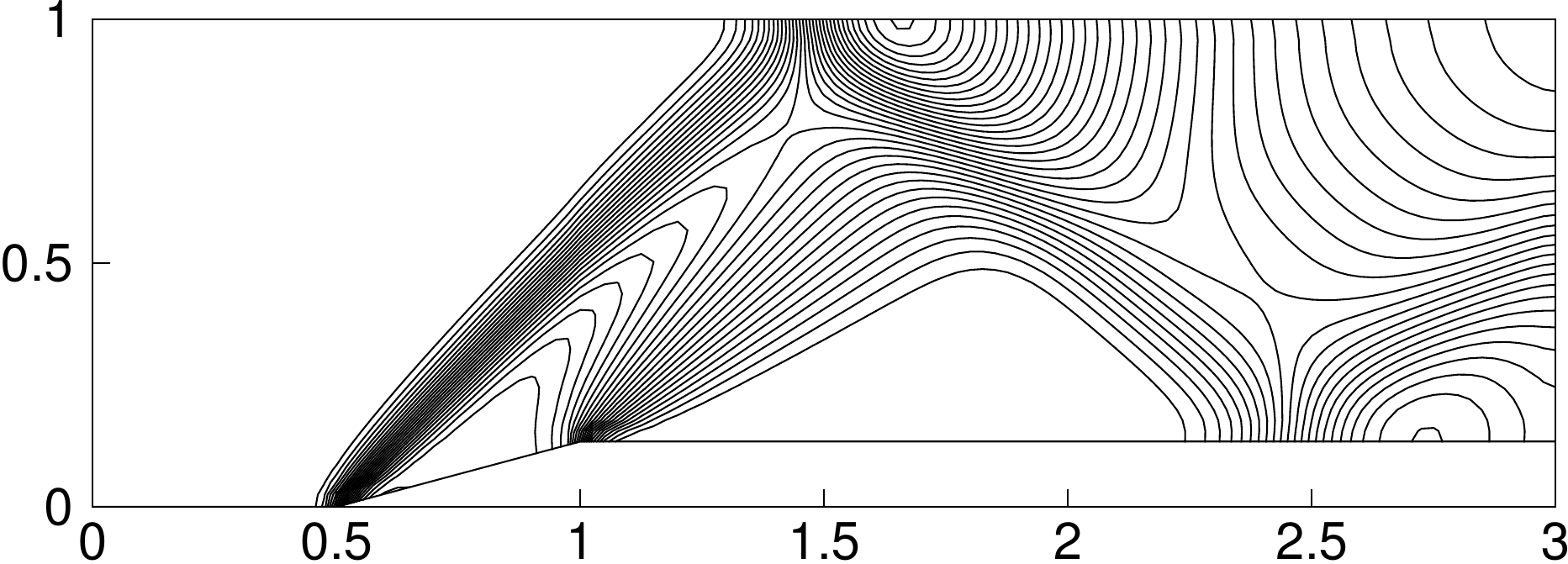}}
 \end{minipage}
 \begin{minipage}{0.48\linewidth}
  \centering
(a)
  {\includegraphics[trim=0.0cm 0.0cm 0cm 0.0cm, clip, width=\textwidth]{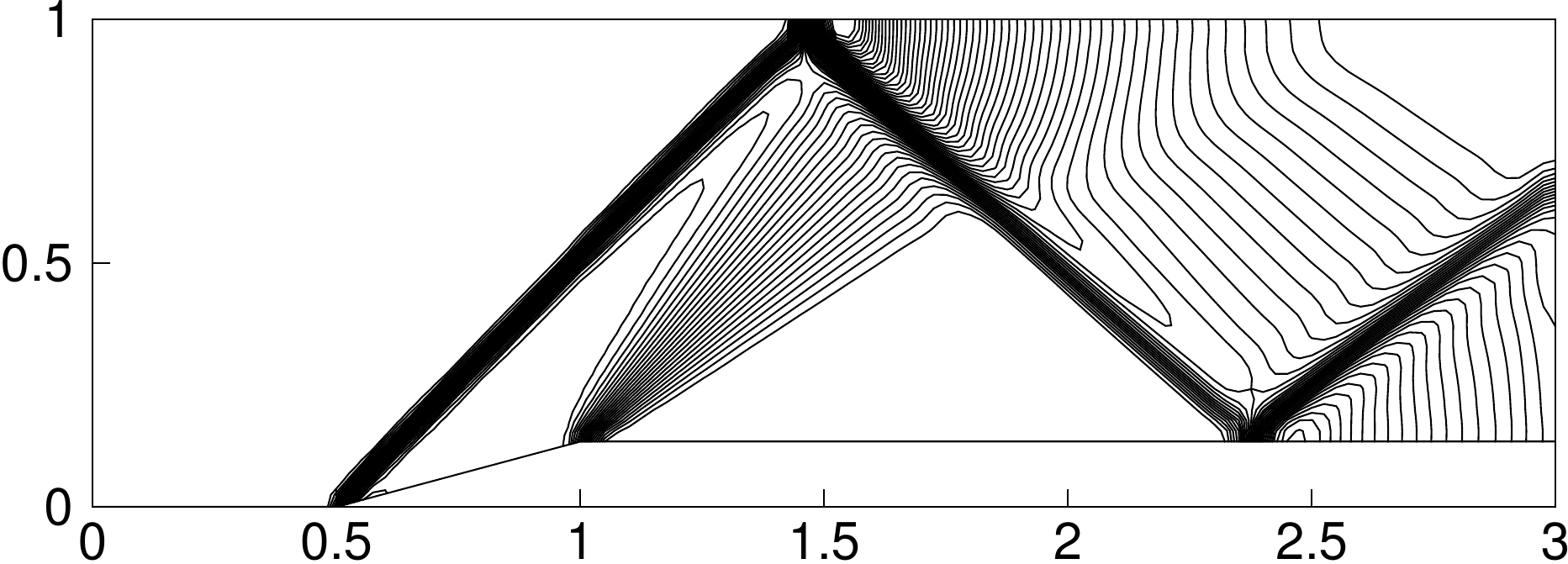}}
\end{minipage}
 \begin{minipage}{0.48\linewidth}
 \centering
(b)
 {\includegraphics[trim=0.0cm 0.0cm 0cm 0.0cm, clip, width=\textwidth]{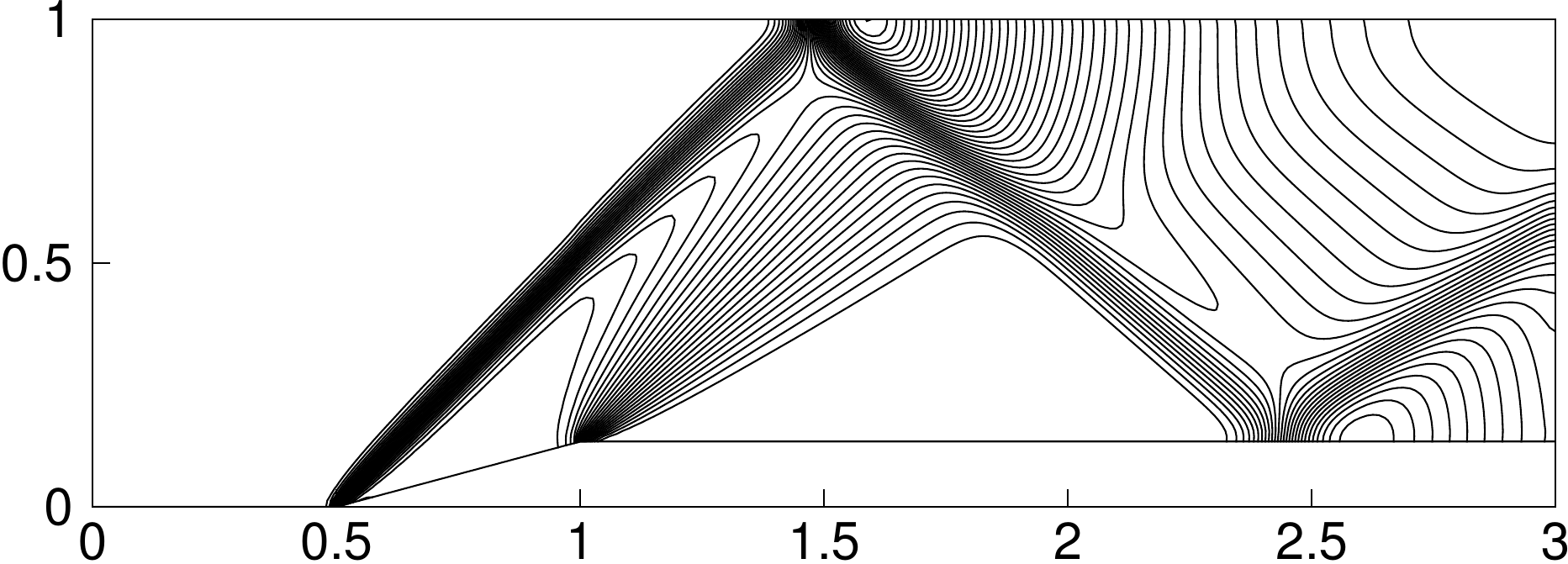}}
 \end{minipage}
 \begin{minipage}{0.48\linewidth}
  \centering
(b)
  {\includegraphics[trim=0.0cm 0.0cm 0cm 0.0cm, clip, width=\textwidth]{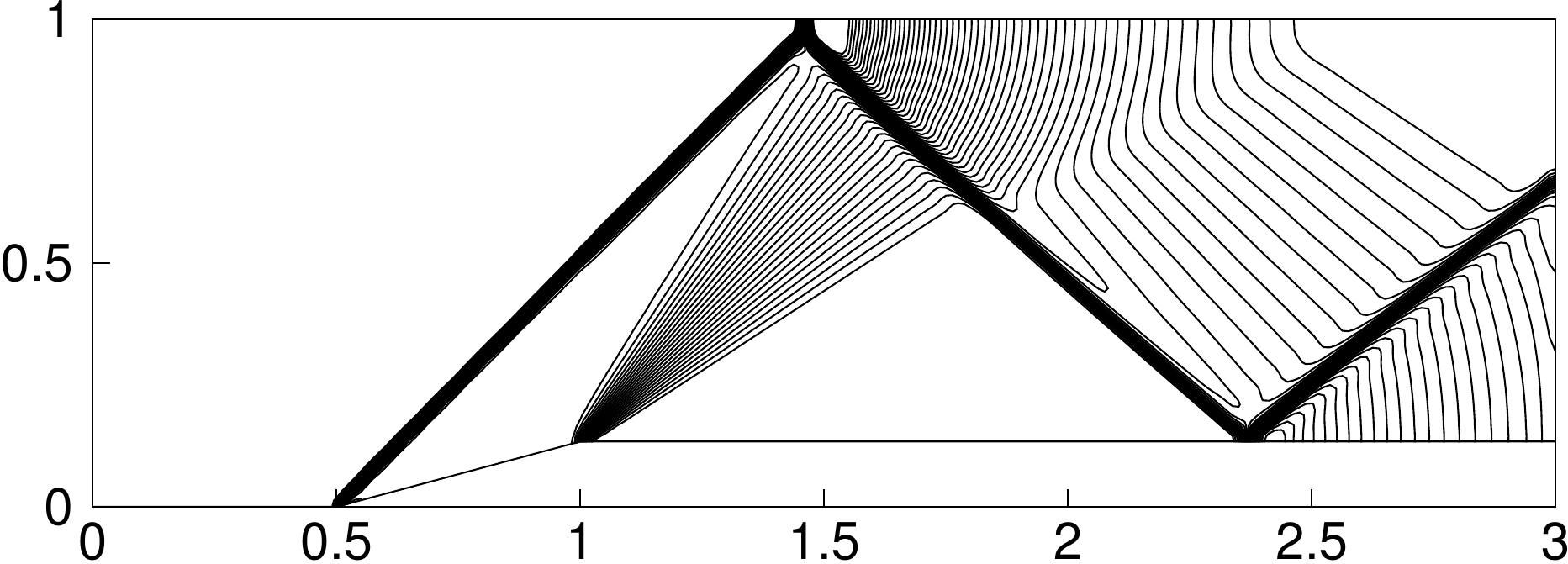}}
\end{minipage}
  \caption{First order results of ZBS-FDS scheme are given on left and second order results are presented on right for ramp reflection problem; pressure contours (1.1: 0.05: 3.8) on the grids: (a)  $120\times40$ and (b) $240\times80$}
 \label{wedge_ZBS-FDS}
 \end{center}
 \end{figure}

 The computational domain of $[0, 3] \times [0, 1]$ is considered for this test problem. Other geometrical features of the problem include a $15^0$ ramp at the lower part of the computational domain. In this two-dimensional steady test case~\cite{Levy_1993}, supersonic flow of a Mach number $M=2$ encounters a fifteen degree ramp to form an oblique shock wave. This shock wave reflects from the upper wall and interacts with the expansion wave generated at the tip of the ramp corner. The so weakened expansion wave again reflects from the top wall and further interacts with the second reflected shock from the ramp surface.  Initial conditions are prescribed at the left boundary, wall boundary conditions are used at the top and on the ramp, and supersonic outflow boundary conditions are imposed at the exit boundary. Both first order and second order results are presented in Figure \ref{wedge_ZBS-FDS}.
\subsubsection{Reflection of a plane shock from wedge}
\begin{figure}[!htbp]
 \begin{center}
 \begin{minipage}{0.49\linewidth}
 \centering
  {\includegraphics[trim=0.0cm 0.0cm 0cm 0.0cm, clip, width=\textwidth]{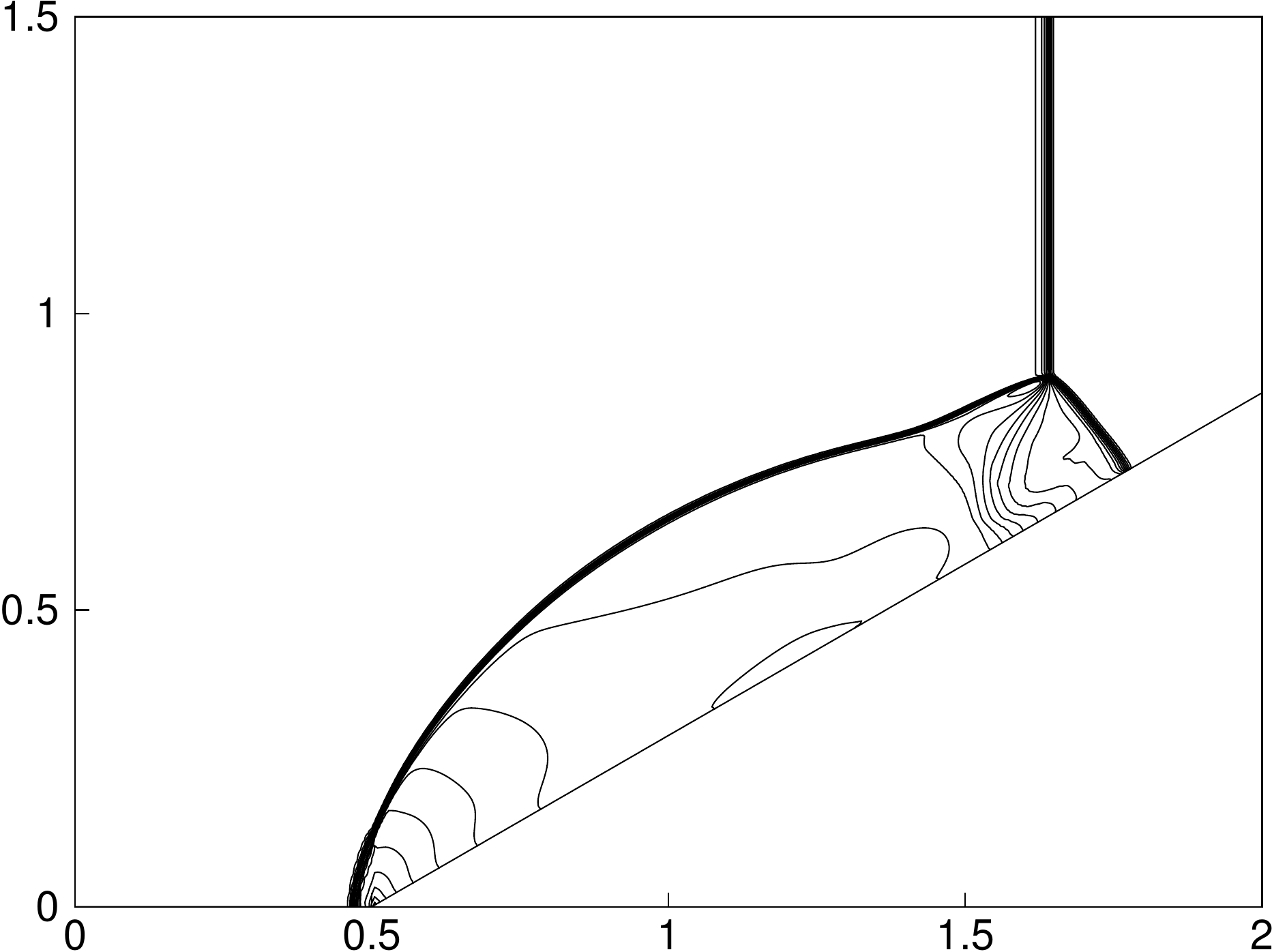}}
 \end{minipage}
 \begin{minipage}{0.49\linewidth}
  \centering
   {\includegraphics[trim=0.0cm 0.0cm 0cm 0.0cm, clip, width=\textwidth]{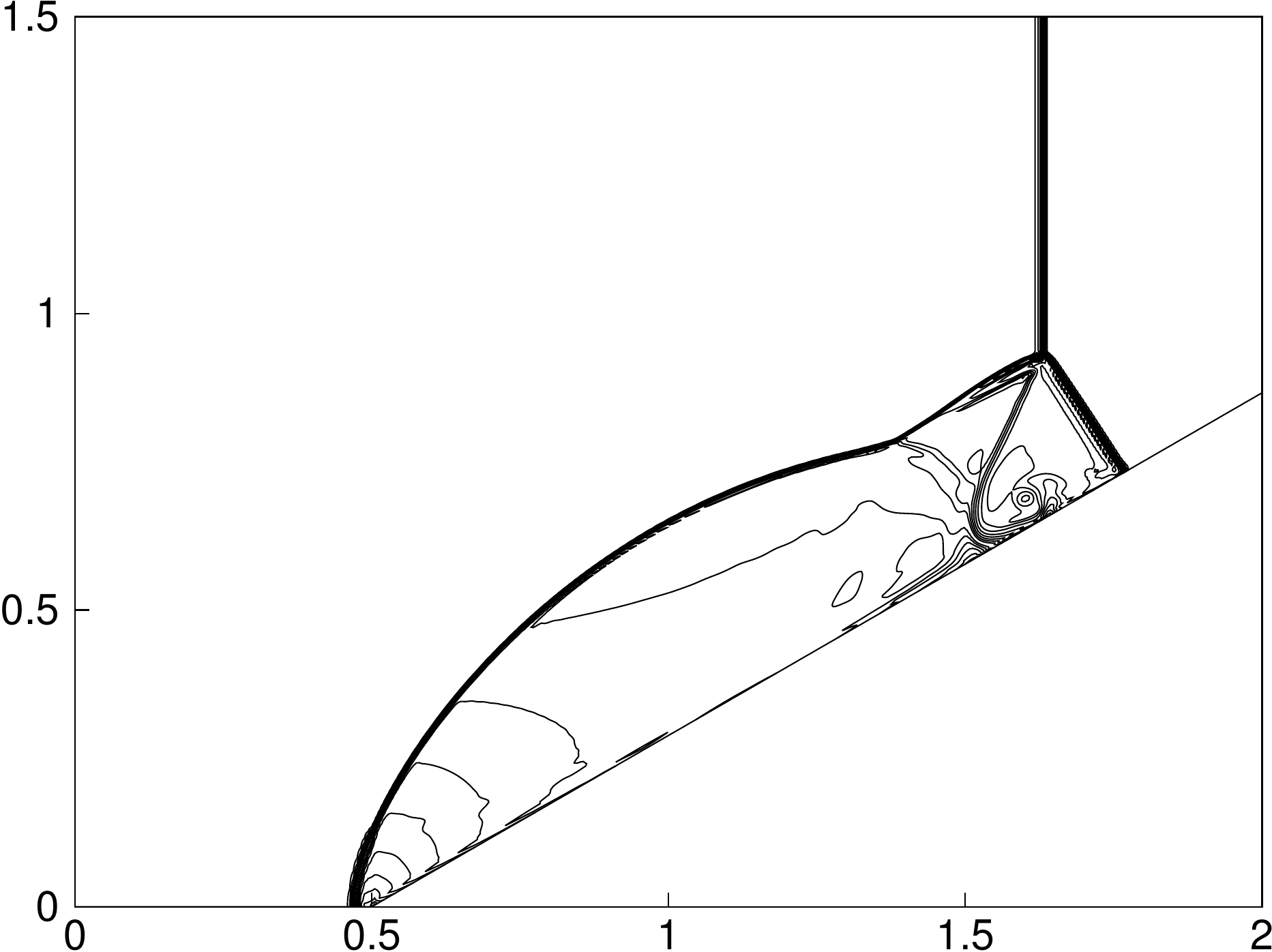}}
\end{minipage}
  \caption{First (left) and Second (right) order results of reflection of a plane shock from a wedge problem with ZBS-FDS scheme on $400\times400$ grid points}
 \label{wr}
 \end{center}
 \end{figure}

In this is a two-dimensional problem in which the reflection of a plane shock wave from a wedge lies in the double-Mach reflection regime, some Riemann solvers are  known to generate kinked Mach stems \cite{Quirk}. In this test case, the kinked Mach stem occurs when a strong normal shock wave moving with Mach $5.5$ encounters the $30^0$ ramp to form Mach reflection and represents a typical shock instability phenomenon.  Three shocks meet to form a triple point and the computational domain considered for this problem is $[0, ~2.0] \times [0, ~1.5]$ with initial shock location at $x_0=0.25$.  All computational results are obtained at time $t = 0.25$. The computational domain to the right of the shock is initialized with a rest fluid of density $1.4$ and pressure $1$.  To the left of the shock, values obtained from the moving shock relations for Mach $5.5$ are used to initialize the domain. Figure \ref{wr} shows the density contours computed with the present scheme and no kink is observed.  
\subsubsection{Hypersonic flow past a half-cylinder}
\begin{figure}[!htbp]
 \begin{center}
 \begin{minipage}{0.18\linewidth}
 \centering
 (a)
 {\includegraphics[trim=0.0cm 0.0cm 0cm 0.0cm, clip, width=\textwidth]{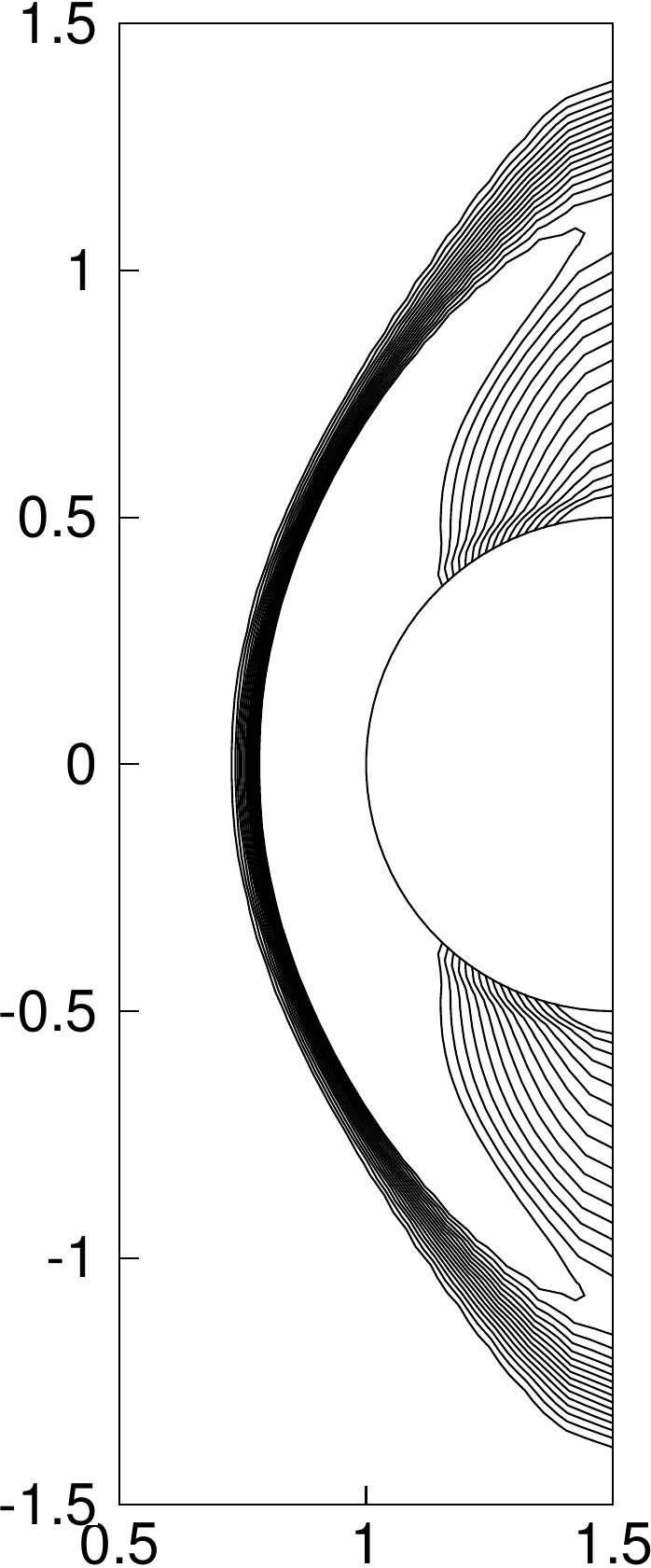}}
 \end{minipage}
 \begin{minipage}{0.18\linewidth}
  \centering
  (b)
  {\includegraphics[trim=0.0cm 0.0cm 0cm 0.0cm, clip, width=\textwidth]{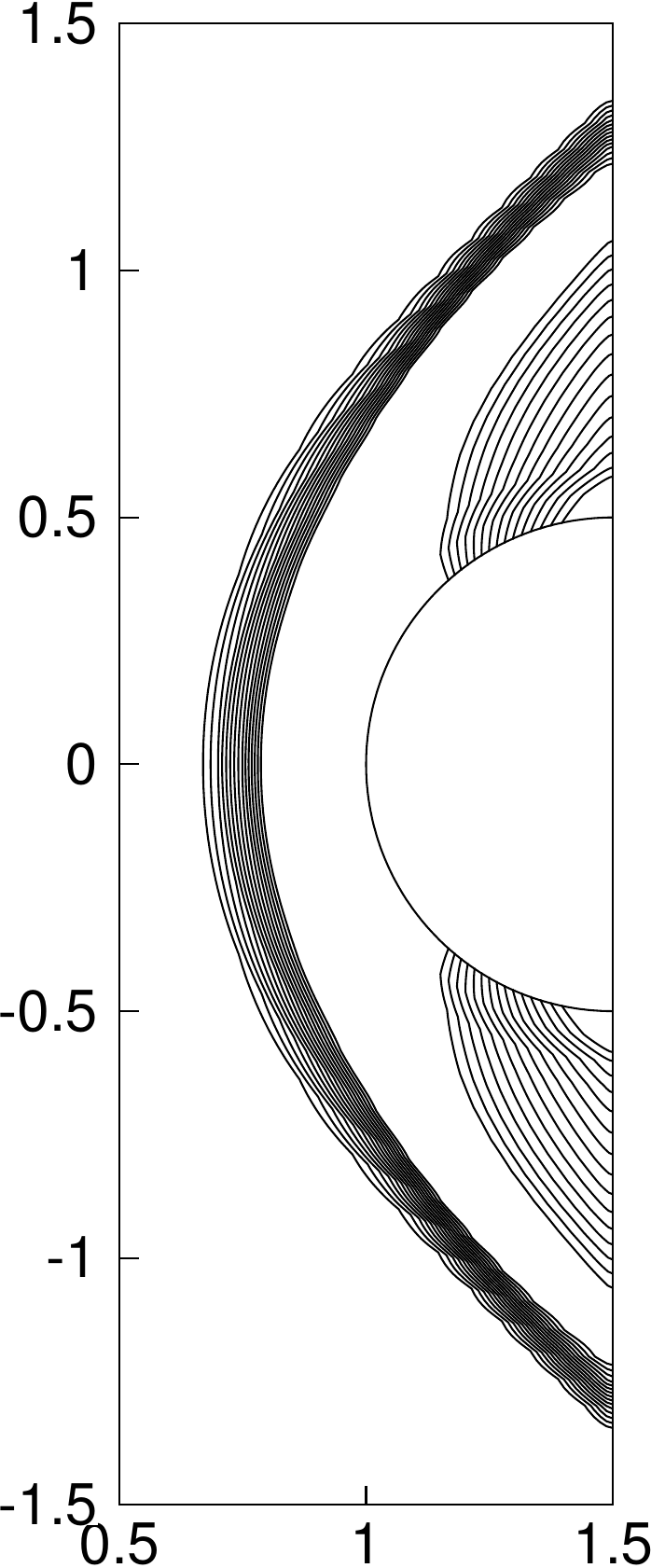}}
\end{minipage}
 \begin{minipage}{0.24\linewidth}
 \centering
Roe
 {\includegraphics[trim=0.0cm 0.0cm 0cm 0.0cm, clip, width=\textwidth]{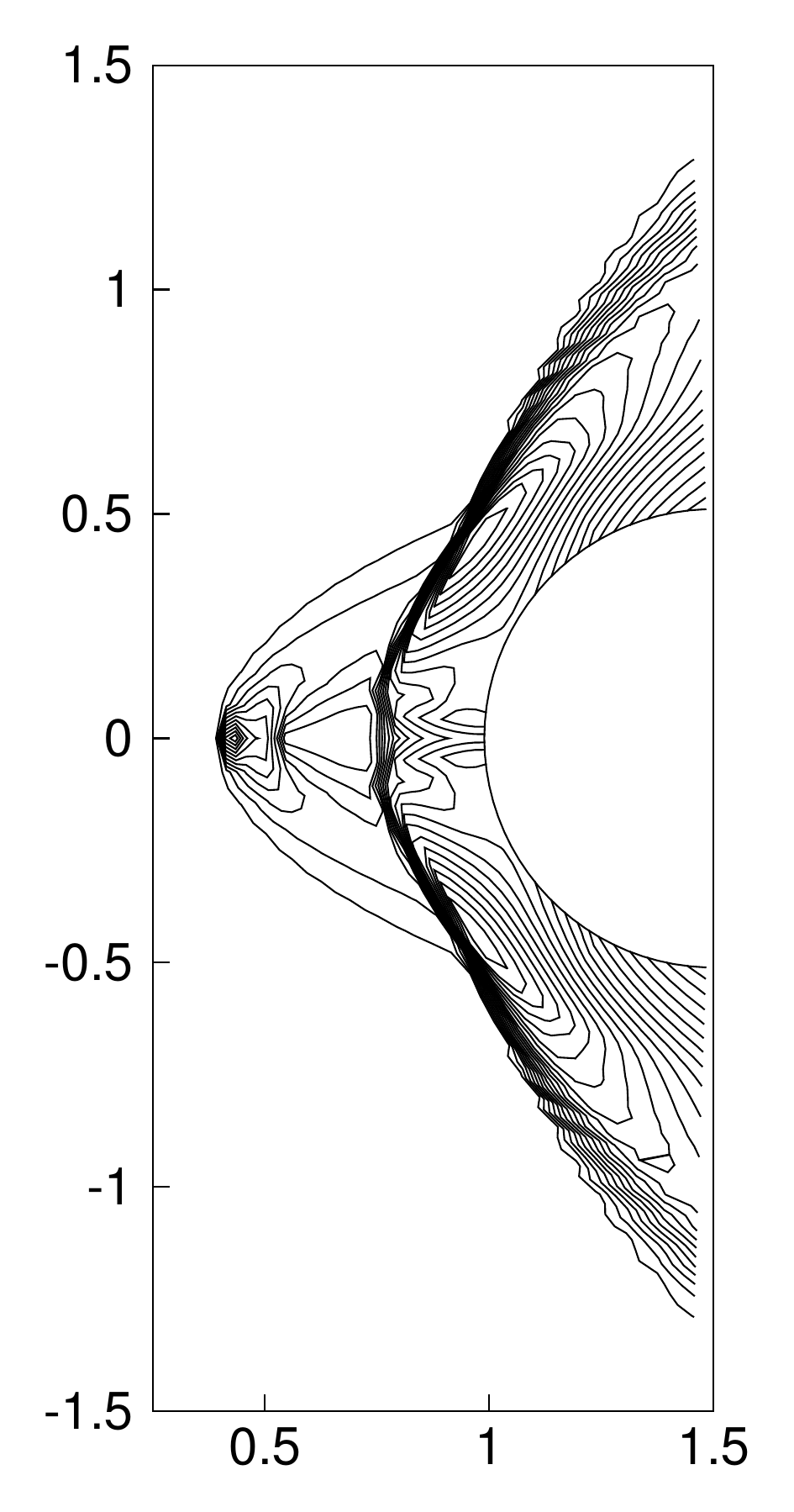}}
 \end{minipage}
 \begin{minipage}{0.18\linewidth}
  \centering
(c)
  {\includegraphics[trim=0.0cm 0.0cm 0cm 0.0cm, clip, width=\textwidth]{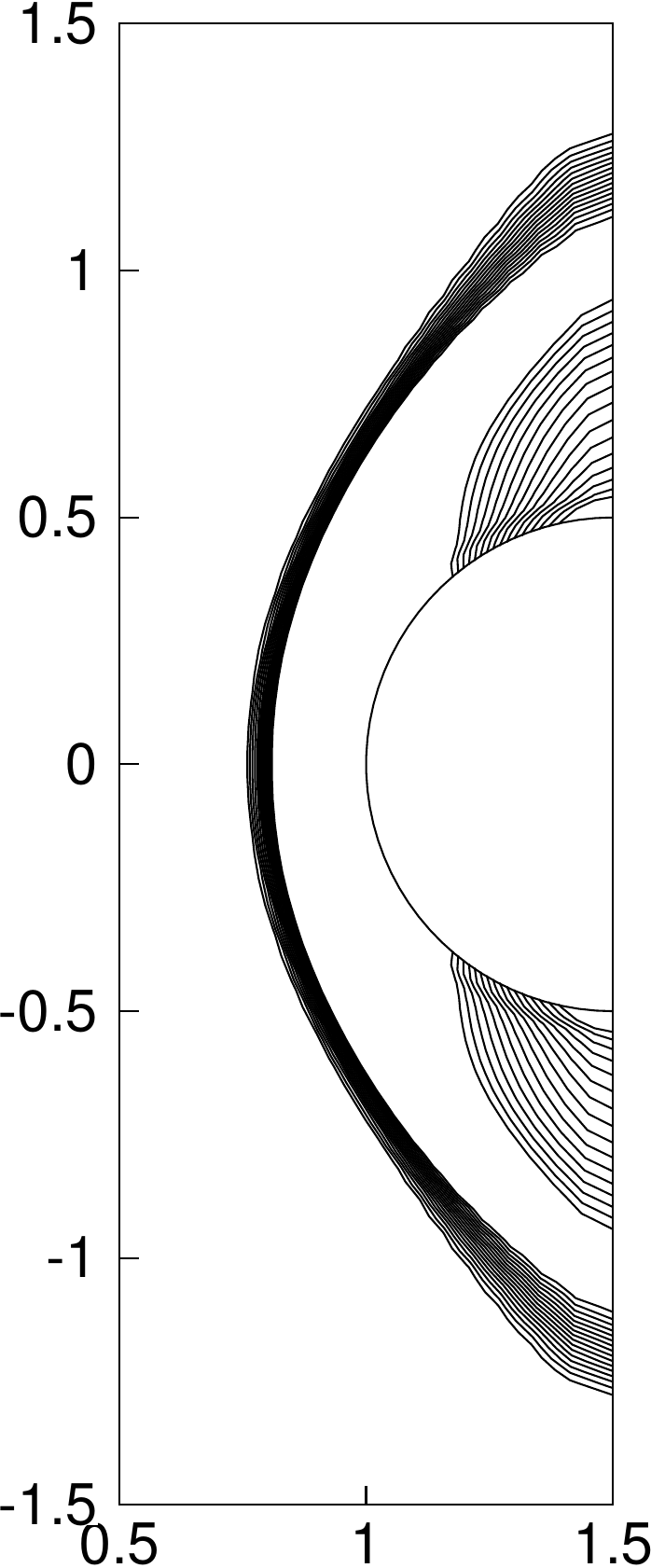}}
\end{minipage}
 \begin{minipage}{0.18\linewidth}
 \centering
(d)
 {\includegraphics[trim=0.0cm 0.0cm 0cm 0.0cm, clip, width=\textwidth]{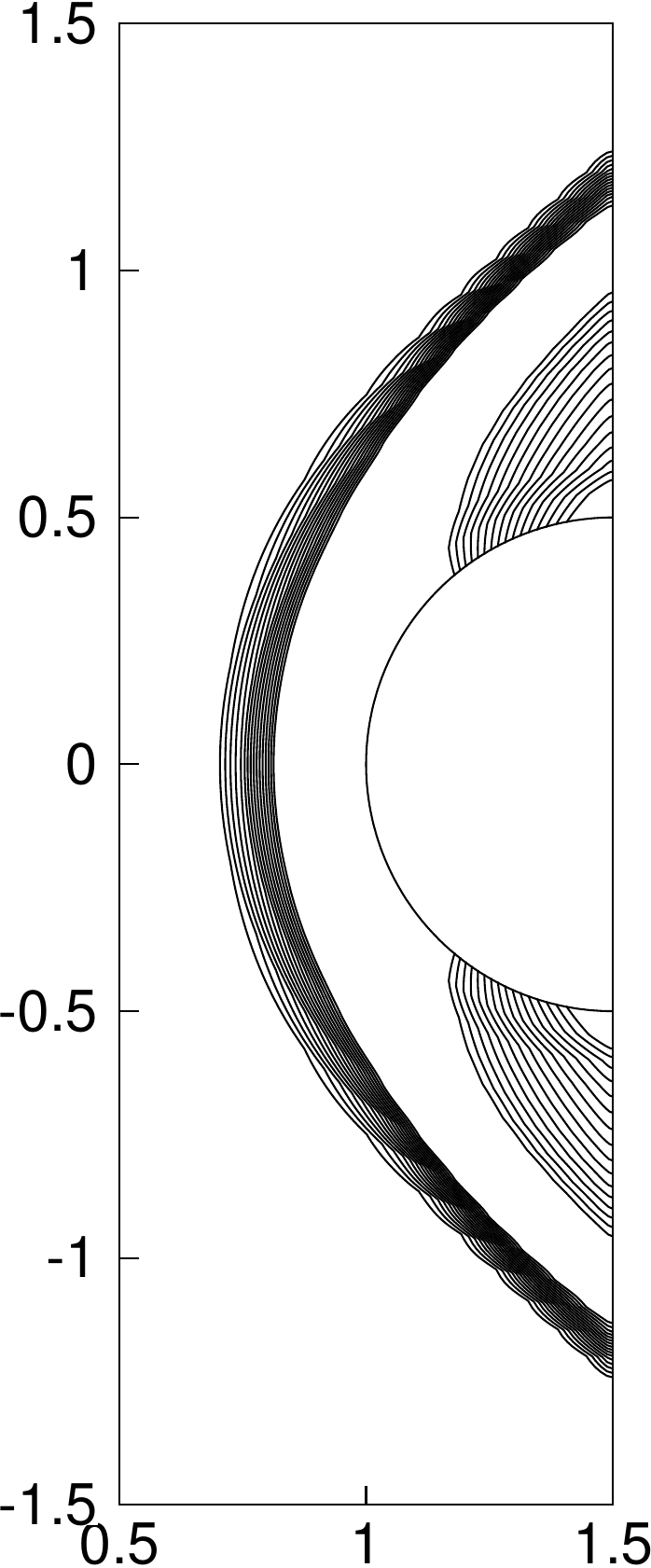}}
 \end{minipage}
 \caption{First order results of ZBS-FDS scheme for half cylinder problem; density contours (2.0: 0.2: 5.0): (a) Mach $6$ on $45\times45$ grid, (b) Mach $6$ on $20\times320$ grid, (c) Mach $20$ on $45\times45$ and (d) Mach $20$ on $20\times320$ grid}
 \label{half_cylinder_ZBS}
 \end{center}
 \end{figure}

The hypersonic flow around a half-cylinder is also a well-known test case to examine the capability of numerical methods in resolving complex flow features accurately without giving shock instabilities.  In this case, the shock instability is known as 'carbuncle shock' which was initially reported by Peery and Imlay \cite{Peery_carbuncle}.  This test is computed for Mach $6$ and Mach $20$ flows on fine and coarse grids in circumferential directions and results are given in Figure \ref{half_cylinder_ZBS}.  Many Riemann solvers generate carbuncle shocks \cite{Quirk,Kim_et_al} in the their numerical solutions.  As an example we presented results of Roe scheme, where carbuncle phenomena
can be seen very clearly.  The present method did not exhibit any such phenomena.  
 
\section{Summary}
In this study, we attempted to develop flux difference split upwind schemes for  convection-pressure splitting frameworks, based on Jordan canonical forms to avoid defective matrices.  FDS solver for Liou and Steffen type splitting is not attractive as pressure subsystem doesn't have any contribution of acoustic signals.  Newly constructed ZBS-FDS and TVS-FDS schemes are tested on various benchmark problems and don't require any entropy fix for sonic point and strong expansion problems.  ZBS-FDS scheme and TVS-FDS scheme perform in similar ways as is evident from several 1-D test cases.  We further extend ZBS-FDS scheme to 2-D Euler system, specifically to those test problems for which several Riemann solvers generate shock instabilities \cite{Quirk,Mandal,Huang_Wu_Yan,Shen_Yan}. The performance of ZBS-FDS scheme for these test cases is impressive and deserves further research.  For example, accuracy of ZBS-FDS scheme in multi-dimensions can be further improved using some good diffusion regulator such as \cite{DR} and other possible direction is to pursue genuinely multi-dimensional modeling, as the eigenvalues of the convection and pressure parts of the fluxes neatly reflect uni-directional and multi-directional information propagation respectively. 
\section{Acknowledgments} 
The authors thank Prof. Michael Junk, Fachbereich Mathematik und Statistik, Universit\"at Konstanz, Germany for very useful discussions.  The authors also thank Indian Institute of Science for supporting this research.   

\end{document}